\documentclass[12pt]{amsbook}
\usepackage[all,cmtip]{xy}

\usepackage{a4wide}
\usepackage{pxfonts}
\usepackage{wasysym}
\usepackage{amssymb,amsxtra,mathrsfs}
\usepackage{eucal}
\usepackage{setspace,upgreek}
\usepackage{mathtools}
\usepackage{url}
\usepackage[hidelinks]{hyperref}

\newcommand{\mb}[1]{\mathbf{#1}}

\newtheorem{theorem}{Theorem}[chapter]
\newtheorem{proposition}[theorem]{Proposition}
\newtheorem{lemma}[theorem]{Lemma}
\newtheorem{corollary}[theorem]{Corollary}

\newtheorem{maintheorem}[theorem]{Main Theorem}

\theoremstyle{definition}
\newtheorem{definition}
[theorem]
{Definition}

\newtheorem{preliminaries}
[theorem]
{Preliminaries}
\newtheorem{notations}
[theorem]
{Notations}
\newtheorem{notationsE}
[theorem]
{Notations with Comments}
\newtheorem{remark}
[theorem]
{Remark}

\newcommand{\on}{
\overset{\|\cdot\|}
{\subseteq}
}
\newcommand{\up}{\upharpoonright}
\newcommand{\w}[1]{\widetilde{#1}}
\newcommand{\ze}{\mathbf{0}} 
\newcommand{\un}{\mathbf{1}} 
\newcommand{\ep}{\varepsilon} 

\newcommand{\Df}{\mathsf{D}}

\newcommand{\F}[2]{\mathfrak{F}_{#1}^{#2}}

\newcommand{\RM}{\mathcal{R}} 
\newcommand{\DM}{\mathcal{D}} 
\newcommand{\p}{\mathcal{P}}

\newcommand{\cc}{\mathcal{C}}
\newcommand{\N}{\mathbb{N}}
\newcommand{\K}{\mathbb{K}}
\newcommand{\R}{\mathbb {R}}
\newcommand{\C}{\mathbb {C}}

\newcommand{\A}{\mathcal{A}}

\newcommand{\n}{\mathcal{N}}

\newcommand{\D}{\mathbf{D}}
\newcommand{\Lf}[2]{\mathfrak{L}_{#1}^{#2}}
\newcommand{\I}[2]{\mathbf{I}_{#1}^{#2}}
\newcommand{\B}{\mathcal{B}}
\newcommand{\bb}{\mathbf{T}\mathbf{M}}
\newcommand{\ov}{\overline}

\newcommand{\ra}{\right\rangle}
\newcommand{\lr}[2]{\left\langle #1,#2\ra}

\begin{document}

\title{Integral Equalities for Functions of Unbounded Spectral Operators in Banach Spaces}

\author{Benedetto Silvestri}

\email{bndsilvestri@proton.me}

\address{440 W Kelly Ave Apt 6 Jackson, WY 83001-8544 USA (via express couriers signature required)}

\urladdr{\url{https://arxiv.org/a/silvestri_b_1}}

\urladdr{https://orcid.org/0000-0002-8042-4932}

\thanks{Research
supported 
by
the
Engineering and Physical 
Sciences Research Council
(EPSRC)
}

\keywords{Unbounded spectral operators in Banach spaces, functional calculus, integration of locally convex space-valued maps.}

\subjclass[2010]
{47B40,
47A60,
46G10}

\maketitle

\tableofcontents

\chapter*{Introduction}

 
The work is dedicated to 
investigating
a
limiting procedure
for
extending ``local'' 
integral
operator
equalities
to the ``global '' ones
in the sense 
explained below,
and
to 
applying it
to 
obtaining 
generalizations
of the
Newton-Leibnitz
formula
for operator-valued
functions for a wide
class of unbounded 
operators.
 
The integral
equalities
considered in the work 
have the following form
\begin{equation}
\label{13561401}
g(R_{F})
\int\,
f_{x}(R_{F})
\,d\,
\mu(x)
=
h(R_{F}).
\end{equation}
They
involve
functions
of the kind
\begin{equation*}
X\ni x\mapsto f_{x}(R_{F})\in B(F),
\end{equation*}
where 
$X$
is a
general
locally compact
space,
$F$ is  
a suitable 
Banach subspace 
of a fixed complex Banach space $G$,
for example
$F=G$.
The integrals are 
with respect to a 
general complex
Radon measure on $X$
and
with respect to the
$\sigma(B(F),\n_{F})-$
topology
\footnote{
Here
$\n_{F}$
is a suitable 
subset of
$B(F)^{*}$,
the topological dual of $B(F)$,
associated
with
the resolution of the identity
of 
$R_{F}$.
}
on $B(F)$. 
$R_{F}$ is a possibly 
unbounded scalar type
spectral operator in 
$F$ 
such that
$
\sigma(R_{F})
\subseteq\sigma(R_{G})
$,
and
for all $x\in X$,
$f_{x}$ 
and $g,h$
are complex-valued
Borelian maps on the
spectrum 
$\sigma(R_{G})$
of $R_{G}$.
 
If
$F\neq G$
we call 
the integral equalities 
\eqref{13561401}
``local'',
while
if
$F=G$
we call them 
``global''.
 
Let $G$ be a complex Banach space 
and
$B(G)$
the Banach algebra of all bounded 
linear operators on $G$.
\emph{Scalar type spectral operators}
in $G$
were defined in \cite{ds}
Definition $18.2.12.$
\footnote{
For the special
case of bounded
spectral operators
on $G$
see \cite{dow}.
}
(see Section \ref{09061650}),
and were created for providing
a general Banach space with a
class of
unbounded linear
operators for which it is possible 
to establish a Borel functional
calculus 
similar to the well-known 
one for unbounded self-adjoint operators
in a Hilbert space.
 
We
start with the following
useful formula 
\footnote{
An important application
of 
this formula
is made 
in proving
the well-known
Stone theorem
for strongly continuous
semigroups of unitary
operators in Hilbert space,
see Theorem $12.6.1.$
of \cite{ds}.
In \cite{dav}
it
has been used for
showing
the equivalence
of
uniform
convergence
in 
strong
operator 
topology
of
a one-parameter
semigroup
depending on a
parameter
and
the
convergence
in 
strong
operator 
topology
of the 
resolvents
of the
corresponding
generators,
Theorem
$3.17.$
}
for the resolvent of $T$
\begin{equation}
\label{03101514}
(T-\lambda\un)^{-1}
=
i
\int_{-\infty}^{0}
e^{-it\lambda}
e^{ it T}
\,d\,t.
\end{equation}
Here 
$\lambda\in\C$
is
such that $Im(\lambda)>0$
and the integral is 
with 
respect to the Lebesgue
measure
and with respect
to the
strong 
operator topology on $B(G)$
\footnote{
Notice that
if
$
\zeta
\doteqdot
-i\lambda
$
and
$
Q
\doteqdot
iT
$,
then the equality
\eqref{03101514}
turns into
$$
(Q+\zeta\un)^{-1}
=
\int_{0}^{\infty}
e^{-t\zeta}
e^{-Q t}
\,d\,t,
$$
which 
is referred 
in 
$IX.1.3.$
of
\cite{kato}
as
the fact 
that
the resolvent
of $Q$
is the 
\emph{Laplace}
trasform
of the semigroup
$e^{-Q t}$.
Applications
of this formula
to perturbation theory
are
in
$IX.2.$
of
\cite{kato}.
}.
It is
known that this formula
holds for
\begin{enumerate}
\item
any
bounded 
operator $T\in B(G)$
on a complex 
Banach space $G$
with real spectrum $\sigma(T)$,
see for example
\cite{laursen};
\item
any 
infinitesimal generator
$T$
of a strongly
continuous semi-group in a Banach space,
see Corollary $8.1.16.$ of \cite{ds},
in particular
for
any unbounded self-adjoint operator
$T:\D(T)\subset H\to H$
in a complex Hilbert space $H$.
\end{enumerate}
 
Next we consider a more general case.
Let
$S$
be
an entire function
and $L>0$,
then 
the 
Newton-Leibnitz
formula
\begin{equation}
\label{15061756A}
R
\int_{u_{1}}^{u_{2}}
\frac{d\,S}{d\,\lambda}(t R)
\,d\,t
=
S(u_{2} R)
-
S(u_{1} R),
\end{equation}
for all
$u_{1},u_{2}\in [-L,L]$
was known
for
any element
$R$ in a Banach algebra $\A$,
where 
$S(tR)$
and
$\frac{d S}{d\lambda}(tR)$
are understood
in the
standard framework of analytic 
functional calculus on Banach
algebras,
while
the
integral
is with
respect to the Lebesgue measure 
in the norm topology on $\A$
see for example 
\cite{Rudin1,Dieud1,Schw}.
If
$E$
is the resolution of the identity of 
$R$
then
for all
$U\in\B(\C)$
\footnote{
$\B(\C)$
is the class of all
Borelian sets of $\C$.
}
$$
\Lf{E}{\infty}(U)
\doteqdot
\left\{
f:\C\to\C
\mid
\|f\chi_{U}\|_{\infty}^{E}
<\infty
\right\}.
$$
Here
$\chi_{U}:\C\to\C$
is equal to $1$ in $U$
and $0$ in $\complement U$
and for all maps
$F:\C\to\C$
$$
\|F\|_{\infty}^{E}
\doteqdot
E-ess\sup_{\lambda\in\C}
|F(\lambda)|
\doteqdot
\inf_{\{\delta\in\B(\C)\mid E(\delta)=\un\}}
\sup_{\lambda\in\delta}
|F(\lambda)|.
$$
See \cite{ds}.
 
We say
(see 
Definition
\ref{13481501})
that
$\n$
is
an
\textbf{
$E-$appropriate set 
}
if
\begin{enumerate}
\item
$
\n
\subseteq 
B(G)^{*}
$
linear subspace;
\item
$\n$ 
separates the points 
of $B(G)$,
namely
$$
(\forall T\in B(G)-\{\ze\})
(\exists\,\omega\in\n)
(\omega(T)\neq 0);
$$
\item
$
(\forall\omega\in\n)
(\forall\sigma\in\B(\C))
$
we have
\begin{equation}
\label{14111501int}
\omega\circ
\RM(E(\sigma))
\in\n
\text{ and }
\omega\circ
\mathcal{L}(E(\sigma))
\in\n.
\end{equation}
\end{enumerate}
Moreover,
we say that
$\n$
is
an
\textbf{
$E-$appropriate set
with the
duality property
}
if
in addition
\begin{equation}
\label{17471801}
\n^{*}
\subseteq 
B(G).
\end{equation}
Here
for any
Banach algebra
$\A$,
so in particular
for
$\A=B(G)$,
we set
$
\RM:\A\to\A^{\A}
$
and
$
\mathcal{L}:\A\to\A^{\A}
$
defined
by
\begin{equation}
\label{16051546int}
\begin{cases}
\RM(T):
\A\ni h\mapsto T h\in\A
\\
\mathcal{L}(T):
\A\ni h\mapsto h T\in\A,
\end{cases}
\end{equation}
for all $T\in\A$.
Notice that
for all
$T,h\in\A$
we have
$
\|
\RM(T)
(h)
\|_{\A}
\leq
\|T\|_{\A}
\|h\|_{\A}
$,
and
$
\|
\mathcal{L}(T)(h)
\|_{\A}
\leq
\|T\|_{\A}
\|h\|_{\A}
$,
so
\begin{equation}
\label{10501634}
\RM(T),
\mathcal{L}(T)
\in
B(\A)
\end{equation}
with
\begin{equation}
\label{13221607}
\|
\RM(T)
\|_{B(\A)} 
\leq
\|T\|_{\A},
\|
\mathcal{L}(T)
\|_{B(\A)} 
\leq
\|T\|_{\A}.
\end{equation}
Since
$
\mathcal{L}
$
and
$
\RM
$
are linear mappings
we can conclude that
\begin{equation}
\label{16051547}
\begin{cases}
\mathcal{L},
\RM
\in
B(\A,B(\A))
\\
\|\RM\|_{B(\A,B(\A))},
\|\mathcal{L}\|_{B(\A,B(\A))}
\leq
1.
\end{cases}
\end{equation}
In
\eqref{17471801}
we mean
$$
(\exists\,Y_{0}\subseteq B(G))
(\n^{*}
=
\{
\hat{A}
\up
\n
\mid
A\in Y_{0}
\}),
$$
where
$
\left(\hat{\cdot}\right):
B(G)
\to
(B(G)^{*})^{*}
$
is the canonical
isometric embedding
of $B(G)$
into
its bidual.
 
In the work the 
following generalizations of 
\eqref{15061756A}
are proved for the case when 
$R:\D\subset G\to G$
is an
unbounded
scalar type spectral 
operator
in a 
complex
Banach space
$G$,
in particular when
$R:\D\subset H\to H$
is an unbounded self-adjoint operator
in a complex Hilbert space $H$.
Under the assumptions that
$
S:U\to\C
$
is an analytic
map
on an 
open
neighbourhood
$U$
of 
the spectrum 
$\sigma(R)$
of
$R$
such 
that
there is
$L>0$
such that
$]-L,L[\cdot U\subseteq U$
and
$$
\w{S_{t}}
\in
\Lf{E}{\infty}(\sigma(R)),\,
\w{
\left(
\frac{d\,S}{d\,\lambda}
\right)_{t}
}
\in
\Lf{E}{\infty}(\sigma(R))
$$
for all 
$
t\in 
]-L,L[
$,
where
$
(S)_{t}(\lambda)
\doteqdot
S(t\lambda)
$
and
$
(\frac{d\,S}{d\,\lambda})_{t}(\lambda)
\doteqdot
\frac{d\,S}{d\,\lambda}(t\lambda)
$
for all 
$t\in]-L,L[$
and
$\lambda\in U$,
while
for any map
$F:U\to\C$
we set
$\w{F}$
the $\ze-$extension
of $F$
to $\C$.
The following statements are proved.
\begin{enumerate}
\item
If
\begin{equation}
\label{15032001}
\int^{*}
\left\|
\w{
\left(
\frac{d\,S}{d\,\lambda}
\right)_{t}
}
\right\|_{\infty}^{E}
\,
d\,t
<\infty
\end{equation}
and for all
$\omega\in\n$
the
map
$
]-L,L[
\ni
t
\mapsto
\omega
\left(
\frac{d\,S}{d\,\lambda}
(tR)
\right)
\in\C
$
is
Lebesgue
measurable,
then
in
Corollary
\ref{20051321ta}
it is proved
that
formula
\eqref{15061756A}
holds
where
the integral 
is 
the weak-integral
\footnote{
See formula
\eqref{17421502}
below.
}
with respect
to the
Lebesgue measure
and
with respect
to the 
$\sigma(B(G),\n)-$topology
for any
$E-$appropriate set 
$\n$
with the duality property.
Moreover
in
Corollary
\ref{20051321taLOC}
it is proved
that
formula
\eqref{15061756A}
also
holds
when
$
\w{
\left(\frac{d\,S}{d\,\lambda}
\right)_{t}
}
\in
\Lf{E}{\infty}(\sigma(R))
$
almost
everywhere
on 
$]-L,L[$
with respect
to the Lebesgue
measure.
\item
In particular
it is proved
that
formula
\eqref{15061756A}
holds
where
the integral 
is 
the weak-integral 
with respect
to the
Lebesgue measure
and
with respect
to the 
sigma-weak operator 
topology,
when $G$ is a Hilbert space
(Corollary
\ref{20051321pd}).
\item
If
in addition to 
\eqref{15032001},
$G$ 
is
a
reflexive
complex
Banach space
then
in
Corollary
\ref{19371401}
it is proved
that
formula
\eqref{15061756A}
holds
where
the integral 
is 
the weak-integral 
with respect
to the
Lebesgue measure
and
with respect
to the 
weak operator
topology.
\item
If
\begin{equation}
\label{16460310}
\sup_{t\in ]-L,L[}
\left\|
\w{
\left(
\frac{d\,S}{d\,\lambda}
\right)_{t}
}
\right\|_{\infty}^{E}
<
\infty, 
\end{equation}
then
in
Theorem
\ref{27051755}
it is proved
that
formula
\eqref{15061756A}
holds
where
the integral 
is 
with respect
to the
Lebesgue measure
and
with respect
to the 
strong
operator
topology.
\item
In
Theorem
\ref{27051053}
it is proved that
if in addition
to the
\eqref{16460310}
$$
\sup_{t\in ]-L,L[}
\|\w{(S)_{t}}\|_{\infty}^{E}
<
\infty,
$$
then for all
$
v\in\D
$
the mapping
$
]-L,L[
\ni
t
\mapsto
S(t R)
v
\in 
G
$
is
differentiable,
and
$
(\forall v\in\D)
(\forall t\in ]-L,L[)
$
\begin{equation}
\label{11400510}
\frac{d S(t R)v}{d t}
=
R
\frac{d S}{d\lambda}
(t R)
v.
\end{equation}
\label{16570310D}
\item
In 
Corollary
\ref{04041425}
formula 
\eqref{03101514}
is deduced from
formula
\eqref{15061756A}
for any
unbounded scalar type spectral 
operator
$T:\D(T)\subset G\to G$
in a complex Banach space
$G$
with real spectrum.
\end{enumerate}
In these statements 
$
\frac{d S}{d\lambda}
(t R)
$
and
$S(tR)$
are
understood
in the framework of the 
Borel 
functional 
calculus
for 
unbounded
scalar type
spectral operators
in $G$.
See 
definition
$18.2.10.$
in 
Vol $3$
of the
Dunford-Schwartz
monograph
\cite{ds}
(also
see
Section \ref{09061650} of the work).
 
In order to prove equality
\eqref{15061756A}
when $R$ is an unbounded 
scalar type spectral operator in $G$,
we procede in two steps.
First of all
we consider the Banach spaces
$G_{\sigma_{n}}\doteqdot E(\sigma_{n})G$
where 
$\sigma_{n}\doteqdot B_{n}(\ze)\subset\C$,
with $n\in\N$,
the bounded operators 
$R_{\sigma_{n}}\doteqdot RE(\sigma_{n})$,
and their restrictions 
$
(R_{\sigma_{n}}\up G_{\sigma_{n}})
$
to 
$G_{\sigma_{n}}$.
Then 
by Key 
Lemma 
\ref{II31031834}
the operators
$
R_{\sigma_{n}}\up G_{\sigma_{n}}
$
are bounded
\emph{scalar type spectral} 
operators
on $G_{\sigma_{n}}$, 
and
for all $x\in G$
\begin{equation}
\label{12382312}
S(R)x
=
\lim_{n\in\N}
S(R_{\sigma_{n}}\up G_{\sigma_{n}})
E(\sigma_{n})x,
\end{equation}
in $G$
and
\begin{equation}
\label{12392312}
(R_{\sigma_{n}}\up G_{\sigma_{n}})
\int_{u_{1}}^{u_{2}}
\frac{d\,S}{d\,\lambda}
(t (R_{\sigma_{n}}\up G_{\sigma_{n}}))
\,d\,t
=
S(u_{2}(R_{\sigma_{n}}\up G_{\sigma_{n}}))
-
S(u_{1}(R_{\sigma_{n}}\up G_{\sigma_{n}})).
\end{equation}
The second and
most important step
it is to set up a 
\emph{limiting}
procedure,
which allows 
by using
the convergence
\eqref{12382312}
to extend
the 
``local''
equality
\eqref{12392312}
to the 
``global''
one
\eqref{15061756A}.
 
As we shall see below such a 
limiting procedure can be 
set up in a very general context.
First 
we wish
to
point 
out 
that
the following 
equalities
for all
$n\in\N$
and 
$t\in ]-L,L[$,
which follow from 
Key Lemma \ref{II31031834}
are essential
for making 
this limiting procedure possible
\begin{equation}
\label{12492312a}
\begin{cases}
\frac{d\,S}{d\,\lambda}
(t R)
E(\sigma_{n})
=
\frac{d\,S}{d\,\lambda}
(t (R_{\sigma_{n}}\up G_{\sigma_{n}}))
E(\sigma_{n}),
\\
S(t R)
E(\sigma_{n})
=
S(t (R_{\sigma_{n}}\up G_{\sigma_{n}}))
E(\sigma_{n}).
\end{cases}
\end{equation}
 
We note that
one 
cannot 
replace in 
\eqref{12392312}
$R_{\sigma_{n}}\up G_{\sigma_{n}}$
with the simpler operator
$R_{\sigma_{n}}$
for the following reason.
Although 
$R_{\sigma_{n}}$
is a bounded operator on $G$
for $n\in\N$
and
$
Rx
=
\lim_{n\in\N}
R_{\sigma_{n}}x
$
in $G$,
in general 
$R_{\sigma_{n}}$
is not a scalar type spectral operator,
hence the expression
$
\frac{d\,S}{d\,\lambda}
(t R_{\sigma_{n}})
$
is not defined
in the 
Dunford-Schwartz Functional
Calculus for scalar 
type spectral operators,
which turns to be mandatory
in the sequel when using
general Borelian maps not
necessarily analytic.
 
Next we 
formulate 
a rather general statement
allowing, by using a limiting procedure, to
pass from ``local'' equalities similar to
\eqref{12392312}
to ``global'' ones similar to 
\eqref{15061756A}.
 
We generalize 
\eqref{15061756A}
in several directions.
We replace 
\begin{itemize}
\item
the operator
$R$ to the left of
the integral
by a function
$g(R)$,
where $g$ is a general Borelian map on 
$\sigma(R)$
\footnote{
The most interesting case is when the
operator $g(R)$ is unbounded.
},
\item
the compact set
$
[u_{1},u_{2}]
$ 
and 
the Lebesgue measure
on it
by
a general locally compact space $X$
and 
a complex
Radon measure on it
respectively,
\item
the
map
$
[u_{1},u_{2}]
\ni
t
\to
\left(\frac{dS}{d\lambda}\right)_{t}
\in
Bor(\sigma(R))
$
by
the map
$
X\ni x\to f_{x}\in Bor(\sigma(R))
$
such that
$
\w{f}_{x}\in\Lf{E}{\infty}(\sigma(R))
$
where
$\w{f}_{x}$ 
is the $\ze-$extension
to $\C$ of $f_{x}$,
and
the map
$
X\ni x\to f_{x}(R)\in B(G)
$
is 
strongly integrable with 
respect to the measure $\mu$;
\footnote{
This
means that
$
X\ni x\to f_{x}(R)v\in G
$
is
integrable with respect to the 
measure
$\mu$ for all $v\in G$,
in the sense of Ch $4$, \S $4$
of Bourbaki \cite{IntBourb},
and
the map
$
G\ni v\mapsto \int f_{x}(R)v\in G
$
is a (linear) bounded operator on $G$.
}
\item
the map
$
S_{u_{2}}-S_{u_{1}}
$
by
a Borelian map
$h$
on $\sigma(R)$
such that
$\w{h}\in\Lf{E}{\infty}(\sigma(R))$.
\end{itemize}
One
of the main
results
of the work
is
Theorem \ref{13051634}
where
we prove
that
if 
$\{\sigma_{n}\}_{n\in\N}$
is
an
$E-$sequence
\footnote{
By definition
this means
that for all
$n\in\N$
$\sigma_{n}\in\B(\C)$,
for all
$n,m\in\N$
$n>m
\Rightarrow
\sigma_{n}
\supseteq
\sigma_{m}
$;
$
\textrm{supp}(E)
\subseteq
\bigcup_{n\in\N}
\sigma_{n}
$;
hence 
we have
$
\lim_{n\in\N}E(\sigma_{n})
=
\un
$
strongly.
},
and
\footnote{
By Key Lemma
\ref{II31031834}
$R_{\sigma_{n}}\up G_{\sigma_{n}}$
is a scalar type spectral operator
in the complex Banach space $G_{\sigma_{n}}$,
but
in contrast 
to the previous
case 
where $\sigma_{n}\doteqdot B_{n}(\ze)$ 
was bounded,
here $\sigma_{n}$
could be unbounded so it may happen
that 
$
G_{\sigma_{n}}
\nsubseteq 
Dom(R)
$
hence
the restriction
$R_{\sigma_{n}}\up G_{\sigma_{n}}$
of
$R_{\sigma_{n}}$
to 
$G_{\sigma_{n}}$
has to be defined
on the set 
$G_{\sigma_{n}}\cap Dom(R)$,
and it could be an unbounded 
operator in $G_{\sigma_{n}}$
}
for all
$n\in\N$
\begin{equation}
\label{10440501}
R_{\sigma_{n}}\up G_{\sigma_{n}}
\doteqdot
RE(\sigma_{n})\up(G_{\sigma_{n}}\cap Dom(R)),
\end{equation}
and for all
$n\in\N$
the following
\emph{local} 
inclusion
\begin{equation}
\label{08530810}
g(R_{\sigma_{n}}\up G_{\sigma_{n}})
\int\,
f_{x}(R_{\sigma_{n}}\up G_{\sigma_{n}})
\,d\,
\mu(x)
\subseteq
h(R_{\sigma_{n}}\up G_{\sigma_{n}})
\end{equation}
holds, 
then 
$
h(R)
\in
B(G)
$
and
the 
\emph{global}
equality holds, i.e.
\begin{equation}
\label{22000310}
g(R)
\int\,
f_{x}(R)
\,d\,
\mu(x)
=
h(R).
\end{equation}
Here
all the integrals 
are with respect to the
strong
operator topology.

Now we can describe
Extension
Theorem
and the 
Newton-Leibnitz
formula for the 
integration
with respect
to the 
$
\sigma(B(G),\n)-
$
topology,
where
$\n$
is a suitable subset 
of $B(G)^{*}$,
which, 
roughly speaking,
is 
the
weakest
among
reasonable
locally
convex 
topologies
on
$B(G)$,
for which
the 
aforementioned
limiting
procedure
can
be
performed.
 
In Section
\ref{SecWEAKINT}
we 
recall the 
definition
of 
scalar essential 
$\mu-$integrability
and the
weak-integral
of maps
defined on $X$
and
with values in a 
Hausdorff
locally convex spaces,
where $\mu$ is a 
Radon measure on a locally
compact space $X$.
 
Here we need just to 
apply
these definitions
to the 
case
of 
$\sigma(B(G),\n)$,
i.e.
the 
weak 
topology
on $B(G)$
defined
by the standard
duality
between
$B(G)$
and
$\n$
where
$\n$
is a subset
of the 
(topological)
dual
$B(G)^{*}$
of 
$B(G)$
such that
it
separates
the points
of
$B(G)$.
 
Thus
$
f:X\to
\lr{B(G)}{\sigma(B(G),\n)}
$
is
by definition
scalarly essentially 
$\mu-$integrable 
or 
equivalently
$
f:X\to B(G)
$
is
scalarly 
essentially 
$\mu-$integrable 
with
respect to the 
measure
$\mu$
and
with respect
to
the
$
\sigma(B(G),\n)
$
topology
on $B(G)$
if
for
all
$\omega\in\n$
the 
map
$
\omega\circ f:X\to\C
$
is
essentially 
$\mu-$integrable
\footnote{
See for the definition
Ch. $5$, 
$\S1$, $n^{\circ}3$,
of \cite{IntBourb}
},
so
we can define
its
\emph{integral}
as 
the following linear 
operator
$$
\n
\ni
\omega
\mapsto
\int
\omega(f(x))
d\,\mu(x)
\in\C.
$$
 
Let
$
f:X\to
\lr{B(G)}{\sigma(B(G),\n)}
$
be
scalarly essentially 
$\mu-$integrable 
and
assume that
\begin{equation}
\label{16021301int}
(\exists\,B\in B(G))
(\forall\omega\in\n)
\left(
\omega
(B)
=
\int
\omega(f(x))
d\,\mu(x)
\right).
\end{equation}
Notice
that
the operator
$B$
is defined by this condition
uniquely.
In this case, by definition
\emph{
$
f:X\to
\lr{B(G)}{\sigma(B(G),\n)}
$
is scalarly essentially 
$(\mu,B(G))-$integrable
}
or
\emph{
$
f:X\to B(G)
$
is scalarly essentially 
$(\mu,B(G))-$integrable
with respect to the 
$\sigma(B(G),\n)-$
topology
}
and
its
\emph{weak-integral
with respect to the 
measure $\mu$
and with respect to the
$\sigma(B(G),\n)-$
topology}
or simply
its
\textbf{weak-integral},
is defined
by
\begin{equation}
\label{17421502}
\int
f(x)
d\,\mu(x)
\doteqdot
B.
\end{equation}
Next we
can state
\textbf{
Theorem
\ref{18051958ta}
},
the main 
result
of the work.
\begin{theorem}
[
\textbf{
$\sigma(B(G),\n)-$
Extension Theorem
}
]
Let
$G$
be
a complex Banach
space,
$X$
a locally compact space
and
$\mu$ a complex
Radon
measure on it.
In addition
let
$R$
be
a 
possibly unbounded 
scalar type spectral operator
in $G$,
$\sigma(R)$
its spectrum,
$E$ its resolution of 
the
identity
and
$\n$
an
$E-$appropriate set.
Let
the map
$
X
\ni 
x
\mapsto
f_{x}
\in
Bor(\sigma(R))
$
be
such that
$
\w{f}_{x}
\in
\Lf{E}{\infty}(\sigma(R))
$\,
$
\mu-
$
locally almost everywhere
on $X$
and the map
$
X
\ni x\mapsto f_{x}(R)
\in 
\lr{B(G)}{\sigma(B(G),\n)}
$
be
scalarly
essentially
$(\mu,B(G))-$integrable.
Finally
let
$
g,
h\in
Bor(\sigma(R))
$
and
$
\w{h}\in\Lf{E}{\infty}(\sigma(R))
$.
 
If
$\{\sigma_{n}\}_{n\in\N}$
is
an
$E-$sequence
and for all
$n\in\N$
\begin{equation}
\label{19211001taint}
g(R_{\sigma_{n}}\up G_{\sigma_{n}})
\int\,
f_{x}(R_{\sigma_{n}}\up G_{\sigma_{n}})
\,d\,
\mu(x)
\subseteq
h(R_{\sigma_{n}}\up G_{\sigma_{n}})
\end{equation}
then
$h(R)
\in
B(G)
$
and
\begin{equation}
\label{16011401taint}
g(R)
\int\,
f_{x}(R)
\,d\,
\mu(x)
=
h(R).
\end{equation}
In
\eqref{19211001taint}
the weak-integral is
with respect to the 
measure $\mu$
and with respect to the
$\sigma(B(G_{\sigma_{n}}),\n_{\sigma_{n}})-$
topology
\footnote{
$\n_{\sigma_{n}}$ is, 
roughly speaking,
the set of the
restrictions
to $B(G_{\sigma_{n}})$
of all the functionals
belonging
to
$\n$.
For the exact 
definition
and properties
see Definition
\ref{10061745}
and
Lemma
\ref{13061030}.
},
while
in
\eqref{16011401taint}
the weak-integral is
with respect to the 
measure $\mu$
and with respect to the
$\sigma(B(G),\n)-$
topology.
\end{theorem}
Notice that
$g(R)$ is 
a
possibly 
\textbf{
unbounded
}
operator
in
$G$.
 
We 
list 
the 
most important
results
that
allow
to prove
Theorem
\ref{18051958ta}:
\begin{enumerate}
\item
Key
Lemma
\ref{II31031834};
\item
``Commutation'' 
property
(Theorem
\ref{18051509ta}):
\begin{equation}
\label{18051301int}
\forall\sigma\in\B(\C)
\,
\left[
\int
f_{x}(R)
d\,\mu(x),
\,
E(\sigma)
\right]
=
\ze;
\end{equation}
\item
``Restriction''
property
(Theorem
\ref{14050121}):
for all
$\sigma\in\B(\C)$
we have 
that
$
f_{x}(R_{\sigma}\up G_{\sigma})
\in
B(G_{\sigma})
$,
$
\mu-
$
locally almost everywhere
on $X$,
$
X\ni x
\mapsto
f_{x}(R_{\sigma}\up G_{\sigma})
\in
\lr{B(G_{\sigma})}
{\sigma(B(G_{\sigma}),\n_{\sigma})}
$
is
scalarly essentially
$(\mu,B(G_{\sigma}))-$integrable,
and
\begin{equation}
\label{restrtaint}
\int
f_{x}(R_{\sigma}\up G_{\sigma})
\,d\mu(x)
=
\int
f_{x}(R)
\,d\mu(x)
\up
G_{\sigma};
\end{equation}
\item
finally the fact that
$$
Dom
\left(
g(R)
\int\,
f_{x}(R)
\,d\,
\mu(x)
\right)
\text{ is dense in $G$. }
$$
\end{enumerate}
 
We remark
that
the reason for introducing
the concept
of an
$E-$appropriate set
is primarily 
for obtaining
the commutation and 
restriction properties.
 
Now
we 
define
\begin{equation}
\label{17051313int}
\n_{st}(G)
\doteqdot
\lr{B(G)}{\tau_{st}(G)}^{*}
=
\mathfrak{L}_{\C}
(\{
\psi_{(\phi,v)}
\mid
(\phi,v)\in G^{*}\times G
\}).
\end{equation}
Here
$\lr{B(G)}{\tau_{st}(G)}^{*}$
is the topological
dual
of 
$B(G)$
with respect
to the strong operator
topology,
$
\psi_{(\phi,v)}:
B(G)
\ni
T
\mapsto
\phi(T v)
\in
\C
$,
while
$\mathfrak{L}_{\C}(J)$
is the complex linear space
generated by the set 
$J\subseteq B(G)^{*}$.
Then
$\sigma(B(G),\n_{st}(G))$
is
the weak operator
topology
on $B(G)$
and
$\n_{st}(G)$
is an
$E-$appropriate set
for
any spectral measure
$E$.
 
Moreover
we set
in the case in which
$G$ is a complex Hilbert space
$$
\n_{pd}(G)
\doteqdot
\text{ predual of }
B(G),
$$
which is 
by definition
the linear space
of all sigma-weakly
continuous linear functionals
on
$B(G)$.
 
Note that
\begin{equation}
\label{18591801}
\n_{pd}(G)^{*}
=
B(G).
\end{equation}
(See statement 
$(iii)$
of
Theorem $2.6.$,
Ch. $2$
of 
\cite{tak},
or
Proposition $2.4.18$
of \cite{bra}).
Here
we mean
that the normed subspace
$\n_{pd}(G)^{*}$
of the bidual
$(B(G)^{*})^{*}$
is 
isometric 
to 
$B(G)$,
through
the canonical 
embedding
of
$B(G)$
into
$(B(G)^{*})^{*}$.
 
Hence
we can apply
the
Extension Theorem \ref{18051958ta}
to the case
$\n\doteqdot
\n_{st}(G)
$,
or
$\n\doteqdot
\n_{pd}(G)
$
and 
use the 
following
additional
property
which is proved
in
Proposition
\ref{15051146}
\begin{equation}
\label{13061209Aint}
(\n_{st}(G))_{\sigma}
=
\n_{st}(G_{\sigma}),
\text{ and }
(\n_{pd}(G))_{\sigma}
=
\n_{pd}(G_{\sigma}).
\end{equation}
 
The reason of introducing
the concept
of
duality property
for
$E-$appropriate set
is primarly for assuring
that
a map
$f:X\to\lr{B(G)}{\sigma(B(G),\n)}$
scalarly 
essentially
$\mu-$integrable
is
also
$(\mu,B(G))-$integrable.
 
As an application
of
this fact 
and of
the Extension theorem 
we obtain
the
Newton-Leibnitz
formula
in \eqref{15061756A}
by replacing
$\A$ with $B(G)$,
$R$ with an unbounded
scalar type spectral operator
in a complex Banach space
$G$,
by considering 
$S$
analytic
in an open
neighbourhood
$U$
of $\sigma(R)$
such
that
$]-L,L[\cdot
U\subseteq U$,
and
the
integral
with respect
to the 
$\sigma(B(G),\n)-$topology,
where
$\n$
is
an 
$E-$appropriate set
with the duality
property
(\textbf{
Corollary
\ref{20051321ta}}).
 
Finally
in a similar
way
we obtain
the corresponding
results
for the cases
of
the
sigma-weak operator
topology
(Corollary
\ref{20051321pd}),
and
for the cases
of
weak operator
topology
(Corollary
\ref{19371401}).
The last
result
is 
a complement
to 
Theorem
\ref{27051755}.
\chapter*{
Summary of the main results
}
Let
$G$
be
a complex Banach space,
$R$
an unbounded scalar type
spectral operator in 
$G$,
for example
an
unbounded self-adjoint
operator in a Hilbert space,
$\sigma(R)$
its spectrum
and
$E$
its
resolution of identity.
The
\textbf{main results}
of the work are
the following ones.
\begin{enumerate}
\item
Extension procedure leading from local
equality 
\eqref{19211001taint}
to global equality
\eqref{16011401taint}
for
integration
with respect to the
$\sigma(B(G),\n)-$topology
(Theorem
\ref{18051958ta}
if
$\n$ is an
$E-$appropriate set
and
Corollary
\ref{17070917TA}
if
$\n$ is an
$E-$appropriate set
with the duality property).
\item
Extension procedure leading from local
equality 
\eqref{19211001taint}
to global equality
\eqref{16011401taint}
for
integration
with respect to the
sigma-weak topology
(
Corollary
\ref{12121201}
and
Theorem
\ref{17070917pd})
and
for
integration
with respect to the
weak operator topology
(Corollary
\ref{18051958}
and
Theorem
\ref{17070917}
or
Theorem \ref{13051634}
and
Corollary \ref{17070901}).
\item
Newton-Leibnitz
formula
\eqref{15061756A}
for 
a suitable
analytic 
map
$S$
for 
integration
with respect
to the
$
\sigma(B(G),\n)-
$
topology,
where
$\n$
is
an 
$E-$appropriate set
with the duality
property
(Corollary
\ref{20051321ta}
and
Corollary
\ref{20051321taLOC});
for 
integration
with respect to the
sigma-weak topology
(Corollary
\ref{20051321pd})
and
for 
integration
with respect to the
weak operator topology
(Corollary
\ref{19371401}
and
Theorem
\ref{27051755}).
\item
Differentiation formula 
\eqref{11400510}
for
a suitable
analytic 
map
$S$
(
Theorem
\ref{19500603}
and
Theorem
\ref{27051053}).
\item
A new proof
for the
resolvent formula
\eqref{03101514}
via 
formula
\eqref{15061756A}
(Corollary
\ref{04041425}).
\end{enumerate}


\mainmatter


\chapter{
Extension
theorem.
The case of the
strong operator topology
}
\label{31050841}
\section{
Key lemma
}
\label{09061650}
\begin{preliminaries}
\label{II01041232}

\textbf{
Integrals of bounded Borelian
functions with 
respect to a vector valued measure.
}
\begin{normalfont}
In 
the sequel
$
G
\doteqdot
\lr{G}{\|\cdot\|_{G}}
$
will be 
a complex 
Banach space.
Denote by
$\Pr(G)$
the class of all
projectors on $G$,
that is the class of $P\in B(G)$ 
such that $P^{2}=P$.
Consider
a Boolean algebra
$\B_{X}$,
see Sec. $1.12$ of \cite{ds}, 
of subsets of a set $X$, 
with respect to the
order relation
defined by
$
\sigma\geq\delta
\Leftrightarrow 
\sigma\supseteq\delta
$
and
complemented by
the operation
$
\sigma'
\doteqdot
\complement
\sigma
$.
In particular
$\B_{X}$
contains $\emptyset$ and $X$ and is closed 
under finite intersection and finite union.
 
The map
$
E:\B_{X}\to B(G)
$
is called
a spectral measure in $G$
on $\B_{X}$,
or simply on $X$ if $X$ is a 
topological space and
$\B_{X}$ is the Boolean algebra
of its Borelian subsets,
if
\begin{enumerate}
\item
$
E(\B_{X})\subseteq\Pr(G)
$;
\item
$
(\forall\sigma_{1},\sigma_{2}\in\B_{X})
(E(\sigma_{1}\cap\sigma_{2})=E(\sigma_{1})E(\sigma_{2}))
$;
\item
$
(\forall\sigma_{1},\sigma_{2}\in\B_{X})
(E(\sigma_{1}\cup\sigma_{2})=E(\sigma_{1})+
E(\sigma_{2})-E(\sigma_{1})E(\sigma_{2}))
$;
\item
$E(X)=\un$;
\item
$E(\emptyset)=\ze$.
\end{enumerate}
(See Definition 
$15.2.1$ of \cite{ds}).
 
If condition $(3)$
is replaced by condition
$$
(3')
(\forall\sigma_{1},\sigma_{2}\in\B_{X}
\mid
\sigma_{1}\cap\sigma_{2}=\emptyset)
(E(\sigma_{1}\cup\sigma_{2})=E(\sigma_{1})+
E(\sigma_{2})),
$$
we obtain an equivalent
definition.
 
Notice
that
if
$E$
is a spectral measure in $G$
on $\B_{X}$,
then 
it
is a
Boolean homomorphism
onto
the
Boolean 
algebra
$E(\B_{X})$
with respect to the
order relation
induced
by that
defined 
in
$\Pr(G)$
by
$
P\geq Z
\Leftrightarrow 
Z=Z P
$
and
complemented by
the operation
$
P'
\doteqdot
(\un-P)
$.
Indeed
for all
$\sigma,\delta\in\B_{X}$
we have
$
\delta
\subseteq
\sigma
\Rightarrow 
E(\delta)
=
E(\delta\cap\sigma)
\doteq
E(\delta)
E(\sigma)
\Leftrightarrow 
E(\delta)
\leq
E(\sigma)
$,
while
$
\un
=
E(\sigma\cup\complement\sigma)
=
E(\sigma)+E(\complement\sigma)
$.
 
A spectral measure $E$ is called 
\emph{
(weakly)
countable
additive
}
if
for all
sequences
$
\{\varepsilon_{n}\}_{n\in\N}\subset\B_{X}
$
of
disjoint sets,
for all
$x\in G$
and for all
$\phi\in G^{*}$
we have 
$$
\phi\left(E(\bigcup_{n\in\N}
\varepsilon_{n})x\right)
=
\sum_{n=1}^{\infty}\phi
\left(E(\varepsilon_{n})x\right).
$$
If $\B_{X}$
is a $\sigma$-field,
i.e.
a Boolean algebra
closed under the operation of 
forming countable 
unions,
we have
by
Corollary
$15.2.4.$
of the 
\cite{ds}
that
$E$
is 
countably additive 
with
respect to the 
strong operator topology,
i.e.
for all
sequence
$
\{\varepsilon_{n}\}_{n\in\N}\subset\B(\C)
$
of disjoint sets 
and
for all $x\in G$
we have
\footnote{
By definition,
see Ch.$3$ of \cite{BourGT},
$
v=
\sum_{n\in\N}E(\varepsilon_{n})x
$
if
$
v=
\lim_{J\in\p_{\omega}(\N)}
\sum_{n\in J}
E(\varepsilon_{n})x
$,
where
$\p_{\omega}(\N)$
is the direct ordered
set of all finite subsets
of $\N$
ordered by inclusion.
}
\begin{equation}
\label{18282801}
E(\bigcup_{n\in\N}\varepsilon_{n})x
=
\sum_{n=1}^{\infty}E(\varepsilon_{n})x
=
\sum_{n\in\N}E(\varepsilon_{n})x.
\end{equation}
Since
$
E(\bigcup_{n\in\N}\varepsilon_{n})
=
E(\bigcup_{n\in\N}\varepsilon_{\rho(n)})
$,
for any 
permutation
$\rho$
of $\N$,
hence
$
\sum_{n=1}^{\infty}E(\varepsilon_{n})x
=
\sum_{n=1}^{\infty}E(\varepsilon_{\rho(n)})x
$
for all $x\in G$,
therefore
by 
Proposition
$9$,\S $5.7.$, Ch. $3$
of \cite{BourGT}
we obtain
the second
equality
in
\eqref{18282801}.
By
$\B(\C)$ 
we denote
the set of the Borelian subsets of $\C$,
and
by
$Bor(U)$
the complex linear space of all 
Borelian complex maps 
defined on a Borelian
subset $U$ of $\C$.
We denote
with
$
\bb
$ 
the 
space
of the totally
$\B(\C)-$
measurable
maps
\footnote{
In \cite{ds}
denoted
by
$B(\C,\B(\C))$,
while
by using the notations
of
\cite{din2} 
and considering
$\C$ as a real Banach space
we have
$
\bb
=
\bb_{\R}(\B(\C))
$.},
which is
the closure 
in the Banach space 
$\lr{\mathbf{B}(\C)}{\|\cdot\|_{\sup}}$ 
of all complex 
bounded functions on $\C$
with
respect to the norm 
$\|g\|_{\sup}
\doteqdot
\sup_{\lambda\in\C}|g(\lambda)|$,
of the linear space generated by the set
$\{\chi_{\sigma}\mid\sigma\in\B(\C)\}$,
where
$\chi_{\sigma}$ 
is the characteristic 
function of the set $\sigma$.
$\lr{\bb}{\|\cdot\|_{\sup}}$ 
is a Banach space,
and
the
space of 
all bounded Borelian complex functions
is in
$\bb$
so dense in it.
Finally
$\lr{\bb}{\|\cdot\|_{\sup}}$
is a 
$C^{*}-$subalgebra, 
in particular
a Banach subalgebra, 
of
$\lr{\mathbf{B}(\C)}{\|\cdot\|_{\sup}}$  
if we define the pointwise
operations 
of product and involution
on $\mathbf{B}(\C)$.
 
Let 
$X$
be 
a complex 
Banach
space
and
$F:\B(\C)\to X$
a
weakly
countably
finite
additive
vector valued 
measure,
see Section $4.10.$
of \cite{ds},
then we can define the 
integral
with respect to $F$,
see 
Section $10.1$ 
of \cite{ds},
which will be
denoted by
$\int_{\C}
f\,d\,F$.
The operator
\begin{equation}
\label{19492902}
\I{\C}{F}:
\bb
\ni
f
\mapsto
\int_{\C}
f\,d\,F
\in
X
\end{equation}
is linear and norm-continuous
\footnote{
Notice
that
if we 
identify
$B(G)$
with
$B(\R,B(G))$
and
recall
that
$
\bb
=
\bb_{\R}(\B(\C))
$,
then
with 
the notations
of
Definition
$24$, \S $1$, Ch. $1$
of
\cite{din2} 
we have
that
$\I{\C}{E}$
is the 
immediate
integral
with respect
to the vector valued
measure
$E:\B(\C)\to B(\R,B(G))$.
}.
We have the
following 
useful property
if $Y$ is a $\C-$Banach space and 
$Q\in B(X,Y)$,
then
\begin{equation}
\label{20522502}
Q\circ\I{\C}{F}=\I{\C}{Q\circ F},
\end{equation}
see 
statement $(f)$
of
Theorem
$4.10.8.$
of the \cite{ds}.
 
If 
$X\doteqdot B(G)$,
the case we are mostly 
interested in,
we have,
as an immediate 
result
of this property 
and the fact 
that 
the map
$
Q_{x}:B(G)\ni A\mapsto A x\in G
$
is linear and continuous
for all
$x\in G$,
that
\begin{equation}
\label{01041859}
(\forall x\in G)
(\forall f\in \bb)
(\I{\C}{F}(f)x=\I{\C}{F^{x}}(f)).
\end{equation}
Here
$
F^{x}:\B(\C)\ni
\sigma
\mapsto
F(\sigma)x
$.
Finally
if $E$ is a spectral measure on $\C$, 
then
$\I{\C}{E}$
is a continuous 
unital
homomorphism
between the 
two
Banach 
algebras
$\lr{\bb}{\|\cdot\|_{\sup}}$,
and
$\lr{B(G)}{\|\cdot\|_{B(G)}}$
and
$\I{\C}{E}(\chi_{supp\, E})=\un$,
see 
\eqref{13422802}
and
Section $(2)$,
Ch $15$
of
\cite{ds}.
 
\textbf{
Borel functional calculus for 
possibly unbounded
scalar type spectral operators in $G$.
}
If 
$T:\DM(T)\subseteq G\to G$ 
is a
possibly unbounded 
linear operator
then we denote
by
$\sigma(T)$
its standard 
spectrum. 
A 
possibly unbounded 
linear operator
$T:\DM(T)\subseteq G\to G$ 
is called a
\textbf{spectral operator in $G$} if
it is closed
and there exists 
a countably additive spectral measure 
$E:\B(\C)\to\Pr(G)$ such that
\begin{description}
\item[i]
for all bounded sets
$\delta\in\B(\C)$
$$
E(\delta)G\subseteq\DM(T);
$$
\item[ii]
$
(\forall\delta\in\B(\C))
(\forall x\in\DM(T))
$
we have
\begin{enumerate}
\item
$(E(\delta)\DM(T)\subseteq\DM(T))$,
\item
$
T E(\delta) x
=
E(\delta) T x
$;
\end{enumerate}
\item[iii]
for all
$\delta\in\B(\C)$
we have
$$
\sigma\left(T\up(\DM(T)\cap E(\delta)G)\right)
\subseteq
\ov{\delta}.
$$
Here
$
\sigma\left(T\up(\DM(T)\cap E(\delta)G)\right)
$
is the spectrum 
of the restriction
of
$T$
to the 
domain
$\DM(T)\cap E(\delta)G$.
\end{description}
(See 
Definition $18.2.1.$ 
of the \cite{ds}).
We call any $E$ with the above
properties
a \textbf{resolution of the identity of $T$}.
Theorem $18.2.5.$ of \cite{ds} states that
the resolution of the identity of a 
spectral operator
is unique.
 
Finally we call support
of a spectral measure $E$ on $\B_{X}$,
the following set
$$
\textrm{supp}\, E
\doteqdot
\bigcap_{\{\sigma\in\B_{X}|E(\sigma)=\un\}}
\ov{\sigma}.
$$
It is easy
to see 
\footnote{
Indeed
let
$
S\doteqdot
\textrm{supp}\, E
$
then
\begin{equation}
\label{13592802}
\complement
S
=
\bigcup_{\{\sigma\in\B_{X}\mid E(\sigma)=\un\}}
\complement
\ov{\sigma}.
\end{equation}
Moreover
$E$
is order-preserving
so
for all
$\sigma\in\B_{X}$
such that
$E(\sigma)=\un$
we have
$
E(\complement\ov{\sigma})
\leq
E(\complement\sigma)
=
\un
-
E(\sigma)
=
\ze
$.
Hence
by the
definition
of the order
$
E(\complement\ov{\sigma})
=
E(\complement\ov{\sigma})
\ze
=
\ze
$.
Therefore
by the Principle of
localization
(Corollary, 
Ch $3$, \S $2$, $n^{\circ} 1$
of \cite{IntBourb})
which holds
also for vector measures
(footnote in
Ch $6$, \S $2$, $n^{\circ} 1$
of \cite{IntBourb})
we deduce by
\eqref{13592802}
that
$
E(\complement S)
=\ze
$.
Finally
$$
E(S)
=
\un
-
E(\complement S)
=
\un.
$$
}
that
\begin{equation}
\label{13422802}
E(\textrm{supp}\, E)
=
\un.
\end{equation}
Notice that
an
unbounded
spectral operator 
$T$
is closed by definition.
Now we will show that 
$T$
is also 
densely defined.
In fact 
if 
$E$
is the resolution of the identity
of 
$T$
and
if 
$\{\sigma_{n}\}_{n\in\N}\subset\B(\C)$ 
is
a
non
decreasing
sequence of 
Borelian
sets
such that
$
\sigma(T)
\subseteq
\bigcup_{n\in\N}\sigma_{n} 
$,
then by the 
strong countable additivity of $E$,
the fact that
$E(\sigma(T))=\un$
we can deduce 
$
\un
=
\lim_{n\in\N}
E(\sigma_{n})
$
in the strong operator topology
of $B(G)$,
see
\eqref{II03041401}.
Now 
we can choose
$\{\sigma_{n}\}_{n\in\N}$ 
such that
$
\sigma_{n}\doteqdot B_{n}(\ze)
\doteqdot
\{\lambda\in\C\mid|\lambda|<n\}
$,
or
$
\sigma_{n}
\doteqdot
W(\ze,2n)
\doteqdot
\{\lambda\in\C
\mid
|Re(\lambda)|<n,
|Im(\lambda)|<n,
\}
$.
But by the
property
$(i)$
of the
Definition $18.2.1.$ of \cite{ds},
we know
that
for all
bounded
sets
$
\sigma\in\B(\C)
$
we have
$
E(\sigma)G\subseteq Dom(T)
$.
Therefore we conclude that
for all
$v\in G$,
$
v=
\lim_{n\in\N}  
E(\sigma_{n})v$
and for all $n\in\N$,
$E(\sigma_{n})v\in Dom(T)$,
so
$Dom(T)$ is dense in $G$.
 
We want to remark
that
for each
possibly unbounded 
spectral operator
$T$ 
in
$G$
by denoting
with
$\sigma(T)$ its spectrum
and with
$E:\B(\C)\to\Pr(G)$ 
its
resolution of the identity,
we deduce
by
Lemma $18.2.25.$
of
\cite{ds}
that
$\sigma(T)$
is closed,
that 
$\textrm{supp}\, E =\sigma(T)$
so
by
\eqref{13422802}
$$
E(\sigma(T))
=
\un.
$$
Now
we will 
give 
the definition
of the Borel functional calculus for 
unbounded
spectral
operators in a complex 
Banach space $G$,
that
is
essentially
the same
as
in
Definition $18.2.10.$ of the \cite{ds}.
\begin{definition}
\label{12470108}
Let $X$ be a set,
$S\subset X$,
$V$ a vector space over $\K\in\{\R,\C\}$
and
$f:S\to V$.
Then we define
$\w{f}^{X}$,
or simply
$\w{f}$
when it 
doesn't
cause confusion,
to be the $\ze-$extension of $f$ to $X$,
i.e.
$\w{f}:X\to V$
such that
$
\w{f}\up S
=
f
$
and
$
\w{f}(x)
=
\ze
$
for all
$x\in (X-S)$,
where
$\ze$
is the
zero vector of $V$.
\end{definition}
\begin{definition}
\label{20282412}
[Borel Functional Calculus of $E$]
Assume that
\begin{enumerate}
\item 
$E:\B(\C)\to\Pr(G)$ 
is a countably additive 
spectral measure
and 
$S$
its support;
\item
$f\in Bor(S)$;
\item 
for all
$\sigma\subseteq\C$ 
we set
$f_{\sigma}:\C\to\C$ 
such that
$f_{\sigma}
\doteqdot
\w{f}
\cdot
\chi_{\sigma}
$;
\item
$\delta_{n}\doteqdot [-n,+n]$
and
$$
f_{n}
\doteqdot
f_{\overset{-1}{|f|}(\delta_{n})}.
$$
Here
$
(\forall\sigma\subseteq\C)
(\forall g:\DM\to\C)
(\overset{-1}{g}(\sigma)
\doteqdot
\{
\lambda\in\DM\mid g(\lambda)
\in
\sigma
\})
$.
\end{enumerate}
Of course 
$
f_{n}\in\bb
$
for all
$n\in\N$
so we can define the following 
operator in $G$
\begin{equation}
\label{II01041439}
\begin{cases}
Dom(f(E))
\doteqdot
\{
x\in G
\mid
\exists
\lim_{n\in\N}\I{\C}{E}(f_{n})x\}\\
(\forall x\in Dom(f(E)))
(f(E)x
\doteqdot
\lim_{n\in\N}\I{\C}{E}(f_{n})x).
\end{cases}
\end{equation}
Here all limits are considered in the space
$G$.
We call the map
$f\mapsto f(E)$
the 
\textbf{Borel functional calculus of 
the spectral measure
$E$}.
\end{definition}
In the case in which $E$ is the resolution
of the identity of a possibly unbounded 
spectral
operator
$T$,
recalling 
Lemma $18.2.25.$
of \cite{ds}
stating
that
$\sigma(T)$
is the support of $E$,
we 
can define
$f(T)\doteqdot f(E)$
for any map
$f\in Bor(\sigma(T))$
and
call the map
$$
Bor(\sigma(T))\ni f\mapsto f(T)
$$
the 
\textbf{Borel functional calculus of 
the operator
$T$}.
\begin{definition}
\label{19382412}
[$18.2.12.$ of \cite{ds}]
A 
\emph{spectral operator of scalar type in $G$}
or
a
\textbf{scalar type spectral operator in $G$}
is
a 
possibly unbounded linear
operator $R$ in $G$
such that
there exists
a countably additive 
spectral measure
$E:\B(\C)\to\Pr(G)$ 
with 
support
$S$
and
the property
$$
R=\imath(E).
$$
Here
$\imath:
S\ni\lambda\mapsto\lambda\in\C$,
and 
$\imath(E)$
is relative to the Borel functional
calculus 
of the
spectral measure
$E$.
We call 
$E$ a 
\textbf{
resolution of the identity
of $R$}.
\end{definition}
Let 
$R$
be
a scalar type spectral operator
in $G$
and $E$
a
resolution of the identity
of $R$,
then
we have the following statements
by \cite{ds}:
\begin{itemize}
\item
$T$ is a spectral operator
in $G$;
\item
$E$ is the resolution of the identity
of $T$ as spectral operator;
\item
$E$ is unique.
\end{itemize}
\begin{definition}
[\cite{ds}]
\label{20212412}
Let $E:\B(\C)\to Pr(G)$ be a countably
additive spectral measure and 
$U\in\B(\C)$, then 
the space of all 
\emph{$E-$essentially bounded}
maps is the following linear space
$$
\Lf{E}{\infty}(U)
\doteqdot
\left\{
f:\C\to\C
\mid
\|f\chi_{U}\|_{\infty}^{E}
<\infty
\right\}.
$$
Here
$\chi_{U}:\C\to\C$
is the characteristic map
of $U$ which is by definition
equal to $1$ in $U$
and $0$ in $\complement U$,
and for each map
$F:\C\to\C$
$$
\|F\|_{\infty}^{E}
\doteqdot
E-ess\sup_{\lambda\in\C}
|F(\lambda)|
\doteqdot
\inf_{\{\delta\in\B(\C)\mid E(\delta)=\un\}}
\sup_{\lambda\in\delta}
|F(\lambda)|.
$$
\end{definition}
For a Borelian map
$f:U\supset\sigma(R)\to\C$,
with
$U\in\B(\C)$,
we define
$f(R)$
to be the 
operator
$
(f\up\sigma(R))(R)
$.
Let 
$g:U\subseteq\C\to\C$ 
be 
a Borelian map.
Then 
$g$ is $E-$essentially bounded
if
$$
E-ess\sup_{\lambda\in U}
|g(\lambda)|
\doteqdot
\|\w{g}\|_{\infty}^{E}
<
\infty.
$$
See
Definition
$17.2.6.$
of \cite{ds}.
One 
formula
arising
by
statement 
$(i)$
of the Spectral Theorem $18.2.11.$ of the \cite{ds},
which will be used many times in the work
is the following:
for all
Borelian complex function
$f:\sigma(R)\to\C$ 
and
for all
$\phi\in G^{*}$
and
$y\in Dom(f(R))$
\begin{equation}
\label{01041949}
\phi\left(f(R)y\right)
=
\int_{\C}
\w{f}\,d\,E_{(\phi,y)}.
\end{equation}
Here
$G^{*}$
is the topological dual of $G$,
that is the normed space of 
all $\C-$linear and continuous functionals 
on $G$ with the $\sup-$norm,
and for all
$\phi\in G^{*}$
and 
$y\in G$
we define
$
E_{(\phi,y)}:\B(\C)
\ni
\sigma
\mapsto
\phi(E(\sigma)y)
\in
\C
$.
Finally
if
$P\in\Pr(G)$
then
$\lr{P(G)}{\|\cdot\|_{P(G)}}
$,
with
$\|\cdot\|_{P(G)}\doteqdot\|\cdot\|_{G}\up P(G)$,
is
a
Banach space.
In fact
let
$\{v_{n}\}_{n\in\N}\subset G$
be
such that
$
v
=
\lim_{n\in\N}
P
v_{n}
$,
in
$\|\cdot\|_{G}$,
so 
$P=P^{2}$ 
being
continuous
we have that
$
P
v
=
\lim_{n\in\N}
P^{2}
v_{n}
=
\lim_{n\in\N}
P
v_{n}
\doteqdot
v
$,
so
$v\in P(G)$,
then
$P(G)$
is closed in 
$\lr{G}{\|\cdot\|_{G}}$ ,
hence
$\lr{P(G)}{\|\cdot\|_{P(G)}}$ 
is a 
Banach space.
If $E:\B_{Y}\to\Pr(G)$
is a spectral measure
in
$G$
on
$\B_{Y}$
and
$\sigma\in\B_{Y}$,
then
we 
shall denote by
$G_{\sigma}^{E}$
or simply
$G_{\sigma}$
the complex
Banach
space
$E(\sigma)G$,
without
expressing
its dependence
by $E$
whenever it does not cause confusion.
In addition
for any
$Q$
possibly unbounded operator
in $G$
we define for all
$\sigma\in\B_{Y}$
the following
possibly unbounded operator
operator
in $G$
$$
Q_{\sigma}
\doteqdot 
Q E(\sigma).
$$
Finally we shall 
denote by
$\B_{b}(\C)$
the subclass 
of all
bounded subsets
of
$\B(\C)$.
\end{normalfont}
\end{preliminaries}
\begin{definition}
\label{12051723}
Let
$F$
be
a
$\C-$Banach space,
$P\in\Pr(F)$
and
$S:Dom(S)\subseteq F\to F$,
then
we define
\begin{equation}
\label{16450710}
S P
\up P(F)
\doteqdot
S P
\up 
(P(F)\cap Dom(SP)).
\end{equation}
Notice that
by
the property
$P^{2}=P$
we have 
$
P(F)\cap Dom(S)
=
P(F)\cap Dom(SP)
$,
and
that
$$
S P
\up P(F)
=
S\up(P(F)\cap Dom(S)).
$$
Moreover
in the case in which
$PS\subseteq SP$
then
$$
SP
\up P(F):
P(F)\cap Dom(S)
\to
P(F).
$$
That 
is
$
SP
\up P(F)
$
is
a linear operator
in
the Banach space
$P(F)$.
Let
$
E:
\B_{Y}
\to
\Pr(G)
$
be
a spectral measure
in
$G$
on
$\B_{Y}$,
$\sigma\in\B_{Y}$
and
$Q$
a
possibly unbounded operator
in $G$
such that
$
E(\sigma)Q
\subseteq 
QE(\sigma)
$,
then
$$
Q_{\sigma}\up G_{\sigma}:
G_{\sigma}\cap Dom(Q)
\to
G_{\sigma}.
$$
In
particular
if
$R$
is
a 
possibly unbounded 
scalar type spectral operator
in 
$G$,
$
E
$
its
resolution of the identity
and
$f\in Bor(\sigma(R))$,
then
by 
statement
$(g)$
of
Theorem
$18.2.11$
of
\cite{ds},
we have 
that
for all
$
\sigma\in\B(\C)
$
$$
E(\sigma)f(R)
\subseteq
f(R) E(\sigma).
$$
Hence
for all
$
\sigma\in\B(\C)
$
\begin{equation}
\label{12022702}
\begin{cases}
R_{\sigma}\up G_{\sigma}
=
R_{\sigma}
\up
(G_{\sigma}\cap Dom(R))
=
R
\up
(G_{\sigma}\cap Dom(R))
\\
f(R)_{\sigma}\up G_{\sigma}
=
f(R)_{\sigma}
\up
(G_{\sigma}\cap Dom(f(R)))
=
f(R)
\up
(G_{\sigma}\cap Dom(f(R)))
\end{cases}
\end{equation}
are
linear operators
in
$G_{\sigma}$.
Finally
$E(\sigma(R))=\un$
implies
$
E(\sigma)
=
E(\sigma\cap\sigma(R))
$
for all
$\sigma\in\B(\C)$
so
by
\eqref{12022702}
\begin{equation}
\label{12022702bis}
\begin{cases}
R_{\sigma}\up G_{\sigma}
=
R_{\sigma\cap\sigma(R)}
\up G_{\sigma\cap\sigma(R)}
\\
f(R)_{\sigma}
\up G_{\sigma}
=
f(R)_{\sigma\cap\sigma(R)}
\up G_{\sigma\cap\sigma(R)}.
\end{cases}
\end{equation}
\end{definition}
\begin{lemma}
[Key Lemma]
\label{II31031834}
Let
$R$
be
a 
possibly unbounded 
scalar type spectral operator
in $G$,
$E$
its
resolution of the identity,
$\sigma(R)$ its spectrum 
and
$f\in Bor(\sigma(R))$. 
Then for all
$\sigma\in\B(\C)$
\begin{enumerate}
\item
$
R_{\sigma} 
\up G_{\sigma}
$
is a scalar type spectral operator 
in $G_{\sigma}$
whose resolution of the identity 
$\w{E}_{\sigma}$
is 
such that
for all
$\delta\in\B(\C)$
$$
\w{E}_{\sigma}(\delta)
=
E(\delta)
\up G_{\sigma}
\in
B(G_{\sigma}),
$$
\item
$$
f(R)_{\sigma} 
\up G_{\sigma}
=
f
\left(
R_{\sigma} 
\up 
G_{\sigma}
\right),
$$
\item for all
$g\in Bor(\sigma(R))$ 
such that
$g(\sigma\cap\sigma(R))$ 
is bounded, 
we have that
$$
g(R)E(\sigma)
=
\I{\C}{E}(\w{g}\cdot\chi_{\sigma})
\in 
B(G).
$$
\end{enumerate}
\end{lemma}
\begin{proof}
Let
$\sigma\in\B(\C)$.
By the fact 
that
$
E(\sigma\cap\delta)
=
E(\delta)
E(\sigma)
=
E(\sigma)
E(\delta)
$
for all
$\delta\in\B(\C)$
and
$
E(\sigma)\up G_{\sigma}
=
\un_{\sigma}
$
the unity 
operator on $G_{\sigma}$,
we
have
for all
$\delta\in\B(\C)$
\begin{equation}
\label{12152702}
\w{E}_{\sigma}(\delta)
=
E(\sigma\cap\delta)\up G_{\sigma}
\in
B(G_{\sigma}).
\end{equation}
In particular
$
\w{E}_{\sigma}:
\B(\C)
\to
B(G_{\sigma})
$,
moreover
$E$
is a countably
additive spectral measure
in $G$,
so
\begin{equation}
\label{11052053}
\w{E}_{\sigma}
\text{
is a countably additive spectral measure in 
$G_{\sigma}$.
}
\end{equation}
By 
Lemma $18.2.2.$
of
\cite{ds}
$\w{E}_{\sigma}$
is the resolution
of identity
of the 
spectral
operator
$R_{\sigma}\up G_{\sigma}$
so
by
Lemma $18.2.25.$
of
\cite{ds}
applied
to
$
R_{\sigma}\up G_{\sigma}
$
\begin{equation}
\label{10132802}
\textrm{supp}\,\w{E}_{\sigma}
=
\sigma(R_{\sigma}\up G_{\sigma}).
\end{equation}
Furthermore
by
\eqref{12022702bis}
and 
$(iii)$
of Definition
$18.2.1.$
of \cite{ds}
we have
$
\sigma(R_{\sigma}\up G_{\sigma})
\subseteq
\ov{\sigma\cap\sigma(R)}
$,
then
by
the equality
$
\ov{
\ov{\sigma}
\cap
\sigma(R)
}
=
\ov{\sigma}
\cap
\sigma(R)
$,
we
deduce 
\begin{equation}
\label{25052019}
\sigma(R_{\sigma}\up G_{\sigma})
\subseteq
\ov{\sigma}\cap\sigma(R)
\subseteq
\sigma(R).
\end{equation}
Hence
\eqref{10132802}
and
\eqref{25052019}
imply
that
the
operator
function
$
f(\w{E}_{\sigma})
$
is
well defined.
For all
$
x\in 
Dom(f(R)_{\sigma}\up G_{\sigma})
$
\begin{alignat*}{2}
(f(R)_{\sigma}\up G_{\sigma})x
&
=
f(R)x
&
\text{ by }
\eqref{12022702}
\\
&
=
\lim_{n\in\N}
\I{\C}{E^{x}}
(\w{f}\cdot\chi_{\overset{-1}{|f|}
(\delta_{n})})
&
\text{ by }
\eqref{II01041439},
\eqref{01041859}
\\
&
=
\lim_{n\in\N}
\I{\C}{\w{E}_{\sigma}^{x}}
(\w{f}\cdot\chi_{\overset{-1}{|f|}
(\delta_{n})})
&
\text{ by 
$x\in G_{\sigma}$,
\eqref{11052053}}
\\
&
=
\lim_{n\in\N}
\I{\C}{\w{E}_{\sigma}}
(\w{f}\cdot\chi_{\overset{-1}{|f|}
(\delta_{n})})
x
&
\text{ by }
\eqref{01041859}
\\
&
=
f(\w{E}_{\sigma})
x.
&
\text{ by }
\eqref{II01041439}
\end{alignat*}
So
$
f(R)_{\sigma}\up G_{\sigma}
\subseteq
f(\w{E}_{\sigma})
$.
For
all 
$
x
\in Dom(f(\w{E}_{\sigma}))
$
\begin{alignat*}{2}
f(\w{E}_{\sigma})
x
&
=
\lim_{n\in\N}
\I{\C}{\w{E}_{\sigma}^{x}}
(\w{f}\cdot\chi_{\overset{-1}{|f|}
(\delta_{n})})
&
\text{ by }
\eqref{II01041439},
\eqref{01041859}
\\
&
=
\lim_{n\in\N}
\I{\C}{E^{x}}
(\w{f}\cdot\chi_{\overset{-1}{|f|}
(\delta_{n})})
\\
&
=
\lim_{n\in\N}
\I{\C}{E}
(\w{f}\cdot\chi_{\overset{-1}{|f|}
(\delta_{n})})
x
&
\text{ by }
\eqref{01041859}
\\
&
=
(f(R)_{\sigma}\up G_{\sigma})x.
&
\text{ by }
\eqref{II01041439},
\eqref{12022702}
\end{alignat*}
So
$
f(\w{E}_{\sigma})
\subseteq
f(R)_{\sigma}\up G_{\sigma}
$,
then
\begin{equation}
\label{09282802}
f(R)_{\sigma}\up G_{\sigma}
=
f(\w{E}_{\sigma})
\end{equation}
Therefore
statement
$(1)$
follows
by setting
$
f=\imath
$,
while
statement
$(2)$
follows
by statement
$(1)$
and
\eqref{09282802}.
Let
$g\in Bor(\sigma(R))$ 
such that
$g(\sigma\cap\sigma(R))$ 
is bounded, 
then 
$$
(\exists\,n\in\N)(\forall m>n)
(
\sigma\cap\overset{-1}{|g|}(\delta_{m})
=
\sigma\cap\sigma(R)
).
$$
Next
$
E(\sigma(R))
=
\un
$,
so
$
E(\sigma)
=
E(\sigma)
E(\sigma(R))
=
E(\sigma(R)\cap\sigma)
$.
Since
$\I{\C}{E}$
is an algebra 
homomorphism, for all
$m\in\N$
\begin{alignat*}{1}
\I{\C}{E}(\w{g}\cdot\chi_{
\overset{-1}{|g|}(\delta_{m})
}
)
E(\sigma)
&
=
\I{\C}{E}(\w{g}\cdot\chi_{\overset{-1}{|g|}(\delta_{m})})
E(\sigma\cap\sigma(R))
\\
&
=
\I{\C}{E}(\w{g}\cdot\chi_{\overset{-1}{|g|}(\delta_{m})})
\I{\C}{E}(\chi_{\sigma\cap\sigma(R)})
\\
&
=
\I{\C}{E}
(
\w{g}\cdot\chi_{\overset{-1}{|g|}(\delta_{m})}
\cdot
\chi_{\sigma\cap\sigma(R)}
)
\\
&
=
\I{\C}{E}
(
\w{g}
\cdot
\chi_{\overset{-1}{|g|}(\delta_{m})\cap\sigma\cap\sigma(R)}
)
\\
&
=
\I{\C}{E}
(
\w{g}
\cdot
\chi_{\overset{-1}{|g|}(\delta_{m})\cap\sigma}
).
\end{alignat*}
This equality implies that
\begin{equation}
\label{II01041435}
(\exists\,n\in\N)(\forall m>n)
(\I{\C}{E}(\w{g}\cdot\chi_{\overset{-1}{|g|}(\delta_{m})})
E(\sigma)
=
\I{\C}{E}
(\w{g}
\cdot
\chi_{\sigma\cap\sigma(R)})
).
\end{equation}
Furthermore
\begin{alignat*}{1}
\I{\C}{E}
(\w{g}\cdot
\chi_{\sigma\cap\sigma(R)})
&
=
\I{\C}{E}
(\w{g}\chi_{\sigma}\chi_{\sigma(R)})\\
&
=
\I{\C}{E}
(\w{g}\chi_{\sigma})
\I{\C}{E}
(\chi_{\sigma(R)})\\
&
=
\I{\C}{E}
(\w{g}\chi_{\sigma})
E(\sigma(R))\\
&
=
\I{\C}{E}
(\w{g}\chi_{\sigma}).
\end{alignat*}
Therefore by 
\eqref{II01041435}
\begin{equation}
\label{II01041438}
(\exists\,n\in\N)(\forall m>n)
(
\I{\C}{E}(\w{g}\cdot\chi_{\overset{-1}{|g|}(\delta_{m})})
E(\sigma)
=
\I{\C}{E}
(\w{g}
\cdot
\chi_{\sigma})
).
\end{equation}
Moreover
by definition
in \eqref{II01041439}
we have for all
$x\in Dom(g(R))$
that
$$
g(R)x
\doteqdot
\lim_{n\to\infty}
\I{\C}{E}
(\w{g}\cdot\chi_{\overset{-1}{|g|}
(\delta_{n})})x
$$
and $Dom(g(R))$ 
is the set of $x\in G$ 
such that such a limit exists;
thus
by 
\eqref{II01041438}
we can conclude that
$E(\sigma)G\subseteq Dom(g(R))$
and
$
g(R)E(\sigma)
=
\I{\C}{E}(\w{g}\cdot\chi_{\sigma})
\in 
B(G)
$,
which is
statement $(3)$.
\end{proof}
\begin{corollary}
\label{22061735}
Let
$R$
be
a 
possibly unbounded 
scalar type spectral operator
in $G$,
and
$f\in Bor(\sigma(R))$.
Then for all
$\sigma\in\B(\C)$
$$
f(R) E(\sigma)
=
f
\left(
R_{\sigma} 
\up G_{\sigma}
\right)
E(\sigma).
$$
Moreover
if
$f(\sigma\cap\sigma(R))$
is
bounded
then
$$
f
\left(
R_{\sigma} 
\up G_{\sigma}
\right)
E(\sigma)
\in
B(G).
$$
\end{corollary}
\begin{proof}
Let $y\in Dom(f(R) E(\sigma))$
then
$E(\sigma)y\in G_{\sigma}\cap Dom(f(R))$
hence by 
\eqref{12022702},
Lemma \ref{II31031834}
\begin{equation*}
f(R)E(\sigma)y
=
(f(R)_{\sigma}\up G_{\sigma})
E(\sigma)y
=
f(R_{\sigma}\up G_{\sigma})
E(\sigma)y.
\end{equation*}
So
$
f(R)E(\sigma)
\subseteq
f(R_{\sigma}\up G_{\sigma})
E(\sigma)
$.
Next
let 
$y\in 
Dom(f(R_{\sigma}\up G_{\sigma})E(\sigma))
$,
then 
$E(\sigma)y\in 
Dom(f(R_{\sigma}\up G_{\sigma}))
$,
hence
by
Lemma \ref{II31031834}
and
\eqref{12022702}
\begin{equation*}
f(R_{\sigma}\up G_{\sigma})
(E(\sigma)y)
=
f(R)E(\sigma)E(\sigma)y
=
f(R)E(\sigma)y.
\end{equation*}
So
$
f(R_{\sigma}\up G_{\sigma})
E(\sigma)
\subseteq
f(R)E(\sigma)
$.
Thus
we obtain 
statement
$(1)$.
Statement
$(2)$
follows
by 
statement
$(1)$
and
statement $(3)$
of
Lemma
\ref{II31031834}.
\end{proof}
\section{
Extension
theorem
for
strong operator 
integral
equalities
}
\begin{notations}
\label{13051406}
\begin{normalfont}
Let $X$ be a locally compact space
and $\mu$ a 
measure on 
$X$ in 
the sense
of the Bourbaki text
\cite{IntBourb}
see $III.7$,
Definition $2$,
that is a continuous linear $\C-$functional on
the $\C-$locally convex space 
$H(X)$
of all compactly supported 
complex 
continuous functions
on $X$, 
with the direct limit topology
(or inductive limit)
of the
spaces
$
H(X;K)
$
with
$K$
running
in
the class of all compact subsets of $X$,
where
$H(X;K)$
is 
the space
of
all
complex
continuous functions 
$f:X\to\C$
such
that
$
\textrm{supp}(f)
\doteqdot
\ov{
\{x\in X\mid f(x)\neq 0\}
}
\subseteq 
K
$
with
the
norm
topology of uniform 
convergence
\footnote{
$H(X;K)$
is isometric
to the
Banach
space
of 
all 
continuous
maps
$g:K\to\C$
equal
to 
$0$
on
$\partial K$,
with
the
norm
topology of uniform convergence
}.
In the work 
any
measure 
$\mu$ 
on $X$
in the sense of \cite{IntBourb}
will be called
complex Radon measure on $X$.
For the definition
of 
\emph{
$\mu-$integrable 
functions 
defined on $X$
and
with
values 
in a 
$\C-$Banach space
$G$
}
see 
$IV.23.$
Definition
$2$
of
\cite{IntBourb},
while
the
\emph{
integral with respect to 
$\mu$
of a
$\mu-$integrable 
function
$f:X\to G$
},
which will be denoted with
$
\int 
f(x)\,d\,\mu(x)
\in
G
$,
is
defined
in
Definition
$1$,
$III.33$
and
Definition
$1$,
$IV.33$
of
\cite{IntBourb}.
For
the definition
of
the total
variation
$|\mu|$,
and
definition
and
properties
of
the upper integral
$
\int^{*}g\,d|\mu|(x)
$
see
Ch. $3-4$
of
\cite{IntBourb}.
We
denote by
$Comp(X)$
the class of the compact 
subsets of $X$
and 
by
$
\F{1}{}(X;\mu)
$
the
seminormed space,
with seminorm
$
\|\cdot\|_{\F{1}{}(X;\mu)}
$,
of all 
maps
$
F:X\to\C
$
such that
$$
\|
F
\|_{\F{1}{}(X;\mu)}
\doteqdot
\int^{*} 
|F(x)|
\,d\,
|\mu|(x)
<\infty.
$$
In this section
it
will be assumed,
unless 
the contrary is stated,
that $X$ is a 
locally compact space 
and
$\mu$ is
a complex Radon
measure 
over 
$X$.
Let 
$B\subseteq X$
be
a
$\mu-$measurable set,
then
by
$\mu-a.e.(B)$
we 
mean
``
almost everywhere
in
$B$
with respect to the 
measure
$\mu$
''.
Let
$f:X\to B(G)$
be a map
$\mu-$integrable 
in the normed space
$B(G)$
(Definition
$2$
Ch. $IV$, \S $3$, $n^{\circ} 4$
of
\cite{IntBourb})
then
we convene
to denote
with the symbol
$$
\oint 
f(x)\,d\,\mu(x)
\in
B(G)
$$
its
integral
in $B(G)$
(Definition
$1$
Ch. $IV$, \S $4$, $n^{\circ} 1$
of
\cite{IntBourb}),
which is uniquely
determined
by the following
property for all
$\phi\in B(G)^{*}$
$$
\phi(\oint 
f(x)\,d\,\mu(x))
=
\int
\phi(f(x))\,d\,\mu(x).
$$
For any 
scalar
type 
spectral operator
$S$
in a complex Banach space $G$
and 
for any
Borelian map
$
f:
U\supseteq
\sigma(S)\to\C
$
we assume that
$f(S)$
is the closed operator
defined in 
\eqref{II01041439}
and recall that
we denote
by
$\w{f}$
the $\ze-$extension
of $f$ to $\C$,
see Definition \ref{12470108}.
\end{normalfont}
\end{notations}
\begin{definition}
[$E-$sequence]
\label{17331001}
Let 
$E:\B_{Y}\to Pr(G)$
be
a 
spectral measure 
in 
$G$
on
$\B_{Y}$
then we say that
$\{\sigma_{n}\}_{n\in\N}$
is an
\emph{$E-$sequence}
if
there
exists
an
$S\in\B_{Y}$
such that
$E(S)=\un$
and
\begin{itemize}
\item
$(\forall n\in\N)(\sigma_{n}\in\B_{Y})$;
\item
$
(\forall n,m\in\N)
(n>m
\Rightarrow
\sigma_{n}
\supseteq
\sigma_{m})
$;
\item
$
S
\subseteq
\bigcup_{n\in\N}
\sigma_{n}
$.
\end{itemize}
\end{definition}
\begin{proposition}
\label{15422702}
Let $E:\B_{Y}\to Pr(G)$
be
a 
countably additive
spectral measure 
in 
$G$
on
a $\sigma-$
field
$\B_{Y}$,
and
$\{\sigma_{n}\}_{n\in\N}$
an
$E-$sequence.
Then
\begin{equation}
\label{II03041401}
\lim_{n\to\infty}
E(\sigma_{n})
=
\un
\quad
\text{ in strong operator topology.}
\end{equation}
\end{proposition}
\begin{proof}
Let 
$
S
\in
\B_{Y}
$
of 
which
in Definition
\ref{17331001}
associated
to the
$E-$sequence
$\{\sigma_{n}\}_{n\in\N}$.
So
$E(S)=\un$
and
$E$ is an
order-preserving
map,
then
$
E(\bigcup_{n\in\N}\sigma_{n})
\geq
E(S)
=
\un
$.
Since
$\un$ is 
a maximal element in $\lr{E(\B_{Y})}{\geq}$
$$
E(\bigcup_{n\in\N}\sigma_{n})
=
\un.
$$
Let us define
$
\eta_{1}
\doteqdot
\sigma_{1}
$,
and
for all
$n\geq 2$,
$\eta_{n}
\doteqdot
\sigma_{n}\cap\complement\sigma_{n-1}$,
so 
for all
$n\in\N$,
$\sigma_{n}
=
\bigcup_{k=1}^{n}
\eta_{k}$,
and
for all
$n\neq m\in\N$,
$\eta_{n}\cap\eta_{m}=\emptyset$,
finally
$
\bigcup_{n\in\N}\eta_{n}
=
\bigcup_{n\in\N}
\left(
\bigcup_{k=1}^{n}
\eta_{k}
\right)
=
\bigcup_{n\in\N}
\sigma_{n}
$.
Therefore
by the countable additivity 
of $E$ 
with respect to
the strong operator topology
\begin{alignat*}{1}
E(\bigcup_{n\in\N}\sigma_{n})
&
=
E(\bigcup_{n\in\N}\eta_{n})
=
\sum_{n=1}^{\infty}
E(\eta_{n})
\\
&
=
\lim_{n\to\infty}
\sum_{k=1}^{n}
E(\eta_{k})
=
\lim_{n\to\infty}
E(\bigcup_{k=1}^{n}\eta_{k})
\\
&
=
\lim_{n\to\infty}
E(\sigma_{n}).
\end{alignat*}
Here all
limits are with respect to 
the strong operator topology,
hence 
the statement.
\end{proof}
\begin{definition}
[Integration 
in the Strong Operator Topology]
\label{13051513}
Let $G_{1},G_{2}$ be two complex Banach spaces,
and 
$
f:
X
\to B(G_{1},G_{2})
$.
Then we say that
\emph{
$f$
is 
$\mu-$
integrable
with respect
to the
strong operator topology
}
if
\begin{enumerate}
\item
for all
$
v\in G_{1}
$
the map
$
X
\ni x
\mapsto 
f(x)v\in G_{2}
$
is
$\mu-$integrable;
\item
if we set
$$
F:
G_{1}
\ni
v
\mapsto
\int f(x)(v)\,d\,\mu(x)
\in
G_{2}
$$
then
$F\in B(G_{1},G_{2})$.
\end{enumerate}
In such 
a case
we 
set
$
\int f(x)\,d\,\mu(x)
\doteqdot
F
$,
in other words
the 
integral
$
\int f(x)\,d\,\mu(x)
$
of $f$ 
with respect to the 
measure
$\mu$
and
the
strong operator topology
is a bounded linear operator from 
$G_{1}$ to $G_{2}$ 
such that
for all $v\in G_{1}$
$$
\left(
\int f(x)\,d\,\mu(x)
\right)(v)
=
\int f(x)(v)\,d\,\mu(x).
$$
\end{definition}
We shall need the following
version of the Minkowski
inequality
\begin{proposition}
\label{13051512}
Let 
$G_{1},G_{2}$ be two complex Banach spaces,
and 
a map
$
f:
X
\to B(G_{1},G_{2})
$
such that
\begin{enumerate}
\item
$
(\forall v\in G_{1})
(\forall \phi\in G_{2}^{*})
$
the complex map
$
X
\ni x
\mapsto\phi(f(x)v)\in\C
$
is $\mu-$measurable;
\item
for all
$v\in G_{1},
K\in Comp(X)$
there is
$H\subset G_{2}$
such that $H$ is countable
and
$f(x)v\in\ov{H}$\,
$\mu-a.e.(K)$;
\item
$
(
X\ni x
\mapsto
\|f(x)\|_{B(G_{1},G_{2})}
)
\in
\F{1}{}(X;\mu)
$,
\end{enumerate}
Then
$f$
is
$\mu-$integrable
with respect
to the strong operator topology
and
we have
$$
\left
\|
\int
f(x)
\,d\,
\mu(x)
\right
\|_{B(G_{1},G_{2})}
\leq
\int^{*}
\|
f(x)
\|_{B(G_{1},G_{2})}
\,d\,
|\mu|(x).
$$
\end{proposition}
\begin{proof}
By
hypothesis $(3)$
we have for all
$v\in G_{1}$
\begin{equation}
\label{II10051823}
\int^{*}
\|
f(x)
v
\|_{G_{2}}
\,d\,
|\mu|(x)
\leq
\|v\|_{G_{1}}
\int^{*}
\|
f(x)
\|_{B(G_{1},G_{2})}
\,d\,
|\mu|(x)
<\infty.
\end{equation}
By
hypothesis
$(1-2)$
and
Corollary
$1$,
$IV.70$
of
\cite{IntBourb},
we 
have for all
$
v\in G_{1}
$
that
the map
$
X\mapsto f(x)v\in G_{2}
$
is
$\mu-$measurable.
Therefore
by
\eqref{II10051823}
and
by
Theorem $5$,
$IV.71$
of
\cite{IntBourb}
we
deduce for all
$
v\in G_{1}
$
that
$
X\mapsto f(x)v\in G_{2}
$
is
$\mu-$integrable.
So in particular
by 
Definition
$1$,
$IV.33$
of
\cite{IntBourb}
for all
$v\in G_{1}$
there is
$\int
f(x)v
\,d\,
\mu(x)
\in
G_{2}$
while
by
Proposition $2$,
$IV.35$
of
\cite{IntBourb}
and the 
\eqref{II10051823}
we obtain for all
$v\in G_{1}$
\begin{equation*}
\left
\|
\int
f(x)v
\,d\,
\mu(x)
\right
\|_{G_{2}}
\leq
\|v\|_{G_{1}}
\int^{*}
\|
f(x)
\|_{B(G_{1},G_{2})}
\,d\,
|\mu|(x)
\end{equation*}
Hence
the statement
follows.
\end{proof}
\begin{remark}
\label{13051626}
As it follows
by the above proof
Proposition
\ref{13051512}
is also valid 
if we replace the hypotheses
$(1-2)$
with
the
following one
\begin{equation}
\forall v\in G_{1}
\quad
\text{the map }\,
X
\ni x
\mapsto
f(x)v
\in
G_{2}
\quad
\text{is
$\mu-$measurable.}
\tag*{(1')}
\end{equation}
\end{remark}
\begin{lemma}
\label{II04041224}
Let
$X, Y, Z$
be
three normed spaces 
over the same field $\K\in\{\R,\C\}$,
$R: Dom(R)\subseteq Y\to Z$ 
a possibly 
unbounded closed 
linear
operator 
and $A\in B(X,Y)$.
Then
$R\circ A:\Df\to Z $
is a closed operator,
where
$\Df\doteqdot Dom(R\circ A)$
\end{lemma}
\begin{proof}
Let
$
\{x_{n}\}_{n\in\N}
\subset 
\Df
\doteqdot
\{
x\in X
\mid A(x)\in Dom(R)
\}
$,
and
$(x,z)\in X\times Z$
such that
$
x
=
\lim_{n\to\infty}
x_{n}
$,
and
$
z
=
\lim_{n\to\infty}
R\circ A
(x_{n})
$.
$A$ 
being 
continuous we have
$
A(x)
=
\lim_{n\to\infty}
A(x_{n})
$,
but 
$
z
=
\lim_{n\to\infty}
R(Ax_{n})
$,
and
$R$ is closed, so 
$
z
=
R(A(x))
\doteqdot
R\circ A(x)
$,
hence
$
(x,z)\in Graph(R\circ A)
$,
which
is just the statement.
\end{proof}
\begin{lemma}
\label{II04041227}
Let 
$X$
be a
normed space 
and
$Y$
a Banach space
over the same field $\K\in\{\R,\C\}$,
finally
$
U:\Df\subseteq X\to Y
$
be
a linear
operator. 
If $U$ is continuous and closed
then
$\Df$ is closed.
\end{lemma}
\begin{proof}
Let
$\{x_{n}\}_{n\in\N}\subset\Df$ 
and
$x\in X$
such that 
$x=\lim_{n\to\infty}x_{n}$.
So
by the continuity of $U$
we have for all
$n,m\in\N$
that
$
\|U(x_{n})-U(x_{m})\|
=
\|U(x_{n}-x_{m})\|
\leq
\|U\|
\|x_{n}-x_{m}\|
$,
hence
$
\lim_{(n,m)\in\N^{2}}
\|U(x_{n})-U(x_{m})\|
=
0
$,
thus
$Y$ 
being 
a Banach space we have that
there is
$y\in Y$
such that
$y=\lim_{n\to\infty}U(x_{n})$.
But $U$ 
is closed,
therefore
$y=U(x)$, 
so
$x\in\Df$,
which
is the statement.
\end{proof}
\begin{theorem}
\label{10541501}
Let
$R$
be
a 
possibly unbounded 
scalar type spectral operator
in $G$,
$\sigma(R)$
its spectrum
and
$E$ its resolution of the identity.
Let 
the map
$
X
\ni 
x
\mapsto
f_{x}
\in
Bor(\sigma(R))
$
be
such that
for all $x\in X$,
$
\w{f_{x}}
\in
\Lf{E}{\infty}(\sigma(R))
$
where 
$
X\ni x\mapsto f_{x}(R)\in B(G)
$
is
$\mu-$integrable
with respect to the 
strong operator topology.
 
Then
\begin{enumerate}
\item
for all
$\sigma\in\B(\C)$
the
map
$
X
\ni
x
\mapsto
f_{x}(R_{\sigma}\up G_{\sigma})
\in
B(G_{\sigma})
$
is
$\mu-$integrable
with
respect
to
the strong operator topology
and
$$
\left
\|
\int
\,
f_{x}(R_{\sigma}\up G_{\sigma})
\,
d\,
\mu(x)
\right
\|_{B(G_{\sigma})}
\leq
\left
\|
\int
\,
f_{x}(R)
\,
d\,
\mu(x)
\right
\|_{B(G)}.
$$
\item
If
$
g,
h\in
Bor(\sigma(R))
$,
$\{\sigma_{n}\}_{n\in\N}$
is
an
$E-$sequence,
and for all
$n\in\N$
\begin{equation}
\label{20101201}
g(R_{\sigma_{n}}\up G_{\sigma_{n}})
\int\,
f_{x}(R_{\sigma_{n}}\up G_{\sigma_{n}})
\,d\,
\mu(x)
\subseteq
h(R_{\sigma_{n}}\up G_{\sigma_{n}}).
\end{equation}
then
\begin{equation}
\label{23241501}
g(R)
\int\,
f_{x}(R)
\,d\,
\mu(x)
\up
\Theta
=
h(R)
\up
\Theta,
\end{equation}
\end{enumerate}
where
$
\Theta
\doteqdot
Dom
\left(
g(R)
\int\,
f_{x}(R)
\,d\,
\mu(x)
\right)
\cap
Dom(h(R))
$
and
all the integrals 
are with respect to the
strong operator topologies.
\end{theorem}
Notice that
$g(R)$ is 
possibly an
\textbf{
unbounded}
operator
in
$G$.
\begin{proof}
Let
$\sigma\in\B(\C)$
then
by
\eqref{25052019}
$$
(\forall\sigma\in\B(\C))
(\sigma(R_{\sigma}\up G_{\sigma})
\subseteq
\ov{\sigma}\cap\sigma(R)
\subseteq
\sigma(R)).
$$
which
implies
that
all the following
operator
functions
$
g(R_{\sigma}\up G_{\sigma})
$,
$
h(R_{\sigma}\up G_{\sigma})
$
and for all
$x\in X$
the
$
f_{x}(R_{\sigma}\up G_{\sigma})
$,
are 
well defined.
By 
the fact
that
$
\{
\delta\in\B(\C)
\mid
E(\delta)=\un
\}
\subseteq
\{
\delta\in\B(\C)
\mid
\w{E}_{\sigma}(\delta)=\un_{\sigma}
\}
$
which follows
by statement
$(1)$
of 
Lemma \ref{II31031834},
we deduce for all
$x\in X$
$$
\|
\w{f}_{x}
\|_{\infty}^{\w{E}_{\sigma}}
\leq
\|
\w{f}_{x}
\|_{\infty}^{E}
=
\|
\w{f}_{x}
\chi_{\sigma(R)}
\|_{\infty}^{E}
<
\infty,
$$
where the last
equality
came
by 
$
\w{f}_{x}
\chi_{\sigma(R)}
=
\w{f}_{x}
$,
while
the boundedness
by the
hypothesis
$
\w{f_{x}}
\in
\Lf{E}{\infty}(\sigma(R))
$.
Thus
$
\w{f_{x}}
\in
\Lf{\w{E}_{\sigma}}{\infty}
(\C)
$
hence
by 
statement 
$(c)$
of
Theorem
$18.2.11.$
of
\cite{ds}
applied to
the scalar type
spectral operator
$
R_{\sigma}\up G_{\sigma}
$
\begin{equation}
\label{15051542}
(\forall\sigma\in\B(\C))
(f_{x}(R_{\sigma}\up G_{\sigma})
\in
B(G_{\sigma})).
\end{equation}
A more direct
way for 
obtaining
\eqref{15051542}
is to use 
statement
$(2)$
of
Lemma
\ref{II31031834}
and
the fact that
$
\w{f_{x}}
\in
\Lf{E}{\infty}(\sigma(R))
$
implies
$
f_{x}(R)
\in
B(G)
$.
For all
$\sigma\in\B(\C)$ 
we claim
that
$
X\ni x
\mapsto
f_{x}(R_{\sigma}\up G_{\sigma})
\in
B(G_{\sigma})
$
is
$\mu$-integrable
with respect
to the
strong operator
topology.
By
Lemma
\ref{II31031834}
we have
for all
$\sigma\in\B(\C)$
and for all
$v\in G_{\sigma}$
\begin{equation}
\label{15051418}
\int^{*}
\|
f_{x}(R_{\sigma}\up G_{\sigma})
v
\|_{G_{\sigma}}
\,d\,
|\mu|(x)
=
\int^{*}
\|
f_{x}(R)
v
\|_{G}
\,d\,
|\mu|(x)
<\infty.
\end{equation}
Here
the
boundedness
comes
by
Theorem
$5$,
$IV.71$
of 
\cite{IntBourb}
applied
to the 
$\mu-$integrable
map
$X\ni x\mapsto f_{x}(R)v\in G$.
By 
Corollary $1$, $IV.70$ 
and
Theorem $5$, $IV.71$ 
of
\cite{IntBourb}
applied,
for any
$v\in G$,
to the $\mu-$integrable map
$X\ni x\mapsto f_{x}(R)v\in G$,
we have for all
$v\in G,K\in Comp(X)$
there is
$H^{v}\subseteq G$
countable 
such that
$(f_{x}(R)v\in\ov{H^{v}},\,\mu-a.e.(K))$.
But
by 
statement $(g)$
of
Theorem $18.2.11.$
of
\cite{ds}
and
$f_{x}(R)\in B(G)$,
we have 
for all
$\sigma\in\B(\C)$,
$
[f_{x}(R),E(\sigma)]=\ze
$,
hence
by
the
previous
equation
and
by
the fact that
$E(\sigma)\in B(G)$,
so it is
continuous,
we obtain
for all
$\sigma\in\B(\C),
v\in G,
K\in Comp(X)$
$$
(\exists\,H^{v}\subseteq G\text{ countable })
(
f_{x}(R)
E(\sigma)
v
=
E(\sigma)
f_{x}(R)
v
\in
\ov{
H_{\sigma}^{v}
},
\,
\mu-a.e.(K)).
$$
Here
$
H_{\sigma}^{v}
\doteqdot
E(\sigma)H^{v}
$.
Therefore
by 
Lemma
\ref{II31031834}
we state that
for all
$
\sigma\in\B(\C),
v\in G_{\sigma},
K\in Comp(X)
$
\begin{equation}
\label{15051543}
(\exists\,H_{\sigma}^{v}
\subset G_{\sigma}\text{ countable })
(
f_{x}(R_{\sigma}\up G_{\sigma})
v
\in
\ov{
H_{\sigma}^{v}
}
\subseteq
G_{\sigma},
\,
\mu-a.e.(K)).
\end{equation}
That
$\ov{
H_{\sigma}^{v}
}
\subseteq
G_{\sigma}
$
follows by the fact that $G_{\sigma}$ is closed in $G$.
Therefore
we
can consider
the closure
$
\ov{
H_{\sigma}^{v}
}
$
as
the 
closure in the Banach
space
$
G_{\sigma}
$.
By
the Hahn-Banach Theorem, 
see Corollary $3$,
$II.23$ of the \cite{BourTVS}, 
for all
$\sigma\in\B(\C)$
\begin{equation}
\label{15051506}
\{\phi\up G_{\sigma}\mid\phi\in G^{*}\}
=
(G_{\sigma})^{*}.
\end{equation}
Moreover
by
Corollary $1$, $IV.70$ 
and
Theorem $5$, $IV.71$ 
of
\cite{IntBourb}
applied,
for any
$v\in G$,
to the $\mu-$integrable map
$
X\ni x\mapsto f_{x}(R)E(\sigma)v\in G
$,
we have for all
$\phi\in G^{*}$
$$
X\ni x\mapsto 
\phi(f_{x}(R)E(\sigma)v)
\in\C
\text{ is $\mu-$measurable.}
$$
Thus
by 
Lemma
\ref{II31031834}
we have for all
$
\sigma\in\B(\C),
v\in G_{\sigma},
\phi\in G^{*}
$
$$
X\ni x\mapsto 
\phi(f_{x}(R_{\sigma}\up G_{\sigma})v)
\in\C
\text{ is $\mu-$measurable.}
$$
Hence 
by \eqref{15051506}
we can state for all
$
\sigma\in\B(\C),
v\in G_{\sigma}
$
that
\begin{equation}
\label{15051507}
(\forall\phi_{\sigma}\in (G_{\sigma})^{*})
(
X\ni x\mapsto 
\phi_{\sigma}
(f_{x}(R_{\sigma}\up G_{\sigma})v)
\in\C
\text{ is $\mu-$measurable.})
\end{equation}
Now
by collecting
\eqref{15051507},
\eqref{15051418}
and
\eqref{15051543}, 
where the closure
$
\ov{
H_{\sigma}^{v}
}
$
is to be intended
how
the 
closure in the Banach
space
$
G_{\sigma}
$,
we can apply 
Corollary $1$, $IV.70$ 
and
Theorem $5$, $IV.71$ 
of
\cite{IntBourb}
to the map
$
X\ni
x
\mapsto
f_{x}(R_{\sigma}\up G_{\sigma})v
\in
G_{\sigma}
$,
in order
to state 
that
\begin{equation}
\label{15051515}
(\forall\sigma\in\B(\C))
(\forall v\in G_{\sigma})
(X
\ni
x
\mapsto
f_{x}(R_{\sigma}\up G_{\sigma})v
\in
G_{\sigma}
\text{ is $\mu-$integrable.})
\end{equation}
This means in particular that 
there
exists 
its integral, so for all
$
\sigma\in\B(\C),
v\in G_{\sigma}
$
\begin{alignat}{2}
\label{15051536}
\left\|
\int
f_{x}(R_{\sigma}\up G_{\sigma})
v
\,d\,\mu(x)
\right\|_{G_{\sigma}}
&
=
\left\|
\int
f_{x}(R)
v
\,d\,\mu(x)
\right\|_{G}
&
\text{ by 
Lemma
\ref{II31031834}
}
\notag\\
&
\leq
\left\|
\int
f_{x}(R)
\,d\,\mu(x)
\right\|_{B(G)}
\|
v
\|_{G_{\sigma}}.
\end{alignat}
Here the inequality
follows
by
the hypothesis that 
$X\ni x\mapsto f_{x}(R)\in B(G)$
is
$\mu-$integrable in 
the strong operator topology.
Therefore
by 
Definition
\ref{13051513}
and
\eqref{15051542},
\eqref{15051515},
\eqref{15051536}
we can conclude that
\begin{equation}
\label{15051539}
\begin{cases}
(\forall\sigma\in\B(\C))
(X
\ni
x
\mapsto
f_{x}(R_{\sigma}\up G_{\sigma})
\in
B(G_{\sigma})
\text{ 
is $\mu-$integr.
in strong operator topology})\\
\left\|
\int
f_{x}(R_{\sigma}\up G_{\sigma})
\,d\,\mu(x)
\right\|_{B(G_{\sigma})}
\leq
\left\|
\int
f_{x}(R)
\,d\,\mu(x)
\right\|_{B(G)}.
\end{cases}
\end{equation}
Which
is the claim 
we wanted to show, 
then 
statement
$(1)$
follows.
Statement $(1)$
proves that 
the assumption
\eqref{20101201}
is well set,
so we are able to 
start the
proof
of the statement $(2)$.
For all
$y\in\Theta$
\begin{alignat*}{2}
&
=g(R)
\int 
f_{x}(R)
\,d\,\mu(x)
y
\\
&
=
\lim_{n\in\N}
E(\sigma_{n})g(R)\int f_{x}(R)\,d\,\mu(x)y
&
\text{
by
\eqref{II03041401}
}
\\
&
=
\lim_{n\in\N}
g(R)
E(\sigma_{n})
\int f_{x}(R)\,d\,\mu(x)y
&
\text{ by 
$(g)$ of
Theorem $18.2.11$ of \cite{ds}
}
\\
&
=
\lim_{n\in\N}
g(R)
E(\sigma_{n})
\int f_{x}(R)y\,d\,\mu(x)
&
\text{ by Definition \ref{13051513}}
\\
&
=
\lim_{n\in\N}
g(R)
E(\sigma_{n})
\int 
E(\sigma_{n})
f_{x}(R)
y
\,d\,\mu(x)
&
\text{
by 
Theorem
$1$,
$IV.35$
of
\cite{IntBourb}
}
\\
&
=
\lim_{n\in\N}
g(R)
E(\sigma_{n})
\int 
f_{x}(R)
E(\sigma_{n})
y
\,d\,\mu(x)
&
\text{
by 
$(g)$ of
Theorem $18.2.11$ of \cite{ds}
}
\\
&
=
\lim_{n\in\N}
g(R_{\sigma_{n}}\up G_{\sigma_{n}})
\int 
f_{x}(R_{\sigma_{n}}\up G_{\sigma_{n}})
E(\sigma_{n})
y
\,d\,\mu(x)
&
\text{
by 
\textbf{Lemma \ref{II31031834}}
}
\\
&
=
\lim_{n\in\N}
g(R_{\sigma_{n}}\up G_{\sigma_{n}})
\int 
f_{x}(R_{\sigma_{n}}\up G_{\sigma_{n}})
\,d\,\mu(x)
E(\sigma_{n})
y
&
\text{
by
statement
$(1)$ 
and
Definition
\ref{13051513}
}
\\
&
=
\lim_{n\in\N}
h(R_{\sigma_{n}}\up G_{\sigma_{n}})
E(\sigma_{n})
y
&
\text{
by
hypothesis
\eqref{20101201}
}
\\
&
=
\lim_{n\in\N}
h(R)
E(\sigma_{n})
y
&
\text{
by
Lemma \ref{II31031834}
}
\end{alignat*}
\begin{alignat*}{2}
&
=
\lim_{n\in\N}
E(\sigma_{n})
h(R)
y
&
\text{
by 
$(g)$ 
of
Theorem $18.2.11$ of \cite{ds}
}
\\
&
\,\,
h(R)
y
&
\text{
by 
\eqref{II03041401}.
}
\end{alignat*}
Therefore
$$
g(R)
\int\,
f_{x}(R)
\,d\,
\mu(x)
\up
\Theta
=
h(R)
\up
\Theta.
$$
\end{proof}
\begin{theorem}
[
\textbf{
Strong
Extension Theorem
}
]
\label{13051634}
Let
$X$
be a locally compact
space,
$\mu$
a complex Radon measure
on $X$,
$R$
be
a 
possibly unbounded 
scalar type spectral operator
in $G$,
$\sigma(R)$
its spectrum
and
$E$ its resolution of the identity.
Let 
the map
$
X
\ni 
x
\mapsto
f_{x}
\in
Bor(\sigma(R))
$
be
such that
for all 
$x\in X$,
$
\w{f_{x}}
\in
\Lf{E}{\infty}(\sigma(R))
$,
where 
the map
$
X\ni x\mapsto f_{x}(R)\in B(G)
$
be
$\mu-$integrable
with respect to the 
strong operator topology.
Finally
let
$
g,
h\in
Bor(\sigma(R))
$
and
$
\w{h}
\in
\Lf{E}{\infty}(\sigma(R))
$.
 
If
$\{\sigma_{n}\}_{n\in\N}$
is
an
$E-$sequence
and for all
$n\in\N$
\begin{equation}
\label{20081201}
g(R_{\sigma_{n}}\up G_{\sigma_{n}})
\int\,
f_{x}(R_{\sigma_{n}}\up G_{\sigma_{n}})
\,d\,
\mu(x)
\subseteq
h(R_{\sigma_{n}}\up G_{\sigma_{n}})
\end{equation}
then
$
h(R)
\in
B(G)
$
and
\begin{equation*}
g(R)
\int\,
f_{x}(R)
\,d\,
\mu(x)
=
h(R).
\end{equation*}
Here
all the integrals 
are with respect to the
strong operator topologies.
\end{theorem}
Notice that
$g(R)$ is 
possibly an
\textbf{
unbounded}
operator
on
$G$.
\begin{proof}
$h(R)\in B(G)$
by Theorem
$18.2.11.$ of \cite{ds}
and the hypothesis 
$
\w{h}
\in
\Lf{E}{\infty}(\sigma(R))
$,
so
by 
\eqref{23241501}
\begin{equation}
\label{15051936}
g(R)
\int\,
f_{x}(R)
\,d\,
\mu(x)
\subseteq
h(R).
\end{equation}
Let us set
\begin{equation}
\label{17502703ST}
(\forall n\in\N)
(\delta_{n}\doteqdot\overset{-1}{|g|}([0,n])).
\end{equation}
We claim that
\begin{equation}
\label{17542703ST}
\begin{cases}
\bigcup_{n\in\N}\delta_{n}
=\sigma(R)\\
n\geq m\Rightarrow\delta_{n}
\supseteq\delta_{m}\\
(\forall n\in\N)(g(\delta_{n})
\text{ is bounded. })\\
\end{cases}
\end{equation}
In addition
being
$|g|\in Bor(\sigma(R))$
we have
$\delta_{n}\in\B(\C)$
for all $n\in\N$,
so
$\{\delta_{n}\}_{n\in\N}$
is an $E-$sequence,
hence by \eqref{II03041401}
\begin{equation}
\label{20271201}
\lim_{n\in\N}E(\delta_{n})=\un
\end{equation}
with respect to the strong operator topology
on $B(G)$.
Indeed
the first equality 
follows since
$\bigcup_{n\in\N}
\delta_{n}
\doteq
\bigcup_{n\in\N}\overset{-1}{|g|}([0,n])
=
\overset{-1}{|g|}
\left(\bigcup_{n\in\N}[0,n]\right)
=
\overset{-1}{|g|}(\R^{+})
=
Dom(g)\doteqdot\sigma(R)
$,
the second
by the fact that
$\overset{-1}{|g|}$
preserves the inclusion,
the third since
$|g|(\delta_{n})\subseteq [0,n]$.
Hence our claim.
By the third statement
of \eqref{17542703ST},
$\delta_{n}\in\B(\C)$ 
and statement
$3$ of Lemma 
\ref{II31031834}
\begin{equation}
\label{19292703ST}
(\forall n\in\N)
(E(\delta_{n})G
\subseteq
Dom(g(R))).
\end{equation}
$
f_{x}(R)
E(\delta_{n})
=
E(\delta_{n})
f_{x}(R)
$,
by
statement $(g)$
of
Theorem
$18.2.11$
of 
\cite{ds},
so for all
$v\in G$
\begin{equation*}
\begin{aligned}
\int\,
f_{x}(R)
\,d\,
\mu(x)
E(\delta_{n})
v
&
\doteq
\int\,
f_{x}(R)
E(\delta_{n})
v
\,d\,
\mu(x)
\\
&
=
\int\,
E(\delta_{n})
f_{x}(R)
v
\,d\,
\mu(x)
=
E(\delta_{n})
\int\,
f_{x}(R)
v
\,d\,
\mu(x),
\end{aligned}
\end{equation*}
where
the last equality
follows
by applying 
Theorem $1$, $IV.35$.
of
\cite{IntBourb}.
Hence for all
$n\in\N$
$$
\int\,
f_{x}(R)
\,d\,
\mu(x)
E(\delta_{n})
G
\subseteq
E(\delta_{n})
G
\subseteq
Dom(g(R)),
$$
where the
last inclusion
is 
by
\eqref{19292703ST}.
Therefore 
$$
(\forall n\in\N)
(\forall v\in G)
\left
(E(\delta_{n})v
\in 
Dom
\left(
g(R)
\int\,
f_{x}(R)
\,d\,
\mu(x)
\right)
\right).
$$
Hence
by
\eqref{20271201}
\begin{equation}
\label{II04041215}
\D
\doteqdot
Dom
\left(
g(R)
\int\,
f_{x}(R)
\,d\,
\mu(x)
\right)
\text{ is dense in }
G.
\end{equation}
Next
$
\int\,
f_{x}(R)
\,d\,
\mu(x)
\in
B(G)
$
and
$g(R)$
is
closed by 
Theorem $18.2.11.$
of \cite{ds},
so by
Lemma
\ref{II04041224}
\begin{equation}
\label{II04041228}
g(R)
\int\,
f_{x}(R)
\,d\,
\mu(x)
\text{ is closed.}
\end{equation}
Moreover
$h(R)\in B(G)$
hence by 
\eqref{15051936}
\begin{equation}
\label{II04041231} 
g(R)
\int\,
f_{x}(R)
\,d\,
\mu(x)
\in
B(\D,G).
\end{equation}
\eqref{II04041228},
\eqref{II04041231} 
and
Lemma 
\ref{II04041227}
allow us to state that
$\D$ is closed in $G$,
thus
by 
\eqref{II04041215}
we have
$$
\D=G.
$$
Therefore by 
\eqref{15051936}
the
statement 
follows.
\end{proof}
Now we shall prove a corollary
of the previous result in which
conditions
are given 
ensuring the strong operator
integrability of the map
$f_{x}(R)$.
\begin{corollary}
\label{17070901}
Let
$R$
be
a 
possibly unbounded 
scalar type spectral operator
in $G$.
Let
$
\{
\sigma_{n}
\}_{n\in\N}
$ 
be an
$E-$sequence
and for all
$x\in X$,
$
f_{x}\in Bor(\sigma(R))
$
such that
$$
(X\ni x
\mapsto
\|\w{f}_{x}\|_{\infty}^{E})
\in
\F{1}{}(X;\mu)
$$
and
$X\ni x\mapsto f_{x}(R)\in B(G)$
satisfies
the conditions $(1-2)$ of 
Proposition
\ref{13051512}, 
(respectively for all
$v\in G$
the map
$
X
\ni x
\mapsto
f_{x}(R)v
\in
G
$
is
$\mu-$measurable).
Finally
let
$
g,h\in
Bor(\sigma(R))
$.
If
we assume
that for all
$n\in\N$
holds
\eqref{20081201}
and
that
$
\w{h}
\in
\Lf{E}{\infty}(\sigma(R))
$,
then
the same conclusions
of
Thm.
\ref{13051634}
hold.
\end{corollary}
\begin{proof}
By
statement $(c)$
of Theorem $18.2.11$ of \cite{ds}
and
Proposition
\ref{13051512}, 
(respectively
Remark \ref{13051626})
the map
$
X\ni x\mapsto f_{x}(R)
\in B(G)
$
is $\mu-$integrable with respect
to the
strong 
operator topology
and
$$
\left\|
\int
f_{x}(R)
\,d\,\mu(x)
\right\|_{B(G)}
\leq
4M
\int^{*}
\|\w{f}_{x}\|_{\infty}^{E}
\,d\,|\mu|(x)
$$
Here
$
M
\doteqdot
\sup_{\sigma\in\B(\C)}
\|E(\sigma)\|_{B(G)}
$.
Therefore
the statement
follows
by Theorem
\ref{13051634}.
\end{proof}
\section{
Generalization
of
the
Newton-Leibnitz
formula
}
The
main result
of this section
is
Theorem
\ref{27051755}
which
generalizes
the Newton-Leibnitz
formula
to the case
of unbounded
scalar type spectral operators
in $G$.
For proving
Theorem
\ref{27051755}
we need
two preliminary
results,
the first is
Theorem 
\ref{19500603},
concerning
the Newton-Leibnitz
formula
for 
any
bounded
scalar type spectral operator
on $G$
and any
analytic map
on an open 
neighbourhood of
its spectrum.
The second,
Theorem
\ref{27051053},
concerns
strong operator 
continuity,
and under additional
conditions
also
differentiability,
for operator maps
of the type
$
K\ni t\mapsto S(tR)\in B(G)
$,
where
$K$
is an open interval
of $\R$,
$S(tR)$
arises
by the Borel functional
calculus
for the unbounded scalar type
spectral operator
$R$ in $G$
and
$S$ is any
analytic map
on an open 
neighbourhood 
$U$
of
$\sigma(R)$
such that
$K\cdot U\subseteq U$.
Let
$Z$
be a non
empty set,
$Y$
a 
$\K-$linear space
($\K\in\{\R,\C\}$),
$U\subseteq Y$,
$K\subseteq\K$
such that
$
K\cdot U
\subseteq
U 
$
and
$F:U\to Z$.
Then 
we set
$
F_{t}:U\to Z
$
such that
$
F_{t}(\lambda)
\doteqdot
F(t\lambda)
$,
for all
$t\in K$
and
$\lambda\in U$.
If
$F,G$
are two
$\C-$Banach spaces,
$A\subseteq F$ 
open
and
$
f:
A\subseteq F\to G
$
a map,
we convene
to denote
the
real Banach
spaces
$
F_{\R}
$
and
$
G_{\R}
$
associated
to
$F$
and
$G$
again
by
$F$
and
$G$
respectively,
and
with the 
symbol
$f$
the
map
$
f^{\R}:
A\subseteq F_{\R}\to G_{\R}
$.

\begin{lemma}
\label{12331303}
Let
$\lr{Y}{d}$ be a metric space,
$U$ an open of $Y$
and
$\sigma$ a compact
such that
$\sigma\subseteq U$.
Then there is
$Q>0$
\begin{equation}
\label{12501303}
K
\doteqdot
\ov{\bigcup_{\{y\in\sigma\}}
\ov{B}_{Q}(y)}
\subseteq U,
\end{equation}
moreover
if
$\sigma$
is 
of finite diameter
then
$K$
is of finite diameter.
\end{lemma}
\begin{proof}
By 
Remark
\S $2.2.$, Ch. $9$
of \cite{BourGT}
we deduce
$$
P
\doteqdot
dist
(\sigma,\complement U)
\neq
0,
$$
where
$
dist(A,B)
\doteqdot
\inf_{\{x\in A, y\in B\}}
d(x,y)
$,
for all 
$A,B\subseteq Y$.
Set
$$
Q
\doteqdot
\frac{P}{2}
$$
then for all
$
y\in\sigma,
x\in \ov{B}_{Q}(y),
z\in\complement U
$
we have
\begin{equation}
\label{23101503}
d(x,z)
\geq
d(z,y)
-
d(y,x)
\geq
\frac{P}{2}
\neq
0.
\end{equation}
Thus
by applying
Proposition
$2$, \S $2.2.$, Ch. $9$
of \cite{BourGT}
$
\ov{B}_{Q}(y)
\cap
\complement U
=
\emptyset
$,
i.e.
$
\ov{B}_{Q}(y)
\subseteq
U
$,
then
$$
A
\doteqdot
\bigcup_{\{y\in\sigma\}}
\ov{B}_{Q}(y)
\subseteq 
U.
$$
Moreover
by
Proposition
$3$, \S $2.2.$, Ch. $9$
of \cite{BourGT}
the map
$
x\mapsto
d(x,\complement U)
$
is 
continuous
on $Y$,
hence
by
\eqref{23101503}
for all
$x\in\ov{A}$
$$
d(x,\complement U)
=
\lim_{n\in\N}
d(x_{n},\complement U)
\geq
\frac{P}{2}
\neq
0,
$$
for all
$\{x_{n}\}_{n\in\N}\subset A$
such that
$x=\lim_{n\in\N}x_{n}$.
Therefore
by
Proposition
$2$, \S $2.2.$, Ch. $9$
of \cite{BourGT}
\eqref{12501303}
follows.
Let
$B\subset Y$
be of finite diameter
then by the continuity
of the map
$d:Y\times Y\to \R^{+}$
it is of finite diameter also
$\ov{B}$.
Indeed
let
$
diam(B)
\doteqdot
\sup_{\{x,y\in B\}}
d(x,y)
$,
if by absurdum
$
\sup_{\{x,y\in\ov{B}\}}
d(x,y)
=
\infty
$
then
\begin{equation}
\label{12392403}
(\exists\,x_{0},y_{0}
\in
\ov{B})
(d(x_{0},y_{0})>diam(B)+1).
\end{equation}
Let
$
\{(x_{\alpha},y_{\alpha})\}_{\alpha\in D}
\subset
B\times B
$
be
a
net
such that
$
\lim_{\alpha\in D}
(x_{\alpha},y_{\alpha})
=
(x_{0},y_{0})
$
limit
in
$
\lr{Y}{d}
\times
\lr{Y}{d}
$.
Thus by the continuity
of $d$
$$
d(x_{0},y_{0})
=
\lim_{\alpha\in D}
d(x_{\alpha},y_{\alpha})
\leq
diam(B)
$$
which 
contradicts
\eqref{12392403},
so
$
\sup_{\{x,y\in\ov{B}\}}
d(x,y)
<
\infty
$.
Therefore
if 
$A$ is of finite diameter
it is so
$K$.
Let
$z_{1},z_{2}\in A$
then
there exist
$
y_{1},y_{2}
\in\sigma
$
such that
$
z_{k}
\in
\ov{B}_{Q}(y_{k})
$,
for
$k\in\{1,2\}$.
Then
$$
d(z_{1},z_{2})
\leq
d(z_{1},y_{1})
+
d(y_{1},y_{2})
+
d(y_{2},z_{2})
\leq
2Q+
diam(\sigma)
<\infty,
$$
where
$
diam(\sigma)
\doteqdot
\sup_{\{x,y\in\sigma\}}
d(x,y)
$.
Hence
$A$
is of finite diameter.
\end{proof}
\begin{theorem}
\label{19500603}
Let
$T\in B(G)$
be a scalar type spectral
operator,
$\sigma(T)$
its spectrum.
Assume that
$
0<L\leq\infty
$,
$U$
is
an open
neighbourhood of
$\sigma(T)$
such that
$
]-L,L[
\cdot
U
\subseteq 
U
$
and
$
F:U\to\C
$
is
an 
analytic
map.
Then for all
$t\in
]-L,L[$
\begin{enumerate}
\item
\begin{equation}
\label{29051755}
F(t T)
=
F_{t}(T);
\end{equation}
\item
\begin{equation}
\label{22350603}
\frac{d\,F(tT)}{d\,t}
=
T
\frac{d\,F}{d\,\lambda}(tT);
\end{equation}
\item
for all
$
u_{1},u_{2}
\in
]-L,L[
$
\begin{equation}
\label{22360603}
T
\oint_{u_{1}}^{u_{2}}
\frac{d\,F}{d\,\lambda}(tT)
d\,t
=
F(u_{2}T)
-
F(u_{1}T).
\end{equation}
\end{enumerate}
Here
$F_{t}(T)$,
(respectively
$
\frac{d\,F}{d\,\lambda}(tT)
$
and
$F(tT)$)
are
the 
operators
arising
by the Borelian
functional
calculus
of the operator
$T$
(respectively
$tT$)
for all
$
t\in
]-L,L[
$.
\end{theorem}
\begin{proof}
$T$
is a bounded operator
on $G$
so 
$\sigma(T)$
is compact.
Let
us denote
by
$
\lr{\cc(\sigma(T))}{\|\cdot\|_{\sup}}
$
the Banach algebra
of all continuous
complex valued 
maps
defined
on
$\sigma(T)$
with the norm
$
\|g\|_{\sup}
\doteqdot
\sup_{\lambda\in\sigma(T)}
|g(\lambda)|
$.
Set
\begin{equation}
\label{13092403}
\begin{cases}
\w{\cc}(\sigma(T))
\doteqdot
\left\{
f:\C\to\C
\mid
f\up\sigma(T)
\in
\cc(\sigma(T)),
\,
f\up\complement
\sigma(T)
=
\ze
\right\},
\\
J:\cc(\sigma(T))
\ni
g
\mapsto
\w{g}
\in
\w{\cc}(\sigma(T)).
\end{cases}
\end{equation}
Notice
that
$\w{\cc}(\sigma(T))$
is an algebra
moreover
$J$
is a surjective
morphism
of algebras
and
$
\sup_{\lambda\in\C}
|J(g)(\lambda)|
=
\|g\|_{\sup}
$
for all
$g\in
\cc(\sigma(T))
$
furthermore
$J(g)\in Bor(\C)$
since
$g\in Bor(\sigma(T))$
and
$\sigma(T)\in\B(\C)$.
Hence
$
\w{\cc}(\sigma(T))
$
is a
subalgebra
of
$
\bb
$,
moreover
$J$
is an 
isometry
between
$
\lr{\cc(\sigma(T))}{\|\cdot\|_{\sup}}
$
and
$
\lr{\w{\cc}(\sigma(T))}{\|\cdot\|_{\sup}}
$.
Thus
$
\lr{\w{\cc}(\sigma(T))}{\|\cdot\|_{\sup}}
$
is a
Banach 
subalgebra
of 
the Banach algebra
$
\lr{\bb}{\|\cdot\|_{\sup}}
$
and
$J$
is
an
isometric
isomorphism
of algebras.
Therefore
by denoting
with
$E$
the
resolution of the
identity
of
$T$,
by
\eqref{19492902}
we have
that
$\I{\C}{E}\circ J$
is a
unital
\footnote{
indeed
by 
setting
$\un:\C\ni\lambda\mapsto 1\in\C$
the unity
element
in
$\bb$
then
$
\I{\C}{E}
\circ
J(\un\up\sigma(T))
=
\I{\C}{E}(\un\cdot\chi_{\sigma(T)})
=
\I{\C}{E}(\un)
\I{\C}{E}(\chi_{\sigma(T)})
=
\un
$.
}
morphism
of algebras
such that
$
\I{\C}{E}\circ J
\in
B(\lr{\cc(\sigma(T))}{\|\cdot\|_{\sup}},B(G))
$.
In the sequel 
we
convene
to denote
for brevity
with
the symbol
$\I{\C}{E}$
the operator
$\I{\C}{E}\circ J$
so 
\begin{equation}
\label{56971303}
\begin{cases}
\I{\C}{E}
\in
B(\lr{\cc(\sigma(T))}{\|\cdot\|_{\sup}},B(G)),
\\
\I{\C}{E}
\text{ is a unital morphism of algebras}
\\
(\forall g\in\cc(\sigma(T)))
(g(T)=\I{\C}{E}(g)).
\end{cases}
\end{equation}
In particular
$
\I{\C}{E}
$
is Fr\'{e}chet
differentiable
with constant
differential
map
equal
to
$\I{\C}{E}$.
In the sequel
we shall
denote
with
$\ze$
the 
zero
element
of
the 
Banach
space
$
\lr{\cc(\sigma(T))}{\|\cdot\|_{\sup}}
$.
Let
$t\in]-L,L[-\{0\}$,
and
$
\imath_{t}
\doteqdot
t\cdot
\imath
$,
where
$
\imath:
\sigma(T)
\ni
\lambda
\mapsto
\lambda
$.
So
$
\imath_{t}(T)
=
\I{\C}{E}(t\cdot\imath)
=
t\I{\C}{E}(\imath)
=
t T
$.
Hence
by 
the general
spectral mapping
theorem
$18.2.21.$
of \cite{ds}
applied
to the map
$\imath_{t}$,
the fact that
$\sigma(T)$
is closed
and the product by 
no zero scalars
in $\C$ is a homeomorphism,
we deduce that
$t T$
is a scalar type
spectral operator
and
$
E_{t}:\B(\C)
\ni
\delta
\mapsto
E(t^{-1}\delta)
$
its resolution of the identity.
Finally
\begin{equation*}
(\forall t\in ]-L,L[)
(\sigma(t T)
=
t\sigma(T)
\subseteq
U),
\end{equation*}
the inclusion
is by hypothesis.
So
$F(t T)$
arising
by 
the Borel functional calculus
of the operator
$t T$
is well defined
and
by
\eqref{56971303}
\begin{equation}
\label{14222403}
\begin{aligned}
F(t T)
=
\I{\C}{E_{t}}(F\up\sigma(t T))
&
\doteq
\I{\C}{E\circ\imath_{t^{-1}}}
(F\up\sigma(t T))
\\
&
=
\I{\C}{E}(F\circ\imath_{t})
\\
&
=
\I{\C}{E}(F_{t}\up\sigma(T))
=
F_{t}(T).
\end{aligned}
\end{equation}
Thus
\eqref{29051755}.
Set
$$
\Delta:
]-L,L[\ni t
\mapsto
(F\circ\imath_{t})
\in
\lr{\cc(\sigma(T))}{\|\cdot\|_{\sup}},
$$
by
the 
third
equality in
\eqref{14222403}
\begin{equation}
\label{18581303}
(\forall t\in]-L,L[)
(F(t T)
=
\I{\C}{E}\circ\Delta(t)).
\end{equation}
We claim
that
$\Delta$
is derivable
(i.e. Fr\'{e}chet differentiable)
and for all
$t\in]-L,L[$
\begin{equation}
\label{19051303}
\frac{d\,\Delta}{d t}(t)
=
\imath
\cdot
\left(\frac{d\,F}{d\lambda}\right)_{t}
\up
\sigma(T).
\end{equation}
Set
$$
\begin{cases}
\cc_{U}(\sigma(T))
\doteqdot
\{
f\in
\cc(\sigma(T))
\mid
f(\sigma(T))
\subseteq
U
\},                                             
\\
\zeta:
]-L,L[\ni t\mapsto\imath_{t}
\in
\cc_{U}(\sigma(T))
\subset
\lr{\cc(\sigma(T))}{\|\cdot\|_{\sup}}
\\
\varUpsilon:
\cc_{U}(\sigma(T))
\ni f
\mapsto
F\circ f
\in
\lr{\cc(\sigma(T))}{\|\cdot\|_{\sup}}.
\end{cases}
$$
Notice
\begin{equation}
\label{16361403}
\Delta=\varUpsilon\circ\zeta,
\end{equation}
moreover
$\zeta$ 
is Fr\'{e}chet differentiable
and for all
$t\in]-L,L[$
\begin{equation}
\label{16351403}
\frac{d\,\zeta}{d t}(t)
=
\imath.
\end{equation}
Next for all
$f\in\cc_{U}(\sigma(T))$
by
Lemma \ref{12331303}
applied to the compact
$
f(\sigma(T))
$,
there is
$Q_{f}>0$
\begin{equation}
\label{12341303}
K_{f}
\doteqdot
\ov{
\bigcup_{\{\lambda\in\sigma(T)\}}
\ov{B}_{Q_{f}}(f(\lambda))
}
\subseteq
U,
\end{equation}
in particular
\begin{equation}
\label{12341303bis}
\ov{B}_{Q_{f}}(f)
\subseteq
\cc_{U}(\sigma(T)).
\end{equation}
Thus
$
\cc_{U}(\sigma(T))
$
is an open set
of
the space
$\lr{\cc(\sigma(T))}{\|\cdot\|_{\sup}}$,
therefore
we can
claim that
$\varUpsilon$ 
is Fr\'{e}chet differentiable
and
its
differential
map
$
\varUpsilon^{[1]}:
\cc_{U}(\sigma(T))
\to
B(\cc(\sigma(T)))
$
is such that for all
$
f\in\cc_{U}(\sigma(T)),
h\in\cc(\sigma(T)),
\lambda\in\sigma(T)
$
\begin{equation}
\label{21391303}
\begin{cases}
\varUpsilon^{[1]}(f)(h)(\lambda)
=
\frac{d\,F}{d\lambda}(f(\lambda))
h(\lambda),
\\
\|\varUpsilon^{[1]}(f)
\|_{B(\cc(\sigma(T)))}
\leq
\|\frac{d\,F}{d\lambda}\circ f\|_{\sup}
\end{cases}
\end{equation}
Let us fix
$f\in\cc_{U}(\sigma(T))$
and
$K_{f}$
as in
\eqref{12341303},
so
by 
the boundedness
of 
$f(\sigma(T))$
and
Lemma \ref{12331303}
$K_{f}$
is compact.
Morever
$
\frac{d\,F}{d\lambda}
$
is continuous
on $U$
therefore
uniformly
continuous
on
the compact
$K_{f}$,
hence
$(\forall\ep>0)
(\exists\,\delta>0)
(\forall h\in
\ov{B}_{Q_{f}}(\ze)
\cap
\ov{B}_{\delta}(\ze))
$
\begin{equation}
\label{23451503}
\sup_{t\in[0,1]}
\sup_{\lambda\in\sigma(T)}
\big|
\frac{d\,F}{d\lambda}
(f(\lambda)+t h(\lambda))
-
\frac{d\,F}{d\lambda}
(f(\lambda))
\big|
\leq
\ep,
\end{equation}
indeed
$
f(\lambda)+t h(\lambda)
\in
K_{f}
$
and
$
|f(\lambda)+t 
h(\lambda)
-f(\lambda)|
\leq
|h(\lambda)|
\leq
\delta
$,
for all
$\lambda\in\sigma(T)$
and
$t\in[0,1]$.
Let us 
identify
for the moment
$\C$ 
as the
$\R-$normed
space
$\R^{2}$,
then
the usual 
product
$
(\cdot):
\C
\times
\C
\to
\C
$
is
$\R-$bilinear,
therefore
the map
$F:U\subseteq\R^{2}\to\R^{2}$
is Fr\'{e}chet differentiable
and for all
$x\in U,h\in\R^{2}$
\begin{equation}
\label{12101403}
F^{[1]}(x)(h)
=
\frac{d\,F}{d\lambda}(x)\cdot h.
\end{equation}
for all
$
h\in
\ov{B}_{Q_{f}}(\ze)
$
\begin{alignat}{2}
\label{00271603}
\sup_{\lambda\in\sigma(T)}
\big|
(F(f(\lambda)+h(\lambda))
-F(f(\lambda))
-
\frac{d\,F}{d\lambda}
(f(\lambda))h(\lambda)
\big|
&
=
\\
\notag
\sup_{\lambda\in\sigma(T)}
\big|
(F(f(\lambda)+h(\lambda))
-F(f(\lambda))
-
F^{[1]}(f(\lambda))(h(\lambda)))
\big|
&
\leq
\\
\notag
\sup_{t\in[0,1]}
\sup_{\lambda\in\sigma(T)}
\|
F^{[1]}(f(\lambda)+th(\lambda))
-
F^{[1]}(f(\lambda))
\|_{B(\R^{2})}
\sup_{\lambda\in\sigma(T)}
|h(\lambda)|
&
=
\\
\notag
\sup_{t\in[0,1]}
\sup_{\lambda\in\sigma(T)}
\left|
\frac{d\,F}{d\lambda}
(f(\lambda)+th(\lambda))
-
\frac{d\,F}{d\lambda}
(f(\lambda))
\right|
\|h\|_{\sup}.
\end{alignat}
Here
in the first
equality
we use
\eqref{12101403},
in the first inequality
an application
of the 
Mean value theorem
applied to the 
Fr\'{e}chet differentiable
map
$F:U\subset\R^{2}\to\R^{2}$,
in the second equality
we use a corollary of
\eqref{12101403}.
Finally
by
\eqref{00271603}
and
\eqref{23451503}
$
(\forall\ep>0)
(\exists\,\delta>0)
(\forall h\in
\ov{B}_{Q_{f}}(\ze)
\cap
\ov{B}_{\delta}(\ze)
-\{\ze\})
$
$$
\frac{
\sup_{\lambda\in\sigma(T)}
\big|
(F(f(\lambda)+h(\lambda))
-F(f(\lambda))
-
\frac{d\,F}{d\lambda}
(f(\lambda))h(\lambda)
\big|
}{\|h\|_{\sup}}
\leq
\ep.
$$
Equivalently
\begin{equation}
\label{16251403}
\lim_{
\begin{subarray}{l}
h\to\ze
\\ 
h\neq\ze
\end{subarray}
}
\frac
{
\|
\varUpsilon(f+h)
-\varUpsilon(f)
-\varUpsilon^{[1]}(f)(h)
\|_{\sup}}
{\|h\|_{\sup}}
=0
\end{equation}
Moreover
$$
\|\varUpsilon^{[1]}(f)(h)\|_{\sup}
\leq
\|
\frac{d\,F}{d\lambda}\circ f
\|_{\sup}
\|
h
\|_{\sup}
$$
then
by
\eqref{16251403}
we proved
the claimed
\eqref{21391303}.
By
\eqref{16361403},
\eqref{16351403}
and
\eqref{21391303}
we deduce
that
$\Delta$
is derivable
in addition for all
$t\in]-L,L[,
\lambda\in\sigma(T)$
\begin{alignat*}{2}
\frac{d\,\Delta}{d t}(t)(\lambda)
&
=
\varUpsilon^{[1]}(\zeta(t))(\imath)(\lambda)
\\
&
=
\frac{d\,F}{d\lambda}(\zeta_{t}(\lambda))
\imath(\lambda)
=
\imath
\left(\frac{d\,F}{d\lambda}\right)_{t}
(\lambda).
\end{alignat*}
Thus
the claimed
\eqref{19051303}.
In conclusion
by
the fact that
$\I{\C}{E}$
is a  
morphism of 
algebras,
\eqref{18581303},
\eqref{56971303}
and
\eqref{19051303}
for all
$t\in]-L,L[$
\begin{alignat*}{1}
\frac{d\,F(tT)}{d t}
&
=
\frac{d}{d t}
\left(\I{\C}{E}\circ\Delta\right)(t)
=
\I{\C}{E}
\left(\frac{d\,\Delta}{d t}(t)\right)
\\
&
=
\I{\C}{E}
\left(\imath\cdot
\left(\frac{d\,F}{d\lambda}\right)_{t}
\up\sigma(T)
\right)
\\
&
=
\I{\C}{E}
\left(\imath\right)
\I{\C}{E}
\left(
\left(\frac{d\,F}{d\lambda}\right)_{t}
\up\sigma(T)
\right)
=
T
\left(\frac{d\,F}{d\lambda}\right)_{t}(T).
\end{alignat*}
Therefore
statement
$(2)$
by statement
$(1)$
applied
to the map
$
\frac{d\,F}{d\lambda}
$.
The
map
$
]-L,L[
\ni
t
\mapsto
\frac{d\,F}{d\lambda}
(tT)
\in B(G)
$
is continuous
by
\eqref{22350603}
(by replacing the map
$F$ 
with
$
\frac{d F}{d\lambda}
$)
hence
it
is 
Lebesgue-measurable
in
$B(G)$.
Let
$
u_{1},u_{2}
\in
]-L,L[
$,
by statement
$(1)$
and
Theorem
$18.2.11.$
of \cite{ds}
\begin{alignat*}{2}
\int_{[u_{1},u_{2}]}^{*}
\left
\|
\frac{d\,F}{d\,\lambda}
(tT)
\right
\|
d
\,
t
&
=
\int_{[u_{1},u_{2}]}^{*}
\left
\|
\left(
\frac{d\,F}{d\,\lambda}
\right)_{t}
(T)
\right
\|
d\,t
\\
&
\leq
4M
\int_{[u_{1},u_{2}]}^{*}
\left
\|
\left(
\frac{d\,F}{d\,\lambda}
\right)_{t}
\up
\sigma(T)
\right
\|_{\sup}
d\,t
\\
&
\leq
4 M D
|u_{2}-u_{1}|
<
\infty,
\end{alignat*}
where
$
M
\doteqdot
\sup_{\delta\in\B(\C)}\|E(\delta)\|
$,
and
$$
D
\doteqdot
\sup_{t\in[u_{1},u_{2}]}
\left
\|
\left(
\frac{d\,F}{d\,\lambda}
\right)_{t}
\up
\sigma(T)
\right
\|_{\sup}
=
\sup_{(t,\lambda)
\in
[u_{1},u_{2}]
\times
\sigma(T)}
\big|
\frac{d\,F}{d\,\lambda}
(t\lambda)
\big|
<
\infty,
$$
indeed
$
[u_{1},u_{2}]
\times
\sigma(T)
$
is compact
and
the map
$
(t,\lambda)
\mapsto
\frac{d\,F}{d\,\lambda}
(t\lambda)
$
is continuous
on
$
]-L,L[
\times
U
$.
Therefore
by
Theorem $5$,
$IV.71$
of
\cite{IntBourb} 
$
]-L,L[
\ni
t
\mapsto
\frac{d\,F}{d\lambda}
(tT)
$
is Lebesgue-integrable 
with respect to the 
norm topology on
$B(G)$,
so in particular
by Definition
$1$,
$IV.33$
of
\cite{IntBourb}
\begin{equation}
\label{09052027aaa}
\exists
\oint_{u_{1}}^{u_{2}}
\frac{d\,F}{d\,\lambda}
(tT)d\,t
\in 
B(G).
\end{equation}
Therefore 
by
\eqref{10501634},
\eqref{09052027aaa},
Theorem $1$,
$IV.35$
of
\cite{IntBourb} 
and
\eqref{22350603}
\begin{equation}
\label{21061013aaa}
T
\oint_{u_{1}}^{u_{2}}
\frac{d\,F}{d\,\lambda}(tT)
d\,t
=
\oint_{u_{1}}^{u_{2}}
T
\frac{d\,F}{d\,\lambda}(tT)
d\,t
=
\oint_{u_{1}}^{u_{2}}
\frac{d\,F(tT)}{d\,t}
dt.
\end{equation}
By 
\eqref{22350603}
the map
$
]-L,L[
\ni
t\mapsto
F(tT)
$,
is
derivable
moreover
its derivative
$
]-L,L[
\ni
t\mapsto
\frac{d\,F(tT)}{d\,t}
$
is
continuous
in
$B(G)$
by 
\eqref{22350603}
and 
the continuity
of
the map
$
]-L,L[
\ni
t\mapsto
\frac{d\,F}{d\lambda}
(tT)
$
in
$B(G)$.
Therefore
$
[u_{1},u_{2}]
\ni
t\mapsto
\frac{d\,F(tT)}{d\,t}
$
is 
Lebesgue integrable
in
$B(G)$,
where 
the integral
has to be understood
as defined in
Ch
$II$
of
\cite{FRV},
see 
Proposition
$3$,
$n^{\circ}3$,
$\S 1$,
Ch
$II$
of
\cite{FRV}.
Finally
the 
Lebesgue integral
for functions with values
in a Banach space
as defined in 
Ch
$II$
of
\cite{FRV},
turns
to be
the 
integral
with respect to the Lebesgue
measure
as defined in
Ch. $IV$, \S $4$, $n^{\circ} 1$ 
of
\cite{IntBourb}
(see 
Ch $III$, \S $3$, $n^{\circ} 3$
and 
example in 
Ch $IV$, \S $4$, $n^{\circ} 4$ 
of 
\cite{IntBourb}).
Thus
statement 
$(3)$
follows
by
\eqref{21061013aaa}.
\end{proof}
\begin{lemma}
\label{27051011}
Let $R$ be a 
possibly
unbounded
scalar type 
spectral operator in $G$, 
$\sigma(R)$ its spectrum,
$E$ its resolution of the identity,
$K\neq\emptyset$ 
and for all
$t\in K$
be
$
f_{t}
\in
Bor(\sigma(R))
$
such that
\begin{equation}
\label{18130502}
N
\doteqdot
\sup_{t\in K}
\|\w{f_{t}}\|_{\infty}^{E}
<
\infty.
\end{equation}
If
$
g
\in 
Bor(\sigma(R))
$
and
$
\{\sigma_{n}\}_{n\in\N}
$
is an
$E-$sequence
then for all
$v\in Dom(g(R))$
$$
\lim_{n\in\N}
\sup_{t\in K}
\left\|
f_{t}(R)
g(R)
v
-
f_{t}(R)
g(R)
E(\sigma_{n})
v
\right\|
=
0.
$$
\end{lemma}
\begin{proof}
By 
statement
$(g)$
of
Theorem
$18.2.11.$ of \cite{ds}
the statement is well set.
Let
$
M
\doteqdot
\sup_{\sigma\in\B(\C)}\|E(\sigma)\|_{B(G)}
$
then
$M<\infty$
by
Corollary
$15.2.4.$
of \cite{ds}.
Hypothesis 
\eqref{18130502}
together 
statement $(c)$
of Theorem $18.2.11.$
of \cite{ds},
imply
that for all $t\in K$, 
$f_{t}(R)\in B(G)$
and
$$
\sup_{t\in K}
\|f_{t}(R)\|_{B(G)}
\leq
4 M N.
$$
Therefore for all
$v\in Dom(g(R))$
we have
\begin{alignat*}{2}
&
\lim_{n\in\N}
\sup_{t\in K}
\left\|
f_{t}(R)
g(R)
v
-
f_{t}(R)
g(R)
E(\sigma_{n})
v
\right\|
\\
&
\leq
\lim_{n\in\N}
\sup_{t\in K}
\left\|
f_{t}(R)
\right\|
\cdot
\left\|
g(R)
v
-
g(R)
E(\sigma_{n})
v
\right\|
\\
&
\leq
4 M N
\lim_{n\in\N}
\left\|
g(R)
v
-
g(R)
E(\sigma_{n})
v
\right\|
\\
&
=
4 M N
\lim_{n\in\N}
\left\|
g(R)
v
-
E(\sigma_{n})
g(R)
v
\right\|
&
\text{ 
by $(g)$
of Theorem
$18.2.11.$ of \cite{ds}
}
\\
&
=
0
&
\text{ by \eqref{II03041401}.}
\end{alignat*}
\end{proof}
\begin{theorem}
[
\textbf{
Strong operator 
derivability
}
]
\label{27051053}
Let $R$ be a 
possibly unbounded
scalar type 
spectral operator in $G$, 
$K\subseteq\R$
an open interval 
of $\R$
and
$U$
an open
neighbourhood
of $\sigma(R)$
such that
$
K\cdot U\subseteq U
$.
Assume
that
$
f:U\to\C
$
is
an
analytic
map
and
$$
\sup_{t\in K}
\|\w{f_{t}}\|_{\infty}^{E}
<
\infty.
$$
Then
\begin{enumerate}
\item
the map
$
K
\ni t
\mapsto
f(tR)
\in B(G)
$
is
continuous in 
the strong operator topology,
\item
if
\begin{equation}
\label{11581903}
\sup_{t\in K}
\left\|
\w{
\left(
\frac{d\,f}{d\,\lambda}
\right)_{t}
}
\right\|_{\infty}^{E}
<
\infty,
\end{equation}
then for all
$v\in Dom(R),t\in K$
$$
\frac{d f(tR)v}{d t}
=
R
\frac{d f}{d\lambda}
(tR)
v
\in 
G.
$$
\end{enumerate}
\end{theorem}
\begin{proof}
Let
$
\{\sigma_{n}\}_{n\in\N}
$
be
an
$E-$sequence
of compact
sets,
then 
by Lemma \ref{27051011}
applied 
for 
the Borelian map
$g:\sigma(R)\ni\lambda\to 1\in\C$,
so
$g(R)=\un$,
and
by
\eqref{29051755}
we have for all
$v\in G$
\begin{equation}
\label{27051611}
\lim_{n\in\N}
\sup_{t\in K}
\left\|
f(tR)
v
-
f(tR)
E(\sigma_{n})
v
\right\|
=
0.
\end{equation}
By
\eqref{29051755}
and
Key
Lemma
\ref{II31031834}
for all
$
n\in\N
$
\begin{equation}
\label{27051615}
\begin{aligned}
f(tR)
E(\sigma_{n})
&
=
f_{t}(R)E(\sigma_{n})
=
f_{t}(R_{\sigma_{n}}\up G_{\sigma_{n}})
E(\sigma_{n})
\\
&
=
f(t(R_{\sigma_{n}}\up G_{\sigma_{n}}))
E(\sigma_{n}).
\end{aligned}
\end{equation}
$\sigma_{n}$
is
bounded
so
by
Key
Lemma 
\ref{II31031834}
$R_{\sigma_{n}}\up G_{\sigma_{n}}$
is a scalar
type spectral operator
such that
$
R_{\sigma_{n}}\up G_{\sigma_{n}}
\in B(G_{\sigma_{n}})$,
moreover
by
\eqref{25052019}
$U$
is
an
open
neighbourhood of
$\sigma(R_{\sigma_{n}}\up G_{\sigma_{n}})$.
Thus
by 
statement $(2)$
of Theorem
\ref{19500603}
the map
$$
K\ni t\mapsto
f(t(R_{\sigma_{n}}\up G_{\sigma_{n}}))
\in
B(G_{\sigma_{n}})
$$
is
derivable,
so in particular 
$\|\cdot\|_{B(G_{\sigma_{n}})}-$continuous. 
Now for all
$
n\in\N,
v_{n}\in G_{\sigma_{n}}
$
define
$
\xi_{v_{n}}:B(G_{\sigma_{n}})
\ni A\mapsto A v_{n}\in G
$,
then
$\xi_{v_{n}}\in B(B(G_{\sigma_{n}}),G)$,
thus
as a
composition of two continuous maps
also
the following
map
\begin{equation}
\label{13221503}
K\ni t\mapsto
\xi_{E(\sigma_{n})v}
\left(
f(t(R_{\sigma_{n}}\up G_{\sigma_{n}}))
\right)
\in
G
\end{equation}
is
$\|\cdot\|_{G}-$continuous,
for all
$n\in\N,v\in G$.
Hence by 
\eqref{27051615}
we 
have for all
$n\in\N$
\begin{equation}
\label{27051623}
K\ni t\mapsto
f(tR)E(\sigma_{n})
\in B(G)
\text{ is strongly continuous.}
\end{equation}
Finally
by
\eqref{27051623}
and
\eqref{27051611}
we
can apply Theorem
$2$, $\S 1.6.$, Ch. $10$
of \cite{BourGT}
to the 
uniform
space
$B(G)_{st}$
whose uniformity
is
generated by the
set of seminorms
defining the strong
operator
topology on $B(G)$.
Thus we
conclude that
$K\ni t\mapsto
f(tR)
\in B(G)$
is strongly continuous,
and
statement $(1)$
follows.
Let
$n\in\N$
and
$v_{n}\in G_{\sigma_{n}}$
so
$
\xi_{v_{n}}
\in 
B(B(G_{\sigma_{n}}),G)
$
thus
$\xi_{v_{n}}$ is Fr\'{e}chet 
differentiable with constant
differential map
$
\xi_{v_{n}}^{[1]}: 
B(G_{\sigma_{n}})
\ni
A
\mapsto
\xi_{v_{n}}
\in
B(B(G_{\sigma_{n}}),G)
$.
Therefore
by
statement
$(2)$
of
Theorem
\ref{19500603}
for all
$
n\in\N,
v\in G
$
the map
in
\eqref{13221503}
is
Fr\'{e}chet differentiable
as 
composition
of two Fr\'{e}chet differentiable maps,
and
its derivative is
for all
$t\in K$
\begin{alignat}{2}
\label{27051635}
\frac{d}{dt}
\left(
f(t(R_{\sigma_{n}}\up G_{\sigma_{n}}))
E(\sigma_{n})
v
\right)
&
=
\xi_{E(\sigma_{n})v}
\left(
\frac{d}{dt}
(f(t(R_{\sigma_{n}}\up G_{\sigma_{n}})))
\right)
\notag\\
&
=
\frac{d}{dt}
(f(t(R_{\sigma_{n}}\up G_{\sigma_{n}})))
E(\sigma_{n})v
\notag\\
&
=
(R_{\sigma_{n}}\up G_{\sigma_{n}})
\frac{df}{d\lambda}
(t(R_{\sigma_{n}}\up G_{\sigma_{n}}))
E(\sigma_{n})v,
&
\text{ by \eqref{22350603}}
\notag\\
&
=
\frac{df}{d\lambda}
(t(R_{\sigma_{n}}\up G_{\sigma_{n}}))
(R_{\sigma_{n}}\up G_{\sigma_{n}})
E(\sigma_{n})v,
&
\text{ by 
$18.2.11.$, \cite{ds}}
\notag\\
&
=
\left(\frac{df}{d\lambda}\right)_{t}
(R_{\sigma_{n}}\up G_{\sigma_{n}})
(R_{\sigma_{n}}\up G_{\sigma_{n}})
E(\sigma_{n})
v
&
\text{ by 
\eqref{29051755}
}
\notag\\
&
=
\left(\frac{df}{d\lambda}\right)_{t}(R)
(R_{\sigma_{n}}\up G_{\sigma_{n}})
E(\sigma_{n})
v
&
\text{ by 
Lemma
\ref{II31031834}
}
\notag\\
&
=
\frac{df}{d\lambda}(tR)
(R_{\sigma_{n}}\up G_{\sigma_{n}})
E(\sigma_{n})
v.
&
\text{ by 
\eqref{29051755}
}
\end{alignat}
Thus
by
\eqref{27051615}
for all
$
n\in\N,
v\in G
$
\begin{equation}
\label{27051647}
\begin{cases}
K\ni t
\mapsto 
f(tR)E(\sigma_{n})
v
\in G
\text{ is differentiable and}\\
K\ni t
\mapsto 
\frac{df}{d\lambda}(tR)
(R_{\sigma_{n}}\up G_{\sigma_{n}})
E(\sigma_{n})
v
\in G
\text{ is its derivative.}
\end{cases}
\end{equation}
By
\eqref{11581903}
we can apply 
Lemma \ref{27051011}
to the maps
$
\left(
\frac{d f}{d\lambda}
\right)_{t}
\up
\sigma(R)
$
and
$g=
\imath:\sigma(R)\ni\lambda\mapsto\lambda\in\C$,
so
$g(R)=R$,
hence
by
\eqref{29051755} 
for all
$v\in Dom(R)$
\begin{equation}
\label{27051652}
\lim_{n\in\N}
\sup_{t\in K}
\left\|
\frac{d f}{d\lambda}
(tR)
R
v
-
\frac{d f}{d\lambda}
(tR)
(R_{\sigma_{n}}\up G_{\sigma_{n}})
E(\sigma_{n})
v
\right\|
=
0.
\end{equation}
Moreover for all
$a\in K$
let
$r_{a}\in\R^{+}$
be
such that
$
B_{r_{a}}(a)\subset K
$
which 
exists
$K$
being 
open,
then
the
equations
\eqref{27051652},
\eqref{27051647}
and
\eqref{27051611}
hold
again
if we replace
$K$
by
$B_{r_{a}}(a)$.
Hence
we can 
apply
Theorem $8.6.3.$
of \cite{Dieud1}
and deduce for all
$
v\in Dom(R)
$
that
the map
$
K\ni
t\mapsto
f(tR)
v\in 
G
$
is
derivable,
and
its derivative map is
$$
K\ni
t\mapsto
\frac{d f}{d\lambda}(tR)
R v
\in 
G.
$$
Finally
for all
$v\in Dom(R)$,
$R
\frac{d f}{d\lambda}
(tR)v=
\frac{d f}{d\lambda}(tR)
Rv$,
by
$
Dom(\frac{d f}{d\lambda}
(tR))
=
G
$
and
the commutativity property of the 
Borel functional calculus
expressed in
statement $(f)$
of Theorem
$18.2.11.$
of
\cite{ds}.
Hence
the statement
follows.
\end{proof}
\begin{corollary}
\label{II03041154}
Let
$R$
be
a 
possibly unbounded 
scalar type spectral operator
in 
$G$,
$U$
an 
open
neighbourhood of
$\sigma(R)$
and
$S:U\to\C$
an analytic map.
Assume 
that there is
$L>0$
such that
$]-L,L[\cdot
U\subseteq U$
and
\begin{enumerate}
\item
$
\w{S_{t}}
\in
\Lf{E}{\infty}(\sigma(R))
$
and
$
\w{
\left(
\frac{d\,S}{d\,\lambda}
\right)_{t}
}
\in
\Lf{E}{\infty}(\sigma(R))
$
for all 
$t\in ]-L,L[$;
\item
$$
\int^{*}
\left\|
\w{
\left(
\frac{d\,S}{d\,\lambda}
\right)_{t}
}
\right\|_{\infty}^{E}
\,
dt
<\infty
$$
(here
the upper
integral
is with respect
to the
Lebesgue measure on
$]-L,L[$);
\item
for all
$v\in G$
the
map
$
]-L,L[
\ni
t
\mapsto
\frac{d\,S}{d\,\lambda}
(tR)v
\in
G
$
is
Lebesgue
measurable.
\end{enumerate}
Then
for all
$
u_{1},u_{2}\in ]-L,L[
$
$$
R
\int_{u_{1}}^{u_{2}}
\frac{d\,S}{d\,\lambda}(t R)
\,d\,t
=
S(u_{2} R)
-
S(u_{1} R)
\in B(G).
$$
Here
the
integral
is
with respect to the Lebesgue measure on
$[u_{1},u_{2}]$
and
with respect
to
the strong operator topology on $B(G)$,
see
Definition
\ref{13051513}.
\end{corollary}
\begin{proof}
Let
$
M
\doteqdot
\sup_{\sigma\in\B(\C)}
\|E(\sigma)\|_{G}
$
and
$\mu$
the
Lebesgue measure on
$
[u_{1},u_{2}]
$,
then
by
\eqref{29051755},
hypotheses,
and
statement $(c)$
of
Theorem $18.2.11$
of \cite{ds}
we
have
\begin{description}
\item[a]
for all
$t\in [u_{1},u_{2}]$,
$S(t R)
\in B(G)$;
\item[b]
for all
$t\in [u_{1},u_{2}]$,
$\frac{d\,S}{d\,\lambda}
(t R)
\in B(G)$;
\item[c]
$
(
[u_{1},u_{2}]
\ni t
\mapsto
\|
\frac{d\,S}{d\,\lambda}
(t R)
\|_{B(G)}
)
\in
\F{1}{}([u_{1},u_{2}];\mu)
$.
\end{description}
So
by
hypothesis
$(3)$,
the
$(c)$
and
Remark
\ref{13051626}
we have that
the map
$$
[u_{1},u_{2}]
\ni
t
\mapsto
\frac{d\,S}{d\,\lambda}
(tR)
\in
B(G)
$$
is
Lebesgue
integrable
with respect to 
the strong operator topology.
This means that, 
except
for
\eqref{20081201},
the hypotheses of 
Theorem \ref{13051634}
hold
for
$
X\doteqdot 
[u_{1},u_{2}]
$,
$
h
\doteqdot 
(S_{u_{2}}
-
S_{u_{1}})
\up\sigma(R)
$,
$
g:
\sigma(R)
\ni
\lambda
\mapsto
\lambda
\in
\C
$
and finally
for
the 
maps
$
f_{t}
\doteqdot
\left(\frac{d\,S}{d\,\lambda}\right)_{t}
\up
\sigma(R)
$,
for all $t\in [u_{1},u_{2}]$.
Next
let
$
\sigma\in\B(\C)
$
be
bounded,
so
by
Key
Lemma 
\ref{II31031834}
$R_{\sigma}\up G_{\sigma}$
is a scalar
type spectral operator
such that
$
R_{\sigma}\up G_{\sigma}
\in B(G_{\sigma})$,
moreover
by
\eqref{25052019}
$U$
is
an
open
neighbourhood of
$\sigma(R_{\sigma}\up G_{\sigma})$.
Thus
we can 
apply
statement $(3)$
of Theorem
\ref{19500603}
to the
Banach space
$G_{\sigma}$,
the analytic
map
$S$
and
to the operator
$R_{\sigma}\up G_{\sigma}$.
In particular
the map
$
[u_{1},u_{2}]
\ni
t
\mapsto
\frac{d\,S}{d\,\lambda}
(t (R_{\sigma}\up G_{\sigma}))
\in
B(G_{\sigma})
$
is
Lebesgue integrable
in
$\|\cdot\|_{B(G_{\sigma})}-$topology,
that
is
in the meaning
of
Definition $2$,
$IV.23$
of
\cite{IntBourb}.
Next
we consider for all
$v\in G_{\sigma}$,
the
following
map
$$
T\in B(G_{\sigma})\mapsto T v\in G_{\sigma}
$$
which
is
linear
and
continuous in the 
norm
topologies.
Thus
by
Theorem
$1$, $IV.35$
of the
\cite{IntBourb},
$
[u_{1},u_{2}]
\ni
t
\mapsto
\frac{d\,S}{d\,\lambda}
(t(R_{\sigma}\up G_{\sigma}))
v
\in
G_{\sigma}
$
is
Lebesgue integrable
for all
$v\in G_{\sigma}$
and 
$$
\int_{u_{1}}^{u_{2}}
\frac{d\,S}{d\,\lambda}(t(R_{\sigma}\up G_{\sigma}))
v
\,d\,t
=
\left(
\oint_{u_{1}}^{u_{2}}
\frac{d\,S}{d\,\lambda}(t(R_{\sigma}\up G_{\sigma}))
\,d\,t
\right)
v.
$$
Therefore 
by Definition
\ref{13051513}
we 
can state that
$
[u_{1},u_{2}]
\ni
t
\mapsto
\frac{d\,S}{d\,\lambda}(t(R_{\sigma}\up G_{\sigma}))
\in
B(G_{\sigma})
$
is
Lebesgue 
integrable
with respect to
the strong operator topology
on $B(G_{\sigma})$
and
\begin{equation}
\label{15052324}
\int_{u_{1}}^{u_{2}}
\frac{d\,S}{d\,\lambda}(t(R_{\sigma}\up G_{\sigma}))
\,d\,t
=
\oint_{u_{1}}^{u_{2}}
\frac{d\,S}{d\,\lambda}(t(R_{\sigma}\up G_{\sigma}))
\,d\,t.
\end{equation}
Here
$
\int_{u_{1}}^{u_{2}}
\frac{d\,S}{d\,\lambda}(t(R_{\sigma}\up G_{\sigma}))
d t
$
is the
integral of 
$\frac{d\,S}{d\,\lambda}(t(R_{\sigma}\up G_{\sigma}))$
with respect to the 
Lebesgue 
measure on
$[u_{1},u_{2}]$
and
the strong operator topology
on $B(G_{\sigma})$.
Furthermore
by
statement $(3)$
of Theorem
\ref{19500603}
\begin{equation*}
(R_{\sigma}\up G_{\sigma})
\oint_{u_{1}}^{u_{2}}
\frac{d\,S}{d\,\lambda}
(t(R_{\sigma}\up G_{\sigma}))
\,d\,t
=
S(u_{2}(R_{\sigma}\up G_{\sigma}))
-
S(u_{1}(R_{\sigma}\up G_{\sigma})).
\end{equation*}
Thus
by
\eqref{15052324}
\begin{equation}
\label{15052333}
(R_{\sigma}\up G_{\sigma})
\int_{u_{1}}^{u_{2}}
\frac{d\,S}{d\,\lambda}
(t(R_{\sigma}\up G_{\sigma}))
\,d\,t
=
S(u_{2}(R_{\sigma}\up G_{\sigma}))
-
S(u_{1}(R_{\sigma}\up G_{\sigma})).
\end{equation}
Which
implies
\eqref{20081201},
by choosing
for example
$
\sigma_{n}
\doteqdot 
B_{n}(\ze)
$,
for all
$n\in\N$.
Therefore
by
Theorem \ref{13051634}
we obtain
the
statement.
\end{proof}
\begin{theorem}
[
\textbf{
Strong operator
Newton-Leibnitz
formula
}
]
\label{27051755}
Let
$R$
be
a 
possibly unbounded 
scalar type spectral operator
in 
$G$,
$U$
an 
open
neighbourhood of
$\sigma(R)$
and
$S:U\to\C$
an analytic map.
Assume 
that there is
$L>0$
such that
$]-L,L[
\cdot 
U\subseteq U$
and
\begin{enumerate}
\item
$
\w{S_{t}}
\in
\Lf{E}{\infty}(\sigma(R))
$
for all
$t\in ]-L,L[$;
\item
$$
\sup_{t\in]-L,L[}
\left\|
\w{
\left(
\frac{d\,S}{d\,\lambda}
\right)_{t}
}
\right\|_{\infty}^{E}
<\infty.
$$
\end{enumerate}
Then
\begin{enumerate}
\item
for all
$
u_{1},u_{2}\in ]-L,L[
$
$$
R
\int_{u_{1}}^{u_{2}}
\frac{d\,S}{d\,\lambda}(t R)
\,d\,t
=
S(u_{2} R)
-
S(u_{1} R)
\in B(G).
$$
Here
the
integral
is
with respect to the 
Lebesgue measure on
$[u_{1},u_{2}]$
and
with respect to
the strong operator topology on $B(G)$.
\item
If
also
$
\sup_{t\in]-L,L[}
\left\|
\w{S_{t}}
\right\|_{\infty}^{E}
<\infty
$,
then for all
$
v\in Dom(R),
t\in ]-L,L[$
$$
\frac{d S(t R)v}{d t}
=
R
\frac{d S}{d\lambda}
(t R)
v.
$$
\end{enumerate}
\end{theorem}
\begin{proof}
By
hypothesis $(2)$
and
statement $(1)$
of
Theorem
\ref{27051053}
for all
$v\in G$
the
map
$
]-L,L[
\ni
t
\mapsto
\frac{d\,S}{d\,\lambda}
(tR)v
\in
G
$
is
continuous.
Thus 
statement
$(1)$
by
Corollary
\ref{II03041154}
and the fact that
continuity implies measurability.
Statement $(2)$
follows
by
statement $(2)$
of
Theorem 
\ref{27051053}.
\end{proof}
\begin{remark}
\label{13431502}
We end this section
by
remarking
that
$
f:X
\to B(G)
$
is $\mu-$integrable
with respect to the strong operator
topology
as defined in
Definition
\ref{13051513},
if and only if
$
f:X
\to B(G)
$
is scalarly
$(\mu,B(G))-$integrable
with respect
to the 
weak operator topology
in the sense explained
in Notations
\ref{16051536}.
In Chapter \ref{SECWEAKINT}
we shall extend
the results
of Chapter
\ref{31050841}
to the case of integration
with respect to the measure
$\mu$
and
with respect to the
$\sigma(B(G),\n)-$topology,
where $\n\subset B(G)^{*}$
is a suitable linear subspace
of the topological dual of
$B(G)$.
\end{remark}
\section{
Application
to
resolvents 
of 
unbounded 
scalar
type
spectral
operators
in a
Banach
space
$G$
}
\begin{corollary}
\label{04041425}
Let
$T$
be a 
possibly 
unbounded
scalar type
spectral operator
in $G$
with real
spectrum
$\sigma(T)$.
Then
\begin{enumerate}
\item
for all
$
\lambda\in\C\mid Im(\lambda)>0
$
\begin{equation}
\label{21440303}
(T-\lambda\un)^{-1}
=
i
\int_{-\infty}^{0}
e^{-it\lambda}
e^{ it T}
\,d\,t
\in B(G).
\end{equation}
\item
for all
$v\in Dom(T),
t\in\R$
$$
\frac{d\,e^{it(T-\lambda\un)}v}{d\,t}
=
i
(T-\lambda\un)
e^{i(T-\lambda\un)t}
v.
$$
\end{enumerate}
\end{corollary}
\begin{remark}
If we set
the map
$S(\lambda)\doteqdot\exp(i\lambda)$
for all
$\lambda\in\C$
then
the 
operator
functions
in Corollary
\ref{04041425}
are so defined
$
e^{ it T}
\doteqdot
S_{t}(T)
$
and
$
e^{it(T-\lambda\un)}
\doteqdot
S_{t}(T-\lambda\un)
$,
in the sense
of the
Borelian
functional
calculus
for the
scalar type
spectral operators
$T$
and
$(T-\lambda\un)$,
respectively,
as
defined
in Definition 
\ref{20282412}.

The integral
in
Corollary
\ref{04041425}
is
with respect 
to the Lebesgue measure
and 
with respect
to
the strong operator topology on $B(G)$.
Meaning
by definition
that
$$
\int_{-\infty}^{0}
e^{-it\lambda}
e^{ it T}
\,d\,t
\in B(G)
$$
such that
for all
$v\in G$
$$
\left(
\int_{-\infty}^{0}
e^{-it\lambda}
e^{ it T}
\,d\,t
\right)
v
\doteqdot
\lim_{u\to-\infty}
\left(
\int_{u}^{0}
e^{-it\lambda}
e^{ it T}
\,d\,t
\right)
v
=
\lim_{u\to-\infty}
\int_{u}^{0}
e^{-it\lambda}
e^{ it T}
v
\,d\,t.
$$
Here the integral
in the
right
side
of the first equality
is
with respect 
to the Lebesgue measure
on $[u,0]$
and
with respect 
to
the strong operator topology on $B(G)$.
\end{remark}
\begin{proof}
Let 
$\lambda\in\C$
and
set
$R\doteqdot (T-\un\lambda)$,
then $R$
is a scalar type spectral operator,
see Theorem $18.2.17.$ of the \cite{ds}.
Let
$\lambda\in\C\mid Im(\lambda)\neq 0$
and
$E$ be the resolution of the identity of $R$,
then
$\sigma(R)=\sigma(T)-\lambda$,
as a corollary
of the well-known
spectral mapping theorem.
Then for all
$t\in\R$
\begin{equation*}
\begin{aligned}
E-ess\sup_{\nu\in\sigma(R)}
|\frac{d S}{d\lambda}(t \nu)|
&
=
E-ess
\sup_{\nu\in\sigma(R)}
|S(t \nu)|
\leq
\\
&
\leq
\sup_{\nu\in\sigma(R)}
|S(t \nu)|
\\
&
=
\sup_{\mu\in\sigma(T)}
|e^{i(\mu-\lambda)t}|
\\
&
=
e^{Im(\lambda)t}.
\end{aligned}
\end{equation*}
Therefore
are verified the hypotheses
of 
Corollary \ref{27051755}
with the position
$R\doteqdot(T-\lambda\un)$,
then
we can state for all
$v\in G,
u\in\R$
that
\begin{equation}
\label{29051954}
i
(T-\lambda\un)
\int_{u}^{0}
e^{it(T-\lambda\un)}
v
\,d\,t
=
v
-
e^{iu(T-\lambda\un)}
v.
\end{equation}
Here
$
e^{it(T-\lambda\un)}
\doteqdot
S_{t}(R)
$.
One should note
an apparent
ambiguity
about
the symbol
$
e^{it(T-\lambda\un)}
$,
standing 
here for the operator
$
S_{t}(R)
=
S(t R)
$,
which
could be
seen also as a Borelian 
function of the operator $T$.
By 
setting
$
g^{[\lambda]}
(\mu)
\doteqdot 
\mu-\lambda
$,
so
$
g^{[\lambda]}
=
\imath
-
\lambda
\cdot
\un
$
with
$
\un:\C\ni\lambda\mapsto 1
$,
considering 
that
by
the composition
rule,
see Theorem 
$18.2.24$ of \cite{ds}, 
we
have
$
S_{t}\circ g^{[\lambda]}(T)
=S_{t}(g^{[\lambda]}(T))
$,
finally
$
R
=
\imath(T)
-
\lambda
\un
(T)
=
\left(
\imath
-
\lambda
\cdot
\un
\right)
(T)
=
g^{[\lambda]}(T)
$, 
we can assert
\begin{equation}
\label{29052028}
\begin{cases}
T-\lambda\un
=
g^{[\lambda]}(T)
\doteqdot
T-\lambda
\\
e^{it(T-\lambda\un)}
\doteqdot
S_{t}(T-\lambda\un)
=
S_{t}\circ g^{[\lambda]}(T)
=
e^{it(T-\lambda)}.
\end{cases}
\end{equation}
Therefore
we can consider
the operator
$
e^{it(T-\lambda\un)}
$
as an
operator
function of $R$
or 
of $T$.
Now
for all
$t\in\R$,
$
\sup_{\mu\in\sigma(T)}
|\exp(i \mu t)|
=
1
$,
therefore
we can deduce by 
statement $(c)$
of Theorem $18.2.11.$
of \cite{ds}
\begin{equation}
\label{29052051}
\sup_{t\in\R}
\|\exp(iT t)\|_{B(G)}
\leq
4M.
\end{equation}
Here
$
M\doteqdot
\sup_{\sigma\in\B(\C)}
\|E(\sigma)\|_{G}
$.
But with the notations before adopted we have 
for all
$\mu\in\C$
that
$
S_{t}\circ g^{[\lambda]}(\mu)
=
\exp(it(\mu-\lambda))
=
\exp(-it\lambda)
S_{t}(\mu)
$,
therefore
by considering
that
$S_{t}(T)=S(tT)$,
see 
\eqref{29051755},
we have
$
S_{t}\circ g^{[\lambda]}(T)
=
\exp(-it\lambda)
S_{t}(T)
=
\exp(-it\lambda)
S(t T)
$.
Thus
by 
\eqref{29052028}
we have for all
$t\in\R,
\lambda\in\C\mid Im(\lambda)>0$
\begin{equation}
\label{29052047}
e^{it(T-\lambda\un)}
=
\exp(-it\lambda)
S(t T)
\doteq
\exp(-it\lambda)
e^{ it T}.
\end{equation}
We have
by
\eqref{29052047}
and
\eqref{29052051}
$$
\lim_{u\to-\infty}
\|
e^{iu(T-\lambda\un)}
\|_{B(G)}
\leq
4M
\lim_{u\to-\infty}
\exp(Im(\lambda)u)
=
0
$$
or equivalently
$
\lim_{u\to-\infty}
e^{iu(T-\lambda\un)}
=
\ze
$
in 
$\|\cdot\|_{B(G)}-$topology.
Hence by \eqref{29051954}
for all
$v\in G$
\begin{equation}
\label{29052109}
v
=
i
\lim_{u\to-\infty}
(T-\lambda\un)
\int_{u}^{0}
e^{it(T-\lambda\un)}
v
\,d\,t
\text{ in }\|\cdot\|_{G}.
\end{equation}
By considering that
$Im(\lambda)\neq 0$
we have
$
\{\mu\in\C\mid g^{[\lambda]}(\mu)=0\}
\cap
\sigma(T)
=
\emptyset
$,
therefore
if 
we denote with 
$F$
the resolution of the identity of 
the spectral operator
$T$,
we have 
$F(\sigma(T))=\un$
so
$
F(\{\mu\in\C\mid g^{[\lambda]}(\mu)=0\})
=
F(\{\mu\in\C\mid g^{[\lambda]}(\mu)=0\}
\cap
\sigma(T))
=
F(\emptyset)
\doteqdot
\ze
$.
Thus
by applying statement $(h)$
of
Theorem $18.2.11.$
of \cite{ds},
we can
assert
that
$$
\exists\,
(T-\lambda)^{-1}
=
\frac{1}{g^{[\lambda]}}
(T)
\doteq
\frac{1}{T-\lambda}.
$$
Finally
$
F-ess\sup_{\mu\in\sigma(T)}
|\frac{1}{g^{[\lambda]}(\mu)}|
\leq
\sup_{\mu\in\sigma(T)}
|\frac{1}{g^{[\lambda]}(\mu)}|
=
\sup_{\mu\in\sigma(T)}
|\frac{1}{\mu-\lambda}|
=
\frac{1}{\inf_{\mu\in\sigma(T)}|(\mu-\lambda)|}
\leq
\frac{1}{|Im(\lambda)|}
<
\infty
$,
so
$$
\frac{1}{g^{[\lambda]}}
(T)
\in
B(G).
$$
Hence 
by the previous equation
and
the
fact
$
T-\lambda
=
T-\lambda\un
$,
see \eqref{29052028},
we can state 
$$
(T-\lambda\un)^{-1}
\in
B(G).
$$
Finally
by following
a
standard argument,
see for example
\cite{laursen},
by 
this one
and
\eqref{29052109}
we can deduce 
for all
$v\in G$
that
\begin{equation*}
\begin{aligned}
(T-\lambda\un)^{-1}
v
&
=
i
\lim_{u\to-\infty}
(T-\lambda\un)^{-1}
(T-\lambda\un)
\int_{u}^{0}
e^{it(T-\lambda\un)}
v
\,d\,t
\\
&
=
i
\lim_{u\to-\infty}
\int_{u}^{0}
e^{it(T-\lambda\un)}
v
\,d\,t.
\end{aligned}
\end{equation*}
So
statement $(1)$
by
\eqref{29052047}.
By \eqref{29052047},
the fact
that
$S_{t}(T)=S(tT)$
and
statement $(2)$
of Theorem
\ref{27051053}
applied
to
the operator
$T$
and
to the
map
$S:\C\ni\mu\mapsto e^{i\mu}$,
we obtain 
statement $(2)$.
\end{proof}
\begin{remark}
An important application
of 
this formula
is made 
in proving
the well-known
Stone theorem
for strongly continuous
semigroups of unitary
operators in Hilbert space,
see Theorem $12.6.1.$
of \cite{ds}.
In \cite{dav}
it
has been used for
showing
the equivalence
of
uniform
convergence
in 
strong
operator 
topology
of
a one-parameter
semigroup
depending on a
parameter
and
the
convergence
in 
strong
operator 
topology
of the 
resolvents
of the
corresponding
generators,
Theorem
$3.17.$.
\par
Notice that
if
$
\zeta
\doteqdot
-i\lambda
$
and
$
Q
\doteqdot
iT
$,
then the equality
\eqref{21440303}
turns into
$$
(Q+\zeta\un)^{-1}
=
\int_{0}^{\infty}
e^{-t\zeta}
e^{-Q t}
\,d\,t,
$$
which 
is referred 
in 
$IX.1.3.$
of
\cite{kato}
as
the fact 
that
the resolvent
of $Q$
is the 
\emph{Laplace}
transform
of the semigroup
$e^{-Q t}$.
Applications
of this formula
to perturbation theory
are
in
$IX.2.$
of
\cite{kato}.
\end{remark}
\chapter{
Extension
theorem.
The
case
of
the
topology
$\sigma(B(G),\n)$
}
\label{SECWEAKINT}
\section
{
Introduction 
}
Let 
$R$ be an
unbounded 
scalar type spectral operator
$R$ in a complex Banach space
$G$
and
$E$
its
resolution of identity.
The main results
of this chapter and of the work
are of two types.
 
The results
of the first type
are
Extension
Theorems
for integration
with respect to the 
$\sigma(B(G),\n)-$topology,
when
$\n$ is an
$E-$appropriate set:
Theorems
\ref{18051958ta}
and
when
$\n$ is an
$E-$appropriate set
with the duality
property:
Corollary
\ref{17070917TA}.
 
As an application we will prove,
by using
\eqref{13061209A},
the Extension theorems for
the integration with respect
to the 
sigma-weak topology:
Corollary
\ref{12121201}
and
Corollary
\ref{17070917pd},
and
for
integration
with respect to the
weak operator topology:
Corollary
\ref{18051958},
and
Corollary
\ref{17070917}.
 
The 
results of the 
second type
are
Newton-Leibnitz
formulas
for
integration with respect to the
$\sigma(B(G),\n)-$topology,
when
$\n$ is an
$E-$appropriate set
with the duality
property:
Corollary
\ref{20051321ta}
and
Corollary
\ref{20051321taLOC};
for 
integration
with respect to the
sigma-weak topology:
Corollary
\ref{20051321pd};
for 
integration
with respect to the
weak operator topology:
Corollary
\ref{19371401}
 
For obtaining
the Extension
Theorem
\ref{18051958ta}
we need
to introduce
the concept
of 
$E-$appropriate
set,
Definition
\ref{13481501},
which
allows
us to
establish
two
important
properties
for the proof
of
Theorem
\ref{18051958ta},
namely
the ``Commutation''
property,
Theorem
\ref{18051509ta},
and
the ``Restriction''
property,
Theorem
\ref{14050121}.
 
Finally
for obtaining
Corollary
\ref{17070917TA}
and
the
Newton-Leibnitz
formula
in
Corollary
\ref{20051321ta}
we have to 
introduce
the
concept
of
an
$E-$appropriate
set $\n$
with the duality
property,
Definition
\ref{13481501},
which
allows
us to
establish
conditions
ensuring that
a
map is
scalarly essentially
$(\mu,B(G))-$integrable
with respect
to
the 
$\sigma(B(G),\n)-$topology,
Theorem
\ref{19052218ta}.
Similar results
for the weak operator topology
are contained
in
Theorem
\ref{19052218}
and
Corollary
\ref{19052220}.
\section{
Existence of the
weak-integral
with respect to the 
$\sigma(B(G),\n)-$
topology
}
\label{SecWEAKINT}

 
In this section
we shall obtain 
a
general
result,
Theorem
\ref{19052218ta}
about 
conditions
ensuring that
a
map is
scalarly essentially
$(\mu,B(G))-$integrable
with respect
to
the 
$\sigma(B(G),\n)-$topology,
where
$\n$
is a 
suitable 
subset
of $B(G)^{*}$.
\begin{notationsE}
\label{16051536}
\begin{normalfont}
Let 
$\K\in\{\R,\C\}$,
$Z$
a linear space over
$\K$
and
$\tau$
a locally convex topology on $Z$,
then
we 
indicate with
$\lr{Z}{\tau}$
the associated
locally convex space
over
$\K$.
We denote
with
$LCS(\K)$
the
class of all the 
locally convex
spaces over
$\K$ 
and
for
any
$\lr{Z}{\tau}\in LCS(\K)$
we set
$\lr{Z}{\tau}^{*}$
for
its topological dual,
that is
the 
$\K-$linear space
of
all
$\K-$linear
continuous
functionals
on $Z$.
 
Let
$Y$ 
be
a linear space over $\K$
and
$U$
a
subspace
of
$Hom(Y,\K)$,
then
we indicate
with
the symbol
$
\sigma(Y,U)
$ 
the weakest 
(locally convex) 
topology
on $Y$ such that
$
U
\subseteq
\lr{Y}{\sigma(Y,U)}^{*}
$,
Def.
$2$,
$II.42$
of 
\cite{BourTVS},
which
coincides
with
the 
locally convex topology
on 
$Y$
generated
by the set of seminorms
$\Gamma(U)$
associated to
$U$
where
$
\Gamma(U)
\doteqdot
\{
q_{\phi}:Y\ni y\mapsto 
|\phi(y)|
\mid
\phi\in U
\}
$.
 
It is not hard
to see that
$
\sigma(Y,U)
$
is
the topology generated
by the set of seminorms
$\Gamma(S)$
for any
$S$
such that
$U=\mathfrak{L}_{\K}(S)$,
where
$\mathfrak{L}_{\K}(S)$
is
the $\K-$linear space generated by 
the set $S$.
 
By
Proposition $2$,
$II.43$
of \cite{BourTVS},
$\sigma(Y,U)$
is a
Hausdorff
topology
if and only if
$U$ separates the points of $Y$,
i.e.
\begin{equation}
\label{13481301}
(\forall T\in Y)
(T\neq\ze
\Rightarrow
(\exists\,\phi\in U)
(\phi(T)\neq 0)).
\end{equation}
Also
by 
Proposition $3$,
$II.43$
of \cite{BourTVS}
\begin{equation}
\label{15361301}
\lr{Y}{\sigma(Y,U)}^{*}
=
U.
\end{equation}
 
Let
$X$
be
a locally compact space
and
$\mu$
a 
$\K-$
Radon
measure
on
$X$,
Definition
$2$,
$\S1$, $n^{\circ}3$,
Ch. $3$, 
of \cite{IntBourb}
where it
is called just
measure.
We denote
with
$|\mu|$
the 
total variation
of
$\mu$,
$\S1$, $n^{\circ}6$,
Ch. $3$, 
of \cite{IntBourb},
and
with
$\int^{*}$
the 
upper integral
with respect to a positive measure,
as for example $|\mu|$,
Definition
$1$,
$\S1$, $n^{\circ}1$,
Ch. $4$, 
of \cite{IntBourb},
With
$\int^{\bullet}$
we denote
the essential upper integral
with respect to a positive measure,
Definition
$1$,
$\S1$, $n^{\circ}1$,
Ch. $5$, 
of \cite{IntBourb}.
\footnote{
In general
$
\int^{\bullet}
\leq
\int^{*}
$,
however
if $X$
is $\sigma-$compact,
in particular compact,
then
$
\int^{\bullet}
=
\int^{*}
$.
}
We
readdress
for
the definition
of
essentially 
$\mu-$integrable 
map
$
f:X\to\K
$,
to
Ch. $5$, 
$\S1$, $n^{\circ}3$,
of \cite{IntBourb}.
 
Let
$
\lr{Y}{\tau}
\in LCS(\K)
$
of
Hausdorff
then
$
f:X\to\lr{Y}{\tau}
$
is
\emph{
scalarly essentially 
$\mu-$integrable 
}
or 
equivalently
$
f:X\to Y
$
is
\emph{
scalarly 
essentially 
$\mu-$integrable 
with
respect to the 
measure
$\mu$
and
with respect
to
the
$
\tau-
$
topology
on $Y$
}
if
for
all
$\omega\in\lr{Y}{\tau}^{*}$
the 
map
$
\omega\circ f:X\to\K
$
is
essentially 
$\mu-$integrable, 
so
we can define
its
\emph{integral}
as 
the following linear 
operator
$$
\lr{Y}{\tau}^{*}
\ni
\omega
\mapsto
\int
\omega(f(x))
\,d\mu(x)
\in\K.
$$
See
Ch. $6$, 
$\S1$, $n^{\circ}1$
for $\K=\R$,
and
for the extension to 
the case $\K=\C$
see
the end
of
$\S2$, $n^{\circ}10$,
of \cite{IntBourb}.
 
Notice that the previous definitions
depend only on
the 
dual space
$
\lr{Y}{\tau}^{*}
$,
hence 
both the concepts of 
scalar essential 
$\mu-$integrability
and
integral
will be invariant
if we replace
$\tau$ with
any other 
Hausdorff 
locally convex 
topology
$\tau_{2}$
on $Y$
compatible
with 
the duality
$(Y,\lr{Y}{\tau}^{*})$,
i.e.
such that
$
\lr{Y}{\tau}^{*}
=
\lr{Y}{\tau_{2}}^{*}
$.
 
Therefore
as a corollary
of the well-known
Mackey-Arens
Theorem,
see 
Theorem $1$,
$IV.2$
of
\cite{BourTVS}
or
Theorem $5$
\S $8.5.$
of 
\cite{jar},
fixed 
a locally convex space
$
\lr{Y}{\tau}
$
and denoted
by
$
\n\doteqdot\lr{Y}{\tau}^{*}
$
its topological dual,
we have that
scalar essential 
$\mu-$integrability
(respectively
integral)
is an
invariant 
property
(respectively functional)
under the variation
of any
Hausdorff
locally convex 
topology
$\tau_{1}$
on $Y$
such that
$$
\sigma(Y,\n)
\leq
\tau_{1}
\leq
\tau(Y,\n).
$$
Here
$a\leq b$
means
$a$ is
weaker than $b$
and
$\tau(Y,\n)$
is the Mackey topology
associated to the 
canonical duality
$(Y,\n)$.
 
Let
$
f:X\to\lr{Y}{\tau}
$
be
scalarly essentially 
$\mu-$integrable 
and
assume that
\begin{equation}
\label{16021301}
(\exists\,B\in Y)
(\forall\omega\in
\lr{Y}{\tau}^{*})
\left(
\omega
(B)
=
\int
\omega(f(x))
\,d\mu(x)
\right).
\end{equation}
Notice
that
by the Hahn-Banach theorem
$\lr{Y}{\tau}^{*}$ separates
the points of $Y$,
so
the element
$B$
is defined by this condition
uniquely.
In this case, by definition
\emph{
$f:X\to\lr{Y}{\tau}$
is scalarly essentially 
$(\mu,Y)-$integrable
}
(or
\emph{
$
f:X\to Y
$
is scalarly essentially 
$(\mu,Y)-$integrable
with respect to the $\tau-$
topology
on $Y$
})
and
its
\emph{weak-integral
with respect
to the measure
$\mu$
and with respect
to the 
$\tau-$topology},
or
briefly
its
\emph{weak-integral},
is defined
by
\begin{equation}
\label{16021301bis}
\int
f(x)
\,d\mu(x)
\doteqdot
B.
\end{equation}
In the work
we shall use
this integral
for the case
$
\lr{Y}{\tau}
\doteqdot
\lr{B(G)}{\sigma(B(G),\n)}
$,
where
$\n$
is a linear
subspace
of
$B(G)^{*}$
which
separates
the points
of $B(G)$.
Notice that
by
\eqref{15361301}
$\lr{B(G)}{\sigma(B(G),\n)}^{*}=\n$.
 
Let
$G$
be
a 
$\K-$normed space,
then
the
strong operator topology
$\tau_{st}(G)$
on
$B(G)$
is defined
to be the locally convex topology
generated
by the following
set
of seminorms
$
\left
\{
q_{v}:
B(G)
\ni 
A
\mapsto
\|A v\|_{G}
\mid
v
\in 
G
\right\}
$.
Hence
$\tau_{st}(G)$
is a
Hausdorff topology,
in fact
a base
of the neighbourhoods
of
$A\in B(G)$
is
the class
of the
sets
$
U_{\ov{v},\varepsilon}(A)
\doteqdot
\{
B\in B(G)
\mid
\sup_{k=1,...,n}
\|(A-B)\ov{v}_{k}\|_{G}
<
\varepsilon
\}
$,
with
$\ov{v}$
running
in
$
\bigcup_{n\in\N}G^{n}
$
and
$
\epsilon
$
in
$\R^{+}-\{0\}$.
So
$B\in
\ov{\{\ze\}}
$,
the closure
of
$\{\ze\}$
in the 
strong operator topology,
if and only if
$
\|Bv\|_{G}<\ep
$,
for all
$\ep\in\R^{+}-\{0\},
v\in G$,
that
is
$B=\ze$.
Hence
$
\ov{\{\ze\}}
=
\{\ze\}
$
and then
$\tau_{st}(G)$
is of Hausdorff.
By
Ch. $6$, $\S1$, $n^{\circ}3$,
of \cite{IntBourb}
\begin{equation}
\label{17051313}
\n_{st}(G)
\doteqdot
\lr{B(G)}{\tau_{st}(G)}^{*}
=
\mathfrak{L}_{\K}
(
\{
\psi_{(\phi,v)}
\mid
(\phi,v)\in G^{*}\times G
\}
).
\end{equation}
Here
$$
\psi_{(\phi,v)}:
B(G)
\ni
T
\mapsto
\phi(T v)
\in
\K.
$$
Here
if
$Z$
is a $\K-$linear space
and
$S\subseteq Z$
then
$
\mathfrak{L}_{\K}(S)
$
is the 
space
of
all
$\K-$linear combinations
of elements
in
$S$.
 
The first 
locally convex space
in
which are mainly interested,
is 
$
\lr{B(G)}{\sigma(B(G),\n_{st}(G))}
$,
for which
by 
\eqref{15361301}
we have
\begin{equation}
\label{15381301}
\lr{B(G)}{\sigma(B(G),\n_{st}(G))}^{*}
=
\n_{st}(G).
\end{equation}
 
Notice that
by what said 
$
\sigma(B(G),\n_{st}(G))
$
is the topology
on $B(G)$
generated
by the 
set of seminorms
associated to the set
$\{
\psi_{(\phi,v)}
\mid
(\phi,v)\in G^{*}\times G
\}
$,
hence
$\sigma(B(G),\n_{st}(G))$
is nothing but 
the usual weak operator topology on $B(G)$. 
 
Notice 
that
by
\eqref{13481301},
and
the Hahn-Banach theorem
applied to $G$
we have that
$\sigma(B(G),\n_{st}(G))$
is a 
topology
of Hausdorff.
 
Let $G$ be a complex Hilbert space.
We define
$$
\n_{pd}(G)
\doteqdot 
B(G)_{*}.
$$
Here
$B(G)_{*}$
is the 
``predual'' 
of the von Neumann algebra
$B(G)$, see for example
Definition
$2.4.17.$
of \cite{bra},
or
Definition
$2.13.$,
Ch. $2$
of \cite{tak},
So
every $\omega\in\n_{pd}(G)$ 
has the
following form,
see Proposition $2.4.6$
of \cite{bra}
or
statement 
$(ii.4)$
of
Theorem $2.6.$,
Ch. $2$
of 
\cite{tak}
\begin{equation}
\label{19121201}
\omega: 
B(G)\ni 
B
\mapsto
\sum_{n=0}^{\infty}
\lr{u_{n}}{B w_{n}}
\in
\C.
\end{equation}
Here
$
\{u_{n}\}_{n\in\N},
\{w_{n}\}_{n\in\N}
\subset
G
$
are
such that
$
\sum_{n=0}^{\infty}
\|u_{n}\|^{2}
<\infty
$
and
$
\sum_{n=0}^{\infty}
\|w_{n}\|^{2}
<\infty
$.
 
We say that
\emph{
$\omega$
is 
determined by 
$
\{u_{n}\}_{n\in\N},
\,
\{w_{n}\}_{n\in\N}
$
}
if
\eqref{19121201}
holds.
Notice that
$\omega$ is well-defined,
indeed for all
$B\in B(G)$
we have
$
\sum_{n=0}^{\infty}
|\lr{u_{n}}{B w_{n}}|^{2}
\leq
\|B\|^{2}
\left(
\sum_{n=0}^{\infty}
\|u_{n}\|^{2}
\right)
\left(
\sum_{n=0}^{\infty}
\|w_{n}\|^{2}
\right)
<
\infty
$,
hence there exists
$\omega(B)$
and
$\omega\in B(G)^{*}$,
so 
\begin{equation}
\label{19141201}
\n_{pd}(G)
\subseteq 
B(G)^{*}.
\end{equation}
 
The second
locally convex space
in
which are mainly interested
is 
$
\lr{B(G)}{\sigma(B(G),\n_{pd}(G))}
$,
for which
by 
\eqref{15361301}
we have
\begin{equation}
\label{15381301pd}
\lr{B(G)}{\sigma(B(G),\n_{pd}(G))}^{*}
=
\n_{pd}(G).
\end{equation}
 
By the fact that
every 
$\omega\in\n_{st}(G)$
is determined by the
$
\{u_{n}\}_{n=1}^{N},
\,
\{w_{n}\}_{n=1}^{N}
$,
for
some
$N\in\N$,
we have that
$\n_{st}(G)\subset\n_{pd}(G)$.
Hence
being
$\sigma(B(G),\n_{st}(G))$
a
topology
of Hausdorff
we can conclude 
by
\eqref{13481301}
that
it is so
also
the 
$\sigma(B(G),\n_{pd}(G))-$
topology.
 
Notice that 
by what said 
$\sigma(B(G),\n_{pd}(G))$
is the topology
on $B(G)$
generated
by the 
set of seminorms
associated to the set
$\n_{pd}(G)$,
hence is nothing but 
the usual 
sigma-weak operator topology
on $B(G)$,
see for example for its definition
Section $2.4.1$ of \cite{bra},
so often we shall refer
to it just as the 
sigma-weak operator topology
on $B(G)$.
 
We want just to remark
that as a corollary
of 
the beforementioned
invariance property
for the weak-integration,
when we change 
the topology
$\tau$ on $Y$
with any other 
Hausdorff topology
compatible with 
it,
we deduce 
by
\eqref{17051313}
that
$
f:X\to B(G)
$
is
scalarly 
essentially 
$(\mu,B(G))-$integrable
with
respect to the 
measure
$\mu$
and
with respect
to
the
$
\sigma(B(G),\n_{st}(G))
$
topology
on $B(G)$,
if and only if
it is so
with respect
to
the
strong topology
$\tau_{st}(G)$
on 
$B(G)$,
and in this case
their weak-integrals
coincide.
Let $\A$
be a $\K-$Banach algebra 
then for all
$
A,B\in\A
$
set
$[A,B]\doteqdot AB-BA$,
while
the map
$
\RM:\A\to B(\A)
$
and
$
\mathcal{L}:
\A\to B(\A)
$,
have been
defined
in
\eqref{16051546int}.
Let
$G$
be a
$\K-$Banach space 
and
$
\n\subseteq 
B(G)^{*}
$
a linear subspace of the
normed space
$B(G)^{*}$,
then
we introduce
the following
notations
\begin{equation*}
\begin{cases}
\n^{*}
\subseteq 
B(G)
\overset{def}{\Leftrightarrow}
(\exists\,Y_{0}\subseteq B(G))
(\n^{*}
=
\{
\hat{A}
\up
\n
\mid
A\in Y_{0}
\});
\\
\n^{*}\on B(G)
\overset{def}{\Leftrightarrow}
\\
(\exists\,Y_{0}\subseteq B(G))
(\forall\phi\in\n^{*})
(\exists\,A\in Y_{0})
((\phi=\hat{A}\up\n)
\wedge
(\|\phi\|_{\n^{*}}=\|A\|_{B(G)})).
\end{cases}
\end{equation*}
Here
$
\left(\hat{\cdot}\right):
B(G)
\to
(B(G)^{*})^{*}
$
is the canonical
isometric embedding
of $B(G)$
into
its bidual.
 
By
statement 
$(iii)$
of
Theorem $2.6.$,
Ch. $2$
of 
\cite{tak},
or
Proposition $2.4.18$
of \cite{bra}
\begin{equation}
\label{18591801TH}
\n_{pd}(G)^{*}
\on
B(G).
\end{equation}
 
Let
$\mb{H}:\B_{Y}\to\Pr(G)$
be
a spectral measure in $G$
on $\B_{Y}$
then 
we continue
to follow the notation
$$
(\forall\sigma\in\B(\C))
(G_{\sigma}
\doteqdot
\mb{H}(\sigma)G),
$$
without
expressing
the dependence
on
$\mb{H}$
everywhere
it does not cause
confusion.
 
\emph{
In this Chapter
we assume
to be fixed 
a 
complex Banach space
$G$,
a
locally compact space $X$
a 
complex
Radon
measure 
$\mu$
on $X$,
a possibly
unbounded
scalar type
spectral operator
$R$
with
spectrum
$\sigma(R)$
and
resolution
of the
the identity
$E$}.
 
For each
map
$
f:U\subset\C\to\C
$
we denote
by
$\w{f}$
the $\ze-$extension
of $f$ to $\C$.
 
Finally
we shall denote
with
$
\F{ess}{}(X;\mu)
$
the 
seminormed space,
with the
seminorm
$
\|\cdot\|_{\F{ess}{}(X;\mu)}
$,
of all 
maps
$H:X\to\C$
such that
$$
\|
H
\|_{\F{ess}{}(X;\mu)}
\doteqdot
\int^{\bullet} 
|H(x)|
\,\,d
|\mu|(x)
<\infty.
$$
By
$\mu-l.a.e.(X)$
we shall mean
``locally almost everywhere
on
$X$
with respect to 
$\mu$''.
Moreover
if 
$
f:X_{0}\to\C
$
is a
map
defined
$\mu-l.a.e.(X)$,
then
we convene
to say
that
$f\in\F{ess}{}(X;\mu)$
if 
there exists
a 
map
$
F:X\to\C
$
such that
$F\up X_{0}=f$
and
$F\in\F{ess}{}(X;\mu)$.
In such a case
we set
\begin{equation}
\label{04250403}
\|
f
\|_{\F{ess}{}(X;\mu)}
\doteqdot
\|
F
\|_{\F{ess}{}(X;\mu)}.
\end{equation}
\eqref{04250403}
is well-defined
since
the definition
is independent
of which
map
$F\in\F{ess}{}(X;\mu)$
extends $f$,
as an application
of statement
$(a)$
of Proposition
$1$, $n^{\circ} 1$,
\S $1$,
Ch. $V$ 
of
\cite{IntBourb}.
Moreover
let
$\lr{Y}{\tau}$
be a locally convex space
and
$
f:X_{0}\to Y
$
a
map
defined
$\mu-l.a.e.(X)$,
then
we 
for brevity
say
that
the map
$
f:X\to\lr{Y}{\tau}
$
is
scalarly essentially
$(\mu,Y)-$integrable
if 
there exists
a 
map
$
F:X\to Y
$
such that
$F\up X_{0}=f$
and
$
F:X\to\lr{Y}{\tau}
$
is
scalarly essentially
$(\mu,Y)-$integrable.
In this case
we define
\begin{equation}
\label{20171902}
\int
f(x)
\,d\mu(x)
\doteqdot
\int
F(x)
\,d\mu(x).
\end{equation}
This does not depend on the choice of 
$F$.
Indeed,
for 
$k\in\{1,2\}$
let
$
F_{k}:X\to Y
$
be
such that
$F_{k}\up X_{0}=f$
and
$
F_{k}:X\to\lr{Y}{\tau}
$
be
scalarly essentially
$(\mu,Y)-$integrable,
then for all
$
\omega\in\lr{Y}{\tau}^{*},
k\in\{1,2\}
$
$$
\omega
\left(
\int 
F_{k}(x)
\,d\mu(x)
\right)
=
\int 
\omega
(F_{k}(x))
\,d\mu(x)
=
\int 
\chi_{X_{0}}(x)
\omega
(F_{k}(x))
\,d\mu(x).
$$
Next
for all
$x\in X$,
$\chi_{X_{0}}(x)
\omega(F_{1}(x))
=
\chi_{X_{0}}(x)
\omega
(F_{2}(x))$,
so 
for all
$\omega\in\lr{Y}{\tau}^{*}$
$$
\omega
\left(
\int 
F_{1}(x)
\,d\mu(x)
\right)
=
\omega
\left(
\int 
F_{2}(x)
\,d\mu(x)
\right).
$$
Then
\eqref{13481301}
yields
$
\int 
F_{1}(x)
\,d\mu(x)
=
\int 
F_{2}(x)
\,d\mu(x)
$.
\end{normalfont}
\end{notationsE}
Now we will show 
some
result
about which
functions
are
scalarly 
essentially
$(\mu,B(G))-$integrable
with respect
to 
the 
$\sigma(B(G),\n)-$topology.
Here
$\n\subseteq B(G)^{*}$,
such that
separates the points of
$B(G)$
and 
$
\n^{*}\subseteq B(G)
$.
Then we apply
these results 
to the case
when $G$ is a Hilbert space
and  
$\n=\n_{pd}(G)$.
\begin{theorem}
\label{19052218ta}
Let 
$G$ be a complex Banach space,
a subspace
$\n\subseteq B(G)^{*}$
be
such that
$\n$
separates the points of 
$B(G)$
and
$$
\n^{*}
\subseteq B(G).
$$
 
Let
$
F:X\to B(G)
$
be
a map
such
that for all
$
\omega
\in 
\n
$
the map
$
\omega\circ F:X\to\C
$
is 
$\mu-$measurable
and
\begin{equation}
\label{20411102}
(X\ni x
\mapsto
\|F(x)\|_{B(G)})
\in
\F{ess}{}(X;\mu).
\end{equation}
 
Then
the map
$
F:X\to
\lr
{B(G)}{\sigma(B(G),\n)}
$
is scalarly essentially
$(\mu,B(G))-$
integrable,
if in addition
$
\n^{*}
\on
B(G)
$
then
its
weak-integral
is such that
\begin{equation}
\label{16170702}
\left
\|
\int
F(x)
\,d\mu(x)
\right
\|_{B(G)}
\leq
\int^{\bullet}
\|F(x)\|_{B(G)}
\,d|\mu|(x).
\end{equation}
\end{theorem}
\begin{proof}
For all 
$\omega\in\n$
we have
$
|\omega(F(x))|
\leq
\|\omega\|
\|F(x)\|_{B(G)}
$,
hence for all
$\omega\in\n$
\begin{equation}
\label{25051022ta}
\int^{\bullet}
|\omega(F(x))|\,d|\mu|(x)
\leq
\|\omega\|
\int^{\bullet}
\|F(x)\|_{B(G)}
\,d|\mu|(x).
\end{equation}
 
Moreover
the map
$
\omega\circ F
$
is 
$\mu-$measurable
by hypothesis,
therefore
by 
\eqref{25051022ta}
and
Proposition
$9$, 
\S $1$,
$n^{\circ} 3$,
Ch.
$5$
of
\cite{IntBourb} 
we have that
$
\omega\circ F
$
is 
essentially
$\mu-$integrable.
 
By this fact we can define the following
map
$$
\Psi:
\n
\ni
\omega
\mapsto
\int
\omega(F(x))
\,d\mu(x)
\in\C
$$
which is linear.
Moreover
for 
any 
essentially
$\mu-$integrable
map
$H:X\to\C$
\begin{equation}
\label{12232102}
\left|
\int
H(x)
\,d\mu(x)
\right|
\leq
\int^{\bullet}
|H(x)|
\,
d|\mu|(x),
\end{equation}
hence
by
\eqref{25051022ta}
\begin{equation}
\label{19421201}
\Psi\in\n^{*}.
\end{equation}
Finally
by 
the duality property
$\n^{*}\subseteq B(G)$
in hypothesis
the statement
follows
by 
\eqref{19421201}
and
\eqref{25051022ta}.
\end{proof}
\begin{remark}
\label{19052218pd}
Let $G$ be a
complex Hilbert
space, then the statement
of Theorem \ref{19052218ta}
holds
if we set
$\n\doteqdot\n_{pd}(G)$.
Indeed 
we have the duality property
\eqref{18591801TH}.
\end{remark}
Now we give 
similar results
for
$\n=\n_{st}(G)$.
\begin{lemma}
\label{19052225}
Let 
$G$
be reflexive,
that is 
$(G^{*})^{*}$
is isometric
to
$G$
through
the 
natural 
injective
embedding
of
any normed space
into 
its bidual.
In addition
let
$B:G^{*}\times G\to\C$
be
a bounded bilinear form,
that
is
$$
(\exists\,C>0)
(\forall(\phi,v)\in G^{*}\times G)
(|B(\phi,v)|
\leq
C
\|\phi\|_{G^{*}}
\|v\|_{G}).
$$
Then
$$
(\exists\,! L\in B(G))
(\forall\phi\in G^{*})
(\forall v\in G)
(
B(\phi,v)
=
\phi(L(v))
$$
and
$
\|L\|_{B(G)}
\leq
\|B\|
$,
where
$
\|B\|
\doteqdot
\sup_{\{(\phi,v)\mid\|\phi\|_{G^{*}},
\|v\|_{G}\leq 1\}}
|B(\phi,v)|
$.
\end{lemma}
\begin{proof} For all
$v\in G $ 
let  
$T(v):G^{*}\ni\phi\mapsto B(\phi,v)\in\C$
so 
$
T(v)\in
(G^{*})^{*}
$
such that
$
\|T(v)\|_{(G^{*})^{*}}
\leq
\|B\|
\cdot
\|v\|_{G}
$.
$G$ is reflexive, hence 
$
(\forall v\in G)
(\exists\,! L(v)\in G)
(\forall\phi\in G^{*})
(\phi(L(v))
=T(v)(\phi))
$,
in addition
$
\|L(v)\|_{G}
=
\|T(v)\|_{(G^{*})^{*}}
\leq
\|B\|
\cdot
\|v\|_{G}
$.
$L$
is
linear
by the linearity
of
$T$
and by 
the fact that
$G^{*}$
separates the points
of $G$
by
the Hahn-Banach theorem.
Thus
$L$ 
is 
linear
and
bounded
and
$
\|L\|_{B(G)}
\leq
\|B\|
$.
This implies
the existence of $L$.
Let now
$L'\in B(G)$
be
another operator
with the same property,
so for all
$\phi\in G^{*},
v\in G$,
$\phi(L(v))
=
\phi(L'(v))$,
thus
by the Hahn-Banach 
theorem 
for all $v\in G$
$L(v)=L'(v)$,
which shows the uniqueness.
\end{proof}
\begin{theorem}
\label{19052218}
Let 
$G$
be 
reflexive,
$
F:X\to B(G)
$
be
a map
such
that
for all
$
(\phi,v)
\in 
G^{*}\times G
$
the map
$
X
\ni
x
\mapsto
\phi(F(x)v)
\in
\C
$
is 
$\mu-$measurable,
finally
assume that
\eqref{20411102}
holds.
Then
the map
$
F:X\to
\lr
{B(G)}{\sigma(B(G),\n_{st}(G))}
$
is scalarly essentially
$(\mu,B(G))-$integrable
and
its
weak-integral
satisfies
\eqref{16170702}.
\end{theorem}
\begin{proof}
We have for all
$\phi\in G^{*},
v\in G,x\in X$
that
$
|\phi(F(x)v)|
\leq
\|\phi\|
\|v\|
\|F(x)\|_{B(G)}
$,
hence
\begin{equation}
\label{25051022}
\int^{\bullet}
|\phi(F(x)v)|
\,
d|\mu|(x)
\leq
\|\phi\|
\|v\|
\int^{\bullet}
\|F(x)\|_{B(G)}
\,
d|\mu|(x)
\end{equation}
 
Furthermore
the map
$
X
\ni
x
\mapsto
\phi(F(x)v)
$
is
$\mu-$measurable
by hypothesis,
therefore
by 
\eqref{25051022}
and
Proposition
$9$, 
\S $1$,
$n^{\circ} 3$,
Ch.
$5$
of
\cite{IntBourb} 
we have that
$
X
\ni
x
\mapsto
\phi(F(x)v)
$
is 
essentially
$\mu-$integrable.
 
By this fact
we can define the following map
$$
B:G^{*}\times G
\ni
(\phi,v)
\mapsto
\int
\phi(F(x)v)
\,d\mu(x)
\in \C,
$$
which is bilinear.
So
by 
\eqref{12232102}
and
\eqref{25051022}
$$
|B(\phi,v)|
\leq
\|\phi\|
\|v\|
\int^{\bullet}
\|F(x)\|_{B(G)}
\,
d|\mu|(x).
$$
Hence
$B$ is a bounded bilinear form
whose norm 
is
such that
$
\|B\|
\leq
\int^{\bullet}
\|F(x)\|_{B(G)}
\,
d|\mu|(x)
$,
then the statement
by 
Lemma 
\ref{19052225}.
\end{proof}
\begin{corollary}
\label{19052220}
Let 
$G$
be 
reflexive,
$
F:X\to B(G)
$
a map
$
\sigma(B(G),\n_{st}(G))-
$
continuous,
i.e.
for all
$
(\phi,v)\in G^{*}\times G
$
the map
$
X
\ni
x
\mapsto
\phi(F(x)v)
\in
\C
$
is continuous,
finally
assume that
\eqref{20411102}
holds.
 
Then
the map
$
F:X
\to
\lr{B(G)}{\sigma(B(G),\n_{st}(G))}
$
is scalarly essentially
$(\mu,B(G))-$integrable
and
its
weak-integral
satisfies
\eqref{16170702}.
\end{corollary}
\begin{proof}
By definition of $\mu-$measurability
we have that
the continuity condition implies
that for all
$
(\phi,v)
\in 
G^{*}\times G
$
the map
$
X
\ni
x
\mapsto
\phi(F(x)v)
\in
\C
$
is 
$\mu-$measurable,
hence the statement
by Theorem
\ref{19052218}.
\end{proof}
\section{
Commutation
and
restriction
properties
}
\label{COMMURESTR}
Let 
$\mb{H}:\B_{Y}\to\Pr(G)$
be a spectral measure 
in $G$
on
$\B_{Y}$,
then
in the sequel we shall 
introduce a special class
of subspaces of $B(G)^{*}$,
the class of all
``$\mb{H}-$appropriate sets'',
which allows one
to show 
two
important 
properties
for
proving
the
main 
Extension
Theorem
\ref{18051958ta}.
These are 
\begin{enumerate}
\item
``Commutation'' 
property:
Theorem
\ref{18051509ta},
for a
general 
$E-$appropriate set
$\n$,
and
Corollary
\ref{18051509}
for
$\n=\n_{pd}(G)$
or
$\n=\n_{st}(G)$;
\item
``Restriction''
property:
Theorem
\ref{14050121}
for a
general 
$E-$appropriate set
$\n$.
\end{enumerate}
\begin{lemma}
\label{15051045}
Let 
$A\in B(G)$ 
such that 
$
AR
\subseteq
RA
$
and
$f\in Bor(\sigma(R))$.
Then 
$$
A f(R)
\subseteq
f(R) A.
$$
\end{lemma}
\begin{proof}
By Corollary $18.2.4.$ of \cite{ds} 
\begin{equation}
\label{16051542}
(\forall\sigma\in\B(\C))([A,E(\sigma)]=\ze).
\end{equation}
By 
\eqref{10501634}
for all
$T\in B(G))$,
$\RM(T),
\mathcal{L}(T)
\in 
B(B(G)))
$,
so
by using the notations in Preliminaries
\ref{II01041232},
we have for all
$n\in\N$
\begin{alignat}{2}
\label{16051616}
\I{\C}{E}(f_{n})
A
&
=
\left(
\mathcal{L}(A)
\circ
\I{\C}{E}
\right)
(f_{n})
\notag\\
&
=
\I{\C}{\mathcal{L}(A)\circ E}(f_{n})
&
\text{
by 
\eqref{20522502},
$
\mathcal{L}(A)
\in 
B(B(G))
$
}
\notag\\
&
=
\I{\C}{\RM(A)\circ E}(f_{n})
&
\text{
by 
\eqref{16051542}
}
\notag\\
&
=
\left(
\RM(A)
\circ
\I{\C}{E}
\right)
(f_{n})
&
\text{
by 
\eqref{20522502},
$
\RM(A)
\in 
B(B(G))
$
}
\notag\\
&
=
A\,
\I{\C}{E}(f_{n}).
\end{alignat}
Let $x\in Dom(f(R))$ then by
\eqref{II01041439},
the fact that
$A\in B(G)$
and 
\eqref{16051616}
$$
A
f(R)
x
=
\lim_{n\to\infty}
\I{\C}{E}
(f_{n})
Ax.
$$
Hence
\eqref{II01041439}
implies
$
Ax\in Dom(f(R))
$
and
$$
f(R)A
x
=
\lim_{n\to\infty}
\I{\C}{E}(f_{n})Ax
=
A f(R)x.
$$
\end{proof}
\begin{lemma}
\label{12051045}
Let 
$
\n
\subseteq
B(G)^{*}
$
be
such that
$
\sigma(B(G),\n)
$
is a
Hausdorff
topology,
$A\in B(G)$,
and
the
map
$
X
\ni 
x
\mapsto
f_{x}
\in
Bor(\sigma(R))
$
be
such that
$
\w{f}_{x}\in\Lf{E}{\infty}(\sigma(R))
$
\,
$\mu-l.a.e.(X)$.
Assume that
\begin{enumerate}
\item
the map
$
X\ni
x
\mapsto
f_{x}(R)
\in
\lr{B(G)}{\sigma(B(G),\n)}
$
is
scalarly essentially
$(\mu,B(G))-$integrable;
\item
$
\phi\circ\RM(A)
\in\n
$
and
$
\phi\circ
\mathcal{L}(A)
\in\n
$,
for all
$
\phi\in\n
$;
\item
$
A R
\subseteq 
R A
$.
\end{enumerate}
Then
$$
\left[
\int
f_{x}(R)
\,d\mu(x),
\,
A
\right]
=
\ze.
$$
\end{lemma}
\begin{proof}
By the
hypothesis
$
\w{f}_{x}
\in
\Lf{E}{\infty}(\sigma(R))
$,\,
$
\mu-l.a.e.(X)
$
and
statement
$(c)$
of
Theorem
$18.2.11.$
of
\cite{ds}
applied to
the scalar type
spectral operator
$
R
$,
we have
$f_{x}(R)
\in
B(G),
\,
\mu-l.a.e.(X)
$.
Let
us
set
$
X_{0}
\doteqdot
\{
x\in X
\mid
f_{x}(R)
\in B(G)
\}
$.
By the hypothesis 
$(1)$
we deduce that there is
$F:X\to B(G)$
such that
\begin{itemize}
\item
$(\forall x\in X_{0})(F(x)=f_{x}(R))$;
\item
$
F:X\to\lr{B(G)}{\sigma(B(G),\n)}
$
is scalarly essentially 
$(\mu,B(G))-$integrable.
\end{itemize}
Thus
by definition
\begin{equation}
\label{20541902}
\int
f_{x}(R)
\,d\mu(x)
\doteqdot
\int
F(x)
\,d\mu(x)
\end{equation}
Notice
that
for all
$x\in X,
\phi\in\n$
\begin{equation}
\label{20511902}
\chi_{X_{0}}(x)\,
\phi
\circ
\mathcal{L}(A)
(F(x))
=
\chi_{X_{0}}(x)\,
\phi
\circ
\RM(A)
(F(x)),
\end{equation}
since
by Lemma
\ref{15051045}
for all
$
x\in X_{0}
$
$$
F(x)A
=
f_{x}(R)A
=
Af_{x}(R)
=
AF(x).
$$
Moreover for all
$\phi\in\n$
\begin{equation}
\label{17311903}
\begin{cases}
\int
\phi\circ\mathcal{L}(A)
\left(
F(x)
\right)
\,d\mu(x)
=
\int
\chi_{X_{0}}(x)\,
\phi\circ\mathcal{L}(A)
\left(
F(x)
\right)
\,d\mu(x),
\\
\int
\phi\circ\RM(A)
\left(
F(x)
\right)
\,d\mu(x)
=
\int
\chi_{X_{0}}(x)\,
\phi\circ\RM(A)
\left(
F(x)
\right)
\,d\mu(x).
\end{cases}
\end{equation}
Indeed
$
\phi\circ\mathcal{L}(A)
\in
\n
$
hence
$
X\ni x
\mapsto
\phi\circ\mathcal{L}(A)
\left(F(x)\right)
$
is 
essentially
$\mu-$integrable
so
by
Proposition
$9$
$n^{\circ} 3$
\S $1$
Ch $5$
of
\cite{IntBourb} 
$$
\int^{\bullet}
\big|
\chi_{X_{0}}(x)\,
\phi\circ\mathcal{L}(A)
\left(F(x)\right)
\big|
d\,|\mu|(x)
\leq
\int^{\bullet}
\big|
\phi\circ\mathcal{L}(A)
\left(F(x)\right)
\big|
d\,|\mu|(x)
<
\infty.
$$
Furthermore
by
Proposition
$6$
$n^{\circ} 2$
\S $5$
Ch $4$
of
\cite{IntBourb} 
$
X\ni x
\mapsto
\chi_{X_{0}}(x)\,
\phi\circ\mathcal{L}(A)
\left(F(x)\right)
$
is
$\mu-$measurable.
Thus
by
Proposition
$9$
$n^{\circ} 3$
\S $1$
Ch $5$
of
\cite{IntBourb} 
the map
$
X\ni x
\mapsto
\chi_{X_{0}}(x)\,
\phi\circ\mathcal{L}(A)
\left(F(x)\right)
$
is
essentially
$\mu-$integrable
and
we obtain
the first
statement
of
\eqref{17311903}
by the fact that
two essentially
$\mu-$integrable
maps
that are
equal 
$\mu-l.a.e.(X)$
have the same integral.
In the same
way it is possible to show
also the second 
statement
of
\eqref{17311903}.
Therefore for all
$\phi\in\n$
\begin{alignat*}{2}
\phi
\left(
\int
f_{x}(R)
\,d\mu(x)
A
\right)
&
=
\phi
\circ
\mathcal{L}(A)
\left(
\int
f_{x}(R)
\,d\mu(x)
\right)
\\
&
=
\phi
\circ
\mathcal{L}(A)
\left(
\int
F(x)
\,d\mu(x)
\right)
&
\text{ by \eqref{20541902}}
\\
&
=
\int
\phi\circ\mathcal{L}(A)
\left(
F(x)
\right)
\,d\mu(x)
&
\text{ 
by 
$
\phi\circ\mathcal{L}(A)
\in
\n
$
}
\\
&
=
\int
\phi\circ\RM(A)
\left(
F(x)
\right)
\,d\mu(x)
&
\text{ 
by 
\eqref{17311903},
\eqref{20511902}
}
\\
&
=
\phi\circ\RM(A)
\left(
\int
F(x)
\,d\mu(x)
\right)
&
\text{ 
by 
$
\phi\circ\RM(A)
\in
\n
$
}
\\
&
=
\phi
\left(
A
\int
f_{x}(R)
\,d\mu(x)
\right).
&
\text{ by \eqref{20541902}}
\end{alignat*}
\hspace{12pt}
Then the statement
by
\eqref{13481301}
\end{proof}
\begin{remark}
\label{15051128}
By
definition
of $\n_{st}(G)$, 
see
\eqref{17051313},
the
hypothesis
$(2)$
of 
Lemma \ref{12051045}
holds for all
$A\in B(G)$
and
for
$\n=\n_{st}(G)$.
Moreover
$\sigma(B(G),\n_{st}(G))$
is a Hausdorff topology on
$B(G)$.
 
Let $G$ be a Hilbert space,
by \eqref{19121201}
we note that
for all 
$
A\in B(G)
$
we
have 
$
\omega
\circ
\mathcal{L}(A)
\in
\n_{pd}(G)
$,
and
$
\omega
\circ
\RM(A)
\in
\n_{pd}(G)
$,
indeed
if
$\omega$
is determined
by
$
\{u_{n}\}_{n\in\N},
\{w_{n}\}_{n\in\N}
$,
then
$
\omega
\circ
\mathcal{L}(A)
$,
(respectively
$\omega
\circ
\RM(A)
$),
is
determined by 
$
\{u_{n}\}_{n\in\N},
\{A w_{n}\}_{n\in\N}
$,
(respectively
$
\{A^{*}u_{n}\}_{n\in\N},
\{w_{n}\}_{n\in\N}
$).
Hence
the
hypothesis
$(2)$
of 
Lemma \ref{12051045}
holds for all
$A\in B(G)$
and
for
$\n=\n_{pd}(G)$.
Furthermore
$\sigma(B(G),\n_{pd}(G))$
is a Hausdorff topology on
$B(G)$.
\end{remark}
\begin{remark}
\label{15051129}
By 
Definition
$18.2.1$
of \cite{ds}
for all
$\sigma\in\B(\C)$,
$E(\sigma)R
\subseteq
R
E(\sigma)$,
thus
hypothesis
$(3)$
of 
Lemma \ref{12051045}
holds
for
$
A
\doteqdot
E(\sigma)
$.
\end{remark}
\begin{definition}
[
\textbf{
$\mb{H}-$appropriate set
}
]
\label{13481501}
Let 
$\mb{H}:\B_{Y}\to\Pr(G)$
be a spectral measure 
in $G$
on
$\B_{Y}$,
see Preliminaries
\ref{II01041232}.
Then 
we define
$\n$
to be
an
\textbf{
$\mb{H}-$appropriate set
},
if
\begin{enumerate}
\item
$
\n
\subseteq 
B(G)^{*}
$
linear subspace;
\item
$\n$ 
separates the points 
of $B(G)$,
namely
$$
(\forall T\in B(G))
(T\neq\ze
\Rightarrow
(\exists\,\phi\in\n)
(\phi(T)\neq 0));
$$
\item
for all
$\phi\in\n,
\sigma\in\B_{Y}$
\begin{equation}
\label{14111501}
\phi\circ
\RM(\mb{H}(\sigma))
\in\n
\quad
\phi\circ
\mathcal{L}(\mb{H}(\sigma))
\in\n.
\end{equation}
\end{enumerate}
Furthermore
$\n$
is
an
\textbf{
$\mb{H}-$appropriate set
with the
duality property
}
if
$\n$
is
an
$\mb{H}-$appropriate set
such that
\begin{equation*}
\n^{*}
\subseteq 
B(G).
\end{equation*}
Finally
$\n$
is
an
$\mb{H}-$appropriate set
with the
isometric
duality property
if
$\n$
is
an
$\mb{H}-$appropriate set
such that
\begin{equation*}
\n^{*}
\on
B(G).
\end{equation*}
\end{definition}
\begin{remark}
\label{14551501}
Some comments about the previous definition.
The separation property
is 
equivalent to require that
$
\sigma(B(G),\n)
$
is
a Hausdorff topology
on $B(G)$,
while
\eqref{14111501}
is equivalent
to require
that
for all
$
\sigma\in\B_{Y}
$
the maps
on $B(G)$,
$
\RM(\mb{H}(\sigma))
$
and
$\mathcal{L}(\mb{H}(\sigma))$
are continuous
with respect to the
$
\sigma(B(G),\n)-
$
topology.
Moreover
the duality property
$
\n^{*}
\subseteq 
B(G)
$
ensures
that 
suitable 
scalarly essentially
$\mu-$integrable
maps
with respect
to the 
$
\sigma(B(G),\n)-
$
topology,
are
$(\mu,B(G))-$integrable,
see Theorem
\ref{19052218ta}.
 
Finally
by 
Remark \ref{15051128}
$\n_{st}(G)$
and
$\n_{pd}(G)$,
in case in which $G$
is a Hilbert space,
are 
$\mb{H}-$appropriate
sets 
for any spectral measure
$\mb{H}$,
furthermore
by 
\eqref{18591801TH},
$\n_{pd}(G)$
is
an
$\mb{H}-$appropriate
set
with
the isometric
duality property.
\end{remark}
\begin{theorem}
[\textbf{Commutation 1}]
\label{18051509ta}
Let
$\n$
be an
$E-$appropriate set,
the map
$
X
\ni 
x
\mapsto
f_{x}
\in
Bor(\sigma(R))
$
be
such 
that
$
\w{f}_{x}\in\Lf{E}{\infty}(\sigma(R))
$
\,
$\mu-l.a.e.(X)$.
Assume
that
the map
$
X
\ni
x
\mapsto
f_{x}(R)
\in
\lr{B(G)}{\sigma(B(G),\n)}
$
is
scalarly essentially
$(\mu,B(G))-$integrable.
Then for all
$\sigma\in\B(\C)$
\begin{equation}
\label{18051301}
\left[
\int
f_{x}(R)
\,d\mu(x),
\,
E(\sigma)
\right]
=
\ze.
\end{equation}
\end{theorem}
\begin{proof}
$\n$
being
an
$E-$appropriate
set
ensures that
hypothesis
$(2)$ of Lemma
\ref{12051045}
is satisfied
for $A\doteqdot E(\sigma)$
for all $\sigma\in\B(\C)$,
so
the statement
by
Remark \ref{15051129}
and
Lemma \ref{12051045}.
\end{proof}
\begin{corollary}
[\textbf{Commutation 2}]
\label{18051509}
\eqref{18051301}
holds
if we 
replace 
$\n$
in
Theorem
\ref{18051509ta}
with
$\n_{st}(G)$
or
with
$\n_{pd}(G)$
and assume that
$G$ is a Hilbert
space.
\end{corollary}
\begin{proof}
By
Remark \ref{14551501}
and
Theorem
\ref{18051509ta}.
\end{proof}
Now we 
are going 
to present some 
results
necessary
for
showing
the 
Restriction
property
in
Theorem
\ref{14050121},
namely
that
the map
$
X\ni x
\mapsto
f_{x}(R_{\sigma}\up G_{\sigma})
\in
\lr{B(G_{\sigma})}
{\sigma(B(G_{\sigma}),\n_{\sigma})}
$
is
scalarly essentially
$(\mu,B(G_{\sigma}))-$integrable,
where
$\n$ is a $E-$appropriate set,
and
$\n_{\sigma}$
could be thought
as
the 
``restriction''
of 
$\n$
to $B(G_{\sigma})$
for all $\sigma\in\B(\C)$.
 
In particular
when 
$\n=\n_{st}(G)$,
respectively
$\n=\n_{pd}(G)$,
we can replace
$\n_{\sigma}$
with
$\n_{st}(G_{\sigma})$,
respectively
$\n_{pd}(G_{\sigma})$,
Proposition
\ref{15051146}.
\begin{lemma}
\label{14050203}
Let 
$\mb{H}:\B_{Y}\to\Pr(G)$
be a spectral measure 
in $G$
on
$\B_{Y}$,
see Preliminaries
\ref{II01041232}.
Then for all
$\sigma\in\B_{Y}$
$G=G_{\sigma}\bigoplus G_{\sigma'}$,
where
$\sigma'\doteqdot\complement\sigma$.
\end{lemma}
\begin{proof}
$
\mb{H}(\sigma)+\mb{H}(\sigma')
=\mb{H}(\sigma\cup\sigma')=\un
$
so
$
\mb{H}(\sigma')=\un-\mb{H}(\sigma)
$
and
$
\mb{H}(\sigma)\mb{H}(\sigma')
=
\mb{H}(\sigma')\mb{H}(\sigma)
=\ze
$.
Hence for all
$v\in G$,
$v=\mb{H}(\sigma)v+\mb{H}(\sigma')v$,
or
$
G
=
G_{\sigma}
+
G_{\sigma'}
$.
But
for any
$\delta\in\B_{Y}$
we have
$
G_{\delta}
=
\{
y\in G\mid
y=\mb{H}(\delta)y
\}
$
then
$
G_{\sigma}
\cap
G_{\sigma'}
=
\{
y\in G
\mid
y=
\mb{H}(\sigma)
\mb{H}(\sigma')
y
\}
=
\{\ze\}
$.
Thus
$
G_{\sigma}
+
G_{\sigma'}
=
G_{\sigma}
\bigoplus
G_{\sigma'}
$.
\end{proof}
\begin{definition}
\label{14050215}
Let 
$\mb{H}:\B_{Y}\to\Pr(G)$
be a spectral measure 
in $G$
on
$\B_{Y}$,
$
\sigma\in\B_{Y}
$
and
$
\sigma'
\doteqdot
\complement
\sigma
$.
Then
Lemma
\ref{14050203}
allows us
to define
the
operator
$
\xi_{\sigma}^{\mb{H}}:
B(G_{\sigma})
\to
B(G)
$,
such
that for all
$T_{\sigma}\in B(G_{\sigma})$
\begin{equation}
\label{18051432}
\xi_{\sigma}^{\mb{H}}
(T_{\sigma})
\doteqdot
T_{\sigma}
\oplus
\ze_{\sigma'}
\in
B(G).
\end{equation}
Whenever
it does
not cause
confusion
we shall
denote
$\xi_{\sigma}^{\mb{H}}$
simply
by
$\xi_{\sigma}$.
Here
$\ze_{\sigma'}\in B(G_{\sigma'})$
is the null element 
of
the space
$B(G_{\sigma'})$,
while
the direct sum of 
two
operators
$
T_{\sigma}
\in
B(G_{\sigma})
$
and
$
T_{\sigma'}
\in
B(G_{\sigma'})
$
is
the following
standard 
definition
$$
(
T_{\sigma}\oplus T_{\sigma'}
):
G_{\sigma}\bigoplus G_{\sigma'}
\ni
(v_{\sigma}\oplus v_{\sigma'})
\mapsto
T_{\sigma}
v_{\sigma}
\oplus 
T_{\sigma'}
v_{\sigma'}
\in
G_{\sigma}\bigoplus G_{\sigma'}.
$$
\end{definition}
\begin{lemma}
\label{13061030}
Let 
$\mb{H}:\B_{Y}\to\Pr(G)$
be a spectral measure 
in $G$
on
$\B_{Y}$,
then for all
$\forall
\sigma\in\B_{Y},
T_{\sigma}\in B(G_{\sigma})$
we have that
\begin{equation}
\label{13060952}
\xi_{\sigma}^{\mb{H}}
(T_{\sigma})
=
T_{\sigma}
\mb{H}(\sigma).
\end{equation}
Hence 
$\xi_{\sigma}^{\mb{H}}$
is well-defined,
injective,
$
\xi_{\sigma}^{\mb{H}}
\in 
B(B(G_{\sigma}),B(G))
$
and
$
\|
\xi_{\sigma}^{\mb{H}}
\|_{B(B(G_{\sigma}),B(G))}
\leq
\|
\mb{H}(\sigma)
\|_{B(G)}
$.
\end{lemma}
\begin{proof}
Let
$\sigma\in\B_{Y}$
then for all
$v\in G$
we have
$$
(T_{\sigma}
\oplus
\ze_{\sigma'})
v
=
(T_{\sigma}
\oplus
\ze_{\sigma'})
(\mb{H}(\sigma)v
\oplus
\mb{H}(\sigma')v)
=
(
T_{\sigma}
\mb{H}(\sigma)v
\oplus
\ze)
=
T_{\sigma}
\mb{H}(\sigma)
v,
$$
then the first part.
Let
$
T_{\sigma}
\in B(G_{\sigma})
$
such that
$
\xi_{\sigma}(T_{\sigma})
=
\ze
$,
then
$
T_{\sigma}
\mb{H}(\sigma)
=
\ze
$,
which implies
that
for all 
$v_{\sigma}\in G_{\sigma}$
we have
$
T_{\sigma}
v_{\sigma}
=
T_{\sigma}
\mb{H}(\sigma)
v_{\sigma}
=
\ze
$.
So
$
T_{\sigma}
=
\ze_{\sigma}
$.
Let
us
consider
$
\mb{H}(\sigma)
\in
B(G,G_{\sigma})
$,
and
$
T_{\sigma}
\in
B(G_{\sigma},G)
$,
so
$
T_{\sigma}
\mb{H}(\sigma)
\in
B(G)
$
and
$
\|
T_{\sigma}
\mb{H}(\sigma)
\|_{B(G)}
\leq
\|
T_{\sigma}
\|_{B(G_{\sigma},G)}
\cdot
\|
\mb{H}(\sigma)
\|_{B(G,G_{\sigma})}
=
\|
T_{\sigma}
\|_{B(G_{\sigma})}
\cdot
\|
\mb{H}(\sigma)
\|_{B(G)}
$.
\end{proof}
Notice that
by \eqref{13060952}
and the fact that
$
B(G_{\sigma})
$
is a Banach space,
it is possible
to show that
$
\xi_{\sigma}(B(G_{\sigma}))
$
is a Banach subspace
of
$B(G)$,
thus
$\xi_{\sigma}$
has a continuous
inverse.
\begin{remark}
\label{15460602}
Let 
$\mb{H}:\B_{Y}\to\Pr(G)$
be a spectral measure 
in $G$
on
$\B_{Y}$,
and 
$\sigma\in\B_{Y}$.
If we
consider
the product space
$
G_{\sigma}\times G_{\sigma'}
$
with the standard 
linearization
and
define
\begin{equation}
\label{18370802}
\begin{cases}
\|(x_{\sigma},x_{\sigma'})\|_{\oplus}
\doteqdot
\|x_{\sigma}+x_{\sigma'}\|_{G},
\\
I:
G_{\sigma}
\times
G_{\sigma'}
\ni
(x_{\sigma},x_{\sigma'})
\mapsto
x_{\sigma}+x_{\sigma'}
\in
G,
\end{cases}
\end{equation}
then
by
$
G
=
G_{\sigma}\bigoplus G_{\sigma'}
$,
see
Lemma
\ref{14050203},
the two
spaces
$
\lr{G_{\sigma}\times G_{\sigma'}}
{\|\cdot\|_{\oplus}}
$
and
$
\lr{G}{\|\cdot\|_{G}}
$
are 
isomorphics,
thus
isometric
and $I$ is an isometry between them.
It is not difficult
to see that
the topology
induced by the norm
$\|\cdot\|_{\oplus}$
is the product topology
on
$
G_{\sigma}
\times
G_{\sigma'}
$
\footnote{
Indeed 
let
$\sigma\in\B_{Y}$
such that
$\mb{H}(\sigma)\neq\ze$,
set
$
M\doteqdot\max\{
\|\mb{H}(\sigma)\|,
\|\mb{H}(\sigma')\|\}
$
and for all
$r>0$
define
$
B_{r}^{\oplus}(\ze)
\doteqdot
\{
(x_{\sigma},x_{\sigma'})
\in G_{\sigma}
\times
G_{\sigma'}
\mid
\|
(x_{\sigma},x_{\sigma'})
\|_{\oplus}
<
r
\}
$.
Thus for all
$\ep>0$
by setting
$
\eta
\doteqdot
\frac{\ep}{2}
$
we have
$
B_{\eta}(\ze_{\sigma})
\times
B_{\eta}(\ze_{\sigma}')
\subset
B_{\ep}^{\oplus}(\ze)
$,
while for all
$\ep_{1},\ep_{2}>0$
by setting
$
\zeta
\doteqdot
\frac{\min\{\ep_{1},\ep_{2}\}}{M}
$
we have 
$
B_{\zeta}^{\oplus}(\ze)
\subset
B_{\ep_{1}}(\ze_{\sigma})
\times
B_{\ep_{2}}(\ze_{\sigma}')
$.
},
which implies
the following property
that
in any case we prefer to show
directly.
\end{remark}
\begin{proposition}
\label{18340802}
Let 
$\mb{H}:\B_{Y}\to\Pr(G)$
be a spectral measure 
in $G$
on
$\B_{Y}$
and
assume
the notations
in \eqref{18370802}
and Definition
\ref{14050215}.
For all
$T_{\sigma}\in B(G_{\sigma})$
and
$T_{\sigma'}\in B(G_{\sigma'})$
set
$$
T_{\sigma}
\times
T_{\sigma'}:
G_{\sigma}
\times
G_{\sigma'}
\ni
(x_{\sigma},x_{\sigma'})
\mapsto
(T_{\sigma}x_{\sigma},T_{\sigma'}x_{\sigma'})
\in
G_{\sigma}
\times
G_{\sigma'}
$$
Then
\begin{equation}
\label{16570802}
\begin{cases}
T_{\sigma}\oplus T_{\sigma'}
=
I
(T_{\sigma}\times T_{\sigma'})
I^{-1}
=
T_{\sigma}\mb{H}(\sigma)
+
T_{\sigma'}\mb{H}(\sigma')
\in
B(G)\\
T_{\sigma}
\times
T_{\sigma'}
=
I^{-1}
(T_{\sigma}\mb{H}(\sigma)
+
T_{\sigma'}\mb{H}(\sigma'))
I
\in
B(G_{\sigma}\times G_{\sigma'}).
\end{cases}
\end{equation}
\end{proposition}
\begin{proof}
$
I(T_{\sigma}\times T_{\sigma'})I^{-1}
(x_{\sigma}\oplus x_{\sigma'})
=
I
(T_{\sigma}x_{\sigma},T_{\sigma'}x_{\sigma'})
=
T_{\sigma}x_{\sigma}\oplus T_{\sigma'}x_{\sigma'}
$,
for all
$x_{\sigma}\in G_{\sigma}$
and
$x_{\sigma'}\in G_{\sigma'}$,
so the first equality.
For all
$x\in G$
\begin{alignat*}{1}
I
(T_{\sigma}\times T_{\sigma'})
I^{-1}
(x)
&
=
I
(T_{\sigma}\times T_{\sigma'})
I^{-1}
(\mb{H}(\sigma)x+\mb{H}(\sigma')x)
\\
&
=
I
(T_{\sigma}\mb{H}(\sigma)x,
T_{\sigma'}\mb{H}(\sigma')x)
\\
&
=
T_{\sigma}\mb{H}(\sigma)x
+
T_{\sigma'}\mb{H}(\sigma')x.
\end{alignat*}
Then 
the second
equality.
The third equality
is by the second
and the fact that
$I$ is an isometry.
\end{proof}
Notice that
by the first statement 
in \eqref{16570802}
we obtain
\eqref{13060952}.
\begin{definition}
\label{10061745}
Let 
$\mb{H}:\B_{Y}\to\Pr(G)$
be a spectral measure 
in $G$
on
$\B_{Y}$
and
$
\n
\subseteq
B(G)^{*}
$.
We define
for all
$\sigma\in\B_{Y},
\psi
\in
\n$
\begin{equation}
\label{18051637}
\begin{cases}
\psi_{\sigma}^{\mb{H}}
\doteqdot
\psi
\circ
\xi_{\sigma}^{\mb{H}}
\in
B(G_{\sigma})^{*}
\\
\n_{\sigma}^{\mb{H}}
\doteqdot
\{
\psi_{\sigma}^{\mb{H}}
\mid
\psi
\in
\n
\},
\end{cases}
\end{equation}
where
$\xi_{\sigma}^{\mb{H}}$
has been defined
in \eqref{18051432}.
We shall
express
$\psi_{\sigma}^{\mb{H}}$
and
$\n_{\sigma}^{\mb{H}}$
simply 
by
$\psi_{\sigma}$
and
$\n_{\sigma}$
respectively,
whenever 
it does not cause confusion.
\end{definition}
\begin{proposition}
\label{13422102}
Let 
$\mb{H}:\B_{Y}\to\Pr(G)$
be a spectral measure 
in $G$
on
$\B_{Y}$,
$
\n
\subseteq
B(G)^{*}
$
such that
$\n$
separates the points
of $B(G)$
and
$\sigma\in\B_{Y}$.
Then
$\n_{\sigma}$
separates
the points
of
$B(G_{\sigma})$.
\end{proposition}
\begin{proof}
Let 
$
T_{\sigma}\in 
B(G_{\sigma})-\{\ze_{\sigma}\}
$,
by Lemma \ref{13061030}
$\xi_{\sigma}$ 
is injective
so 
$
\xi_{\sigma}(T_{\sigma}) 
\neq
\ze
$.
But
$\n$
separates the points
of $B(G)$,
so
there is
$\psi\in\n$
such that
$\psi(\xi_{\sigma}(T_{\sigma}))\neq 0$.
\end{proof}
\begin{theorem}
[\textbf{Restriction}]
\label{14050121}
Let
$\n$
be an
$E-$appropriate set,
the map
$
X
\ni 
x
\mapsto
f_{x}
\in
Bor(\sigma(R))
$
be
such that
$
\w{f}_{x}
\in
\Lf{E}{\infty}(\sigma(R))
$
\,
$\mu-l.a.e.(X)$
Assume
that
the map
$
X\ni
x
\mapsto
f_{x}(R)
\in
\lr{B(G)}{\sigma(B(G),\n)}
$
is
scalarly essentially
$(\mu,B(G))-$integrable.
Then for all
$\sigma\in\B(\C)$
the map
$
X\ni x
\mapsto
f_{x}(R_{\sigma}\up G_{\sigma})
\in
\lr{B(G_{\sigma})}
{\sigma(B(G_{\sigma}),\n_{\sigma})}
$
is
scalarly essentially
$(\mu,B(G_{\sigma}))-$integrable
and
\begin{equation}
\label{restrta}
\int
f_{x}(R_{\sigma}\up G_{\sigma})
\,d\mu(x)
=
\int
f_{x}(R)
\,d\mu(x)
\up
G_{\sigma}.
\end{equation}
\end{theorem}
\begin{proof}
Let
$\sigma\in\B(\C)$
then
\eqref{25052019}
implies
that for all
$x\in X$
the
operator
$
f_{x}(R_{\sigma}\up G_{\sigma})
$
is
well-defined.
By the
hypothesis
$
\w{f}_{x}
\in
\Lf{E}{\infty}(\sigma(R))
$,\,
$
\mu-l.a.e.(X)
$
and
statement
$(c)$
of
Theorem
$18.2.11.$
of
\cite{ds}
applied to
the scalar type
spectral operator
$
R
$,
we have
$f_{x}(R)
\in
B(G),
\,
\mu-l.a.e.(X)
$.
Let
us
set
$$
X_{0}
\doteqdot
\{
x\in X
\mid
f_{x}(R)
\in B(G)
\},
$$
thus
by
statement
$(2)$
of
Lemma
\ref{II31031834}
we
obtain
\begin{equation}
\label{15051542WEAK}
(\forall x\in X_{0})
(f_{x}(R_{\sigma}\up G_{\sigma})
\in
B(G_{\sigma})).
\end{equation}
Hence
$
f_{x}(R_{\sigma}\up G_{\sigma})
\in
B(G_{\sigma})
$,
$
\mu-l.a.e.(X)
$.
So 
by 
Proposition
\ref{13422102}
and
\eqref{13481301}
it
is well-defined
the statement
that
$
X\ni x
\mapsto
f_{x}(R_{\sigma}\up G_{\sigma})
\in
\lr{B(G_{\sigma})}
{\sigma(B(G_{\sigma}),\n_{\sigma})}
$
is
scalarly essentially
$(\mu,B(G_{\sigma}))-$integrable.
By 
hypothesis we deduce that
there is
$F:X\to B(G)$
such that
\begin{itemize}
\item
$(\forall x\in X_{0})(F(x)=f_{x}(R))$;
\item
$
F:X\to\lr{B(G)}{\sigma(B(G),\n)}
$
is scalarly essentially 
$(\mu,B(G))-$integrable.
\end{itemize}
Thus
by 
\eqref{20171902}
\begin{equation}
\label{20011902}
\int
f_{x}(R)
\,d\mu(x)
\doteqdot
\int
F(x)
\,d\mu(x)
\end{equation}
Now for all
$\sigma\in\B(\C)$
let us define the 
map
$
F^{\sigma}:X\to B(G_{\sigma})
$
such that for all
$x\in X$
$$
F^{\sigma}(x)
\doteqdot
E(\sigma)
F(x)\up G_{\sigma}.
$$
By \eqref{15051542WEAK}
we can claim that
\begin{enumerate}
\item
$
(\forall x\in X_{0})
(F^{\sigma}(x)=f_{x}(R_{\sigma}\up G_{\sigma}))
$;
\item
the map
$
F^{\sigma}:X\to\lr{B(G_{\sigma})}
{\sigma(B(G_{\sigma}),\n_{\sigma})}
$
is scalarly essentially 
$(\mu,B(G_{\sigma}))-$integrable,
moreover
\begin{equation}
\label{19571902}
\int
F^{\sigma}(x)
\,d\mu(x)
=
\int
f_{x}(R)\,d\mu(x)
\up
G_{\sigma}.
\end{equation}
\end{enumerate}
Then
the statement
will
follow
by setting
according
\eqref{20171902}
$$
\int
f_{x}(R_{\sigma}\up G_{\sigma})
\,d\mu(x)
\doteqdot
\int
F^{\sigma}(x)
\,d\mu(x).
$$
For all
$x\in X_{0}$
\begin{alignat*}{2}
F^{\sigma}(x)
&
=
E(\sigma)
f_{x}(R)\up G_{\sigma}
\\
&
=
f_{x}(R)
E(\sigma)
\up G_{\sigma}
&
\text{
by
$[f_{x}(R),E(\sigma)]=\ze$
}
\\
&
=
f_{x}(R_{\sigma}\up G_{\sigma})
&
\text{
by
Key
Lemma
\ref{II31031834}.
}
\end{alignat*}
Hence $(1)$ of our claim
follows.
For all
$\psi\in\n,
x\in X$
\begin{alignat}{2}
\label{13580902}
\psi
\circ
\mathcal{L}(E(\sigma))
\circ
\RM(E(\sigma))
\left(
F(x)
\right)
&
\doteq
\psi
\left(
E(\sigma)
F(x)
E(\sigma)
\right)
\notag
\\
&
=
\psi_{\sigma}
\left(
E(\sigma)
F(x)
\up
G_{\sigma}
\right)
\notag
\\
&
\doteq
\psi_{\sigma}
\left(
F^{\sigma}(x)
\right).
\end{alignat}
Here in the second equality
we deduce
by 
Lemma
\ref{13061030}
that
for all
$T\in B(G)$
we have
$
\xi_{\sigma}
\left(E(\sigma)T\up G_{\sigma}
\right)
=
E(\sigma)T E(\sigma)
$.
$
F:X\to\lr{B(G)}{\sigma(B(G),\n)}
$
is
scalarly essentially
$\mu-$integrable,
and for all
$\psi\in\n$,
$\psi
\circ
\mathcal{L}(E(\sigma))
\circ
\RM(E(\sigma))
\in\n$,
hence
by
\eqref{13580902}
the map
$$
F^{\sigma}:X\to\lr{B(G_{\sigma})}
{\sigma(B(G_{\sigma}),\n_{\sigma})}
\text{
is
scalarly essentially
$\mu-$integrable.
}
$$
Now
by
\eqref{18051301}
we
have for all
$
v\in
G_{\sigma}
$
\begin{equation}
\label{13322102}
\int
f_{x}(R)\,d\mu(x)
v
=
\int
f_{x}(R)\,d\mu(x)
E(\sigma)
v
=
E(\sigma)
\int
f_{x}(R)\,d\mu(x)
v
\in
G_{\sigma},
\end{equation}
moreover
$
\int
f_{x}(R)\,d\mu(x)
\in
B(G)
$
so
\begin{equation*}
\int
f_{x}(R)\,d\mu(x)
\up
G_{\sigma}
\in
B(G_{\sigma}).
\end{equation*}
 
Therefore
for all
$\psi\in\n$
\begin{alignat*}{2}
&
\psi_{\sigma}
\left(
\int
f_{x}(R)\,d\mu(x)
\up G_{\sigma}
\right)
\\
&
=
\psi_{\sigma}
\left(
E(\sigma)
\int
f_{x}(R)\,d\mu(x)
\up G_{\sigma}
\right)
\hspace{5pt}
\text{
by 
\eqref{13322102}
}
\\
&
=
\psi
\left(
E(\sigma)
\int
f_{x}(R)\,d\mu(x)\,
E(\sigma)
\right)
\hspace{5pt}
\text{
by Lemma
\ref{13061030}
}
\\
&
\doteq
\psi
\circ
\mathcal{L}(E(\sigma))
\circ
\RM(E(\sigma))
\left(
\int
f_{x}(R)\,d\mu(x)
\right)
\\
&
\doteq
\psi
\circ
\mathcal{L}(E(\sigma))
\circ
\RM(E(\sigma))
\left(
\int
F(x)\,d\mu(x)
\right)
\hspace{5pt}
\text{
by 
\eqref{20011902}
}
\\
&
=
\int
\psi
\circ
\mathcal{L}(E(\sigma))
\circ
\RM(E(\sigma))
\left(
F(x)
\right)
\,d\mu(x)
\hspace{5pt}
\text{
by 
$
\psi
\circ
\mathcal{L}(E(\sigma))
\circ
\RM(E(\sigma))
\in
\n
$
}
\\
&
=
\int
\psi_{\sigma}
\left(
F^{\sigma}(x)
\right)
\,d\mu(x)
\hspace{5pt}
\text{
by
\eqref{13580902}.
}
\end{alignat*}
Hence
\eqref{19571902}
by 
\eqref{16021301}
and
\eqref{16021301bis}
and
the statement
follows.
\end{proof}
\begin{proposition}
\label{15051146}
For all
$\sigma\in\B(\C)$
\begin{equation}
\label{13061209A}
(\n_{st}(G))_{\sigma}
=
\n_{st}(G_{\sigma})
\text{ and }
(\n_{pd}(G))_{\sigma}
=
\n_{pd}(G_{\sigma});
\end{equation}
\end{proposition}
\begin{proof}
By the
Hahn-Banach theorem 
\begin{equation}
\label{18051655}
\left(G_{\sigma}\right)^{*}
=
\{
\phi
\up
G_{\sigma}
\mid
\phi
\in
G^{*}
\}.
\end{equation}
Then we have
\begin{alignat*}{1}
(\n_{st}(G))_{\sigma}
&
\doteq
\mathfrak{L}_{\C}
(\{
\psi_{(\phi,v)}
\circ
\xi_{\sigma}
\mid
(\phi,v)
\in
G^{*}
\times
G
\})
\\
&
=
\mathfrak{L}_{\C}(
\{
\psi_{(\phi\up G_{\sigma},w)}
\mid
(\phi,w)
\in
G^{*}
\times
G_{\sigma}
\})
\\
&
=
\mathfrak{L}_{\C}
(\{
\psi_{(\rho,w)}
\mid
(\rho,w)
\in
(G_{\sigma})^{*}
\times
G_{\sigma}
\})
\\
&
\doteq
\n_{st}(G_{\sigma}).
\end{alignat*}
Here
in the third
equality we used 
\eqref{18051655},
while
in
the
second
equality
we
considered
that
for all
$(\phi,v)\in G^{*}\times G$
and
for all
$T_{\sigma}\in B(G_{\sigma})$
\begin{alignat}{2}
\label{13061209}
\psi_{(\phi,v)}
\circ\xi_{\sigma}
(T_{\sigma})
&
=
\phi(T_{\sigma}E(\sigma)v)
&
\text{
by
\eqref{13060952}
}
\notag
\\
&
=
\left(\phi\up G_{\sigma}\right)
(T_{\sigma}E(\sigma)v)
\notag
\\
&
=
\psi_{(\phi\up G_{\sigma},E(\sigma)v)}
(T_{\sigma}).
\end{alignat}
Let $G$ be a complex Hilbert space
then
\begin{equation*}
\begin{aligned}
(\n_{pd}(G))_{\sigma}
&=
\\
\left
\{
\left(
\sum_{n=0}^{\infty}
\psi_{(u_{n}^{\dagger},w_{n})}
\right)
\circ
\xi_{\sigma}
\Bigg|
\{u_{n}\}_{n\in\N},
\,
\{w_{n}\}_{n\in\N}
\subset
G,
\sum_{n=0}^{\infty}
\|u_{n}\|_{G}^{2}
<\infty,
\sum_{n=0}^{\infty}
\|w_{n}\|_{G}^{2}
<\infty
\right
\}
&
=
\\
\left
\{
\sum_{n=0}^{\infty}
(\psi_{(u_{n}^{\dagger},w_{n})}
\circ
\xi_{\sigma})
\Bigg|
\{u_{n}\}_{n\in\N},
\,
\{w_{n}\}_{n\in\N}
\subset
G,
\sum_{n=0}^{\infty}
\|u_{n}\|_{G}^{2}
<\infty,
\sum_{n=0}^{\infty}
\|w_{n}\|_{G}^{2}
<\infty
\right
\}
&
=
\\
\left
\{
\sum_{n=0}^{\infty}
\psi_{(u_{n}^{\dagger}
\up G_{\sigma},E(\sigma)w_{n})}
\Bigg|
\{u_{n}\}_{n\in\N},
\,
\{w_{n}\}_{n\in\N}
\subset
G,
\sum_{n=0}^{\infty}
\|u_{n}\|_{G}^{2}
<\infty,
\sum_{n=0}^{\infty}
\|w_{n}\|_{G}^{2}
<\infty
\right
\}
&
=
\\
\left
\{
\sum_{n=0}^{\infty}
\psi_{(E(\sigma)^{*}u_{n})^{\dagger}
\up G_{\sigma},
E(\sigma)w_{n})}
\Bigg|
\{u_{n}\}_{n\in\N},
\,
\{w_{n}\}_{n\in\N}
\subset
G,
\sum_{n=0}^{\infty}
\|u_{n}\|_{G}^{2},
\sum_{n=0}^{\infty}
\|w_{n}\|_{G}^{2}
<\infty
\right
\}
&
=
\\
\left
\{
\sum_{n=0}^{\infty}
\psi_{(a_{n}^{\dagger},b_{n})}
\Bigg|
\{a_{n}\}_{n\in\N},
\,
\{b_{n}\}_{n\in\N}
\subset
G_{\sigma},
\sum_{n=0}^{\infty}
\|a_{n}\|_{G_{\sigma}}^{2}
<\infty,
\sum_{n=0}^{\infty}
\|b_{n}\|_{G_{\sigma}}^{2}
<\infty
\right
\}
&
=
\\
\n_{pd}(G_{\sigma}&)
\end{aligned}
\end{equation*}
Here 
for any
Hilbert space
$F$
we set
$u^{\dagger}\in F^{*}$
such that
$u^{\dagger}(v)\doteqdot\lr{u}{v}$
for all $u,v\in F$,
and
the series 
in the first
equality
is
converging
with respect to the strong operator
topology
on $B(G)^{*}$,
while
all
the others
are
converging
with respect to the strong operator
topology
on $B(G_{\sigma})^{*}$.
 
The first equality follows
by
\ref{19121201},
the third 
is by 
\eqref{13061209},
the forth
by the fact that
$E(\sigma)\up G_{\sigma}=\un_{\sigma}$
the identity operator on $G_{\sigma}$.
Now we shall
show
the fifth
equality.
Notice that
$$
\sum_{n=0}^{\infty}
\|E(\sigma)w_{n}\|_{G_{\sigma}}^{2}
\doteq
\sum_{n=0}^{\infty}
\|E(\sigma)w_{n}\|_{G}^{2}
\leq
\|E(\sigma)\|^{2}
\sum_{n=0}^{\infty}
\|w_{n}\|_{G}^{2}
<\infty.
$$
While
the fact
that
$
\dagger:
H\to H^{*}
$
is a
semilinear
isometry,
we have
for all
$
n\in\N
$
that
exists
only one
$a_{n}\in G_{\sigma}$
such that
$
a_{n}^{\dagger}
=
\left(
E(\sigma)^{*}u_{n}
\right)^{\dagger}
\up
G_{\sigma}
$
moreover
$$
\|
a_{n}
\|_{G_{\sigma}}
=
\|
\left(
E(\sigma)^{*}u_{n}
\right)^{\dagger}
\up
G_{\sigma}
\|_{G_{\sigma}^{*}}.
$$
Next
\begin{alignat*}{2}
\|
\left(
E(\sigma)^{*}u_{n}
\right)^{\dagger}
\up
G_{\sigma}
\|_{G_{\sigma}^{*}}
&
=
\sup_{
\{v\in G_{\sigma}\mid
\|v\|_{G_{\sigma}}
\leq
1\}
}
|\lr{E(\sigma)^{*}u_{n}}{v}|
\\
&
=
\sup_{
\{v\in G_{\sigma}\mid
\|v\|_{G_{\sigma}}
\leq
1\}
}
|\lr{u_{n}}{v}|
\leq
\sup_{
\{v\in G\mid
\|v\|_{G}
\leq
1\}
}
|\lr{u_{n}}{v}|
\\
&
=
\|
u_{n}^{\dagger}
\|_{G^{*}}
=
\|
u_{n}
\|_{G}.
\end{alignat*}
Hence
$
\sum_{n=0}^{\infty}
\|a_{n}\|_{G_{\sigma}}^{2}
\leq
\sum_{n=0}^{\infty}
\|u_{n}\|_{G}^{2}
<\infty
$
and
the fifth
equality
follows.
\end{proof}
\section{Extension theorem
for
integral
equalities
with respect
to the 
$\sigma(B(G),\n)-$
topology}
\label{SecEXTEQ}
In 
the present 
section
wil shall prove
the 
Extension
Theorems
for integration
with respect to the 
$\sigma(B(G),\n)-$topology,
when
$\n$ is an
$E-$appropriate set:
Theorems
\ref{18051958ta}
and
when
$\n$ is an
$E-$appropriate set
with the duality
property:
Corollary
\ref{17070917TA}.
As an application
we shall consider
the cases
of the
sigma-weak 
topology:
Corollary
\ref{12121201}
and
Corollary
\ref{17070917pd};
and weak operator topology:
Corollary
\ref{18051958},
and
Corollary
\ref{17070917}.
In this section 
it will be 
adopted all the notations
defined in Section
\ref{SecWEAKINT}.
\begin{theorem}
\label{18051958taPRE}
Let
$\n$
be
an
$E-$appropriate set
and
$\{\sigma_{n}\}_{n\in\N}$
be
an
$E-$sequence
(see Definition \ref{17331001})
and
the map
$
X
\ni 
x
\mapsto
f_{x}
\in
Bor(\sigma(R))
$
be
such that
$
\w{f}_{x}
\in
\Lf{E}{\infty}(\sigma(R))
$\,
$
\mu-l.a.e.(X)
$.
Let
$
X
\ni x\mapsto f_{x}(R)
\in 
\lr{B(G)}{\sigma(B(G),\n)}
$
be
scalarly
essentially
$(\mu,B(G))-$integrable
and
$
g,
h\in
Bor(\sigma(R))
$.
 
If for all
$n\in\N$ 
\begin{equation}
\label{19181001ta}
g(R_{\sigma_{n}}\up G_{\sigma_{n}})
\int\,
f_{x}(R_{\sigma_{n}}\up G_{\sigma_{n}})
\,d\,
\mu(x)
\subseteq
h(R_{\sigma_{n}}\up G_{\sigma_{n}})
\end{equation}
then 
\begin{equation}
\label{16241402ta}
g(R)
\int\,
f_{x}(R)
\,d\,
\mu(x)
\up
\Theta
=
h(R)
\up
\Theta.
\end{equation}
In 
\eqref{19181001ta}
the weak-integral is
with respect to the 
measure $\mu$
and with respect to the
$\sigma(B(G_{\sigma_{n}}),\n_{\sigma_{n}})-$
topology,
while
in
\eqref{16241402ta}
$$
\Theta
\doteqdot
Dom
\left(
g(R)
\int\,
f_{x}(R)
\,d\,
\mu(x)
\right)
\cap
Dom(h(R))
$$
and
the weak-integral is
with respect to the 
measure $\mu$
and with respect to the
$\sigma(B(G),\n)-$
topology.
\end{theorem}
\begin{proof}
\eqref{19181001ta}
is well set
since
Theorem
\ref{14050121}.
 
By 
\eqref{II03041401}
for all
$y\in\Theta$
\begin{alignat}{2}
\label{19051223ta}
g(R)
\int 
f_{x}(R)
\,d\,\mu(x)
\,
y
&
=
\lim_{n\in\N}
E(\sigma_{n})
g(R)
\int 
f_{x}(R)
\,d\,\mu(x)
\,
y
\notag
\\
\intertext{
by
statement $(g)$
of Theorem $18.2.11$
of \cite{ds}
and
\textbf{
\eqref{18051301}
}
}
&=
\lim_{n\in\N}
g(R)
\int 
f_{x}(R)
\,d\,\mu(x)
\,
E(\sigma_{n})
y
\notag
\\
\intertext{
by
\textbf{\eqref{restrta}}
and
Lemma
\ref{II31031834}
applied to 
$g(R)$
}
&=
\lim_{n\in\N}
g(R_{\sigma_{n}}\up G_{\sigma_{n}})
\int 
f_{x}(R_{\sigma_{n}}\up G_{\sigma_{n}})
\,d\,\mu(x)
\,
E(\sigma_{n})
y
\notag
\\
\intertext{
by 
hypothesis
\eqref{19181001ta}
}
&=
\lim_{n\in\N}
h(R_{\sigma_{n}}\up G_{\sigma_{n}})
\,
E(\sigma_{n})
y
\notag
\\
\intertext{
by 
Lemma
\ref{II31031834}
and
statement $(g)$
of Theorem $18.2.11$
of \cite{ds}
}
&=
\lim_{n\in\N}
E(\sigma_{n})
h(R)
\,
y
\notag
\\
&=
h(R)
\,
y.
\end{alignat}
In the last equality
we considered
\eqref{II03041401}.
\end{proof}
\begin{maintheorem}
[
\textbf{
$\sigma(B(G),\n)-$
Extension Theorem
}
]
\label{18051958ta}
Let
$R$
be
a 
possibly unbounded 
scalar type spectral operator
in $G$,
$E$
its resolution of the
identity,
$\n$
an
$E-$appropriate set.
Let
the map
$
X
\ni 
x
\mapsto
f_{x}
\in
Bor(\sigma(R))
$
be
such that
$
\w{f}_{x}
\in
\Lf{E}{\infty}(\sigma(R))
$\,
$
\mu-l.a.e.(X)
$
and the map
$
X
\ni x\mapsto f_{x}(R)
\in 
\lr{B(G)}{\sigma(B(G),\n)}
$
be
scalarly
essentially
$(\mu,B(G))-$integrable.
Finally
let
$
g,
h\in
Bor(\sigma(R))
$
and
$
\w{h}\in\Lf{E}{\infty}(\sigma(R))
$.
 
If
$\{\sigma_{n}\}_{n\in\N}$
is
an
$E-$sequence
and for all
$n\in\N$
\begin{equation}
\label{19211001ta}
g(R_{\sigma_{n}}\up G_{\sigma_{n}})
\int\,
f_{x}(R_{\sigma_{n}}\up G_{\sigma_{n}})
\,d\,
\mu(x)
\subseteq
h(R_{\sigma_{n}}\up G_{\sigma_{n}})
\end{equation}
then
$h(R)
\in
B(G)
$
and
\begin{equation}
\label{16011401ta}
g(R)
\int\,
f_{x}(R)
\,d\,
\mu(x)
=
h(R).
\end{equation}
In
\eqref{19211001ta}
the weak-integral is
with respect to the 
measure $\mu$
and with respect to the
$\sigma(B(G_{\sigma_{n}}),\n_{\sigma_{n}})-$
topology,
while
in
\eqref{16011401ta}
the weak-integral is
with respect to the 
measure $\mu$
and with respect to the
$\sigma(B(G),\n)-$
topology.
\end{maintheorem}
Notice that
$g(R)$ is 
a
possibly 
\textbf{
unbounded
}
operator
in
$G$.
\begin{proof}
Theorem
$18.2.11.$ of \cite{ds}
and
hypothesis
$
\w{h}\in\Lf{E}{\infty}(\sigma(R))
$
imply
that
$h(R)\in B(G)$,
so 
by 
\eqref{16241402ta}
we can deduce
\begin{equation}
\label{W15051936ta}
g(R)
\int\,
f_{x}(R)
\,d\,
\mu(x)
\subseteq
h(R).
\end{equation}
Let us set
\begin{equation}
\label{17502703ta}
(\forall n\in\N)
(\delta_{n}\doteqdot\overset{-1}{|g|}([0,n])).
\end{equation}
We claim that
\begin{equation}
\label{17542703ta}
\begin{cases}
\bigcup_{n\in\N}\delta_{n}=\sigma(R)\\
n\geq m\Rightarrow\delta_{n}\supseteq\delta_{m}\\
(\forall n\in\N)(g(\delta_{n})\text{ is bounded. })\\
\end{cases}
\end{equation}
Since
$|g|\in Bor(\sigma(R))$
we have
$\delta_{n}\in\B(\C)$
for all $n\in\N$,
so
$\{\delta_{n}\}_{n\in\N}$
is an $E-$sequence,
hence by \eqref{II03041401}
\begin{equation}
\label{10511201ta}
\lim_{n\in\N}E(\delta_{n})=\un
\end{equation}
with respect to the strong operator topology
on $B(G)$.
 
Indeed
the first equality 
follows
by
$\bigcup_{n\in\N}
\delta_{n}
\doteq
\bigcup_{n\in\N}\overset{-1}{|g|}([0,n])
=
\overset{-1}{|g|}
\left(\bigcup_{n\in\N}[0,n]\right)
=
\overset{-1}{|g|}(\R^{+})
=
Dom(g)\doteqdot\sigma(R)
$,
the second
by the fact that
$\overset{-1}{|g|}$
preserves the inclusion,
the third 
by the inclusion
$|g|(\delta_{n})\subseteq [0,n]$.
Hence our claim.
By the 
third statement 
of \eqref{17542703ta},
$\delta_{n}\in\B(\C)$
and statement
$3$ of Lemma 
\ref{II31031834}
\begin{equation}
\label{19292703ta}
(\forall n\in\N)
(E(\delta_{n})G
\subseteq
Dom(g(R))).
\end{equation}
By
\eqref{18051301}
and
\eqref{19292703ta}
for all
$n\in\N$
$$
\int\,
f_{x}(R)
\,d\,
\mu(x)
E(\delta_{n})
G
\subseteq
E(\delta_{n})
G
\subseteq
Dom(g(R)).
$$
Therefore
$$
(\forall n\in\N)
(\forall v\in G)
\left
(E(\delta_{n})v
\in 
Dom
\left(
g(R)
\int\,
f_{x}(R)
\,d\,
\mu(x)
\right)
\right).
$$
Hence
by
\eqref{10511201ta}
\begin{equation}
\label{WII04041215ta}
\D
\doteqdot
Dom
\left(
g(R)
\int\,
f_{x}(R)
\,d\,
\mu(x)
\right)
\text{ is dense in }
G.
\end{equation}
But
$
\int\,
f_{x}(R)
\,d\,
\mu(x)
\in
B(G)
$
and
$g(R)$
is
closed 
by 
Theorem $18.2.11.$
of \cite{ds},
so by
Lemma
\ref{II04041224}
we have 
\begin{equation}
\label{WII04041228ta}
g(R)
\int\,
f_{x}(R)
\,d\,
\mu(x)
\text{ is closed.}
\end{equation}
But we know that
$h(R)\in B(G)$
so by 
\eqref{W15051936ta}
we 
deduce
\begin{equation}
\label{WII04041231ta} 
g(R)
\int\,
f_{x}(R)
\,d\,
\mu(x)
\in
B(\D,G).
\end{equation}
The
\eqref{WII04041228ta},
\eqref{WII04041231ta} 
and
Lemma \ref{II04041227}
allow us to state that
$\D$ is closed in $G$,
therefore
by \eqref{WII04041215ta}
$$
\D=G.
$$
Hence by 
\eqref{W15051936ta}
we can conclude that 
the
statement 
holds.
\end{proof}
Now we shall show a corollary
of the previous theorem,
in which
we give conditions 
on the maps $f_{x}$
ensuring
that 
$f_{x}(R)\in B(G)$,
and
that
$
X\ni x\mapsto f_{x}(R)\in B(G)
$
is
scalarly essentially
$(\mu,B(G))-$integrable
with respect to the 
$
\sigma(B(G),\n)-
$
topology.
\begin{corollary}
[
\textbf{
$\sigma(B(G),\n)-$
Extension Theorem.
Duality case.
}
]
\label{17070917TA}
Let
$\n$
be
an
$E-$appropriate
set with 
the
duality
property
and
$
X\ni x\mapsto f_{x}\in Bor(\sigma(R))
$.
Assume
that there is
$X_{0}\subseteq X$
such that
$\complement X_{0}$
is 
$\mu-$locally
negligible
and
$
\w{f}_{x}
\in
\Lf{E}{\infty}(\sigma(R))
$
for all
$x\in X_{0}$,
moreover
let there exist
$F:X\to B(G)$
extending
$
X_{0}\ni x\mapsto
f_{x}(R)\in B(G)
$
such that for all
$
\omega\in\n
$
the
map
$
X
\ni 
x
\mapsto 
\omega(F(x))
\in 
\C
$
is
$\mu-$measurable
and
\begin{equation}
\label{14541002}
(X\ni x
\mapsto
\|F(x)\|_{B(G)})
\in
\F{ess}{}(X;\mu).
\end{equation}
 
If
$g,h\in Bor(\sigma(R))$
such that
$
\w{h}\in\Lf{E}{\infty}(\sigma(R))
$
and
$\{\sigma_{n}\}_{n\in\N}$
is
an
$E-$sequence
such that for all
$n\in\N$
holds
\eqref{19211001ta}
then
the statement
of
Theorem
\ref{18051958ta}
holds.
Moreover
if
$\n$
is
an
$E-$appropriate
set with 
the
isometric
duality
property
$$
\left
\|
\int
f_{x}(R)
\,d\mu(x)
\right
\|_{B(G)}
\leq
\int^{\bullet}
\|f_{x}(R)\|_{B(G)}
\,d|\mu|(x).
$$
\end{corollary}
\begin{proof}
By
the duality property
of hypothesis,
and
Theorem
\ref{19052218ta}
the map
$
X
\ni
x
\mapsto
f_{x}(R)
\in
\lr{B(G)}{\sigma(B(G),\n)}
$
is scalarly essentially
$(\mu,B(G))-$integrable.
Hence the first part of the
statement by
Theorem
\ref{18051958ta}.
The inequality
follows
by
\eqref{16170702},
\eqref{20171902}
and
\eqref{04250403}.
\end{proof}
Now we will give the 
corollaries of the previous two
results
in the cases in which
$\n=\n_{st}(G)$
or
$\n=\n_{pd}(G)$
and 
$G$
be
a Hilbert space.
\begin{corollary}
\label{18051958}
The statement of 
Theorem
\ref{18051958ta}
holds
if 
$\n$
is replaced by
$\n_{st}(G)$
and
$\n_{\sigma_{n}}$
is replaced by
$\n_{st}(G_{\sigma_{n}})$,
for all $n\in\N$.
\end{corollary}
\begin{proof}
By Remark
\ref{14551501}
we know that
$\n_{st}(G)$
is 
an
$E-$appropriate
set,
therefore
the statement 
by
\eqref{13061209A}
and
Theorem
\ref{18051958ta}.
\end{proof}
\begin{corollary}
\label{12121201}
The statement of 
Theorem
\ref{18051958ta}
holds
if $G$ is a complex Hilbert space,
$\n$
is replaced by
$\n_{pd}(G)$
and
$\n_{\sigma_{n}}$
is replaced by
$\n_{pd}(G_{\sigma_{n}})$,
for all $n\in\N$.
\end{corollary}
\begin{proof}
By Remark
\ref{14551501}
we know that
$\n_{pd}(G)$
is in particular
an
$E-$appropriate
set,
therefore
the statement 
by
\eqref{13061209A}
and
Theorem
\ref{18051958ta}.
\end{proof}
\begin{theorem}
[
\textbf{
Sigma-weak
Extension Theorem}
]
\label{17070917pd}
Let 
$G$
be
a Hilbert space,
then
the statement
of
Corollary
\ref{17070917TA}
holds
if we set
$
\n
\doteqdot
\n_{pd}(G)
$
and
$
\n_{\sigma_{n}}
\doteqdot
\n_{pd}(G_{\sigma_{n}})
$
for all
$n\in\N$.
\end{theorem}
\begin{proof}
By
Remark
\ref{14551501}
$\n_{pd}(G)$
is an 
$E-$appropriate set with
the isometric
duality property,
so
we obtain
the statement
by
Corollary
\ref{17070917TA}
and
by
\eqref{13061209A}.
\end{proof}
\begin{corollary}
[
\textbf{
Weak
Extension Theorem}
]
\label{17070917}
Let 
$G$
be
reflexive,
then
the statement
of
Corollary
\ref{17070917TA}
holds
if we set
$
\n
\doteqdot
\n_{st}(G)
$
and
$
\n_{\sigma_{n}}
\doteqdot
\n_{st}(G_{\sigma_{n}})
$
for all
$n\in\N$.
\end{corollary}
\begin{proof}
By
Theorem
\ref{19052218}
we have that
the map
$
X
\ni
x
\mapsto
f_{x}(R)
\in
\lr{B(G)}{\sigma(B(G),\n_{st}(G))}
$
is scalarly essentially
$(\mu,B(G))-$integrable.
Hence the first part of the
statement by
Corollary
\ref{18051958}.
While the inequality
follows
by
\eqref{16170702},
\eqref{20171902}
and
\eqref{04250403}.
\end{proof}
\begin{remark}
\label{19311801}
In the case 
in which
$G$ is an Hilbert space
we 
can
obtain
Corollary
\ref{17070917}
as an application
of the duality
property of the predual
$\n_{pd}(G)$.
Indeed
as we know
$
\n_{st}(G)
\subset
\n_{pd}(G)
$,
hence
by the
Hahn-Banach theorem
for all
$\Psi_{0}
\in
\n_{st}(G)^{*}$
there exists
$\Psi
\in
\n_{pd}(G)^{*}$
such that
$\Psi\up\n_{st}(G)=\Psi_{0}$,
thus
by
the 
duality
property
$\n_{pd}(G)^{*}=B(G)$
we 
obtain
$
(\forall
\Psi_{0}
\in
\n_{st}(G)^{*})
(\exists\,
B\in B(G))
(\forall\omega\in\n_{st}(G))
(\Psi_{0}(\omega)=\omega(B))
$,
which ensures
that the
weak-integral
with respect
to
the measure
$\mu$
and
with respect
to the weak operator topology
of the map
$
X\ni x
\mapsto
f_{x}(R)
\in 
B(G)
$
belongs
to
$B(G)$.
\end{remark}
\begin{remark}
\label{29021937}
Let
$\Df\subset G$
be a
linear
subspace of
$G$
and
$E:\B(\C)\to\Pr(G)$
be a countably
additive spectral measure,
then
by
\eqref{20522502}
for
all
$f\in\bb$,
$\phi\in G^{*}$
and
$v\in\Df$
that
\begin{equation}
\label{02050103}
\left|
\phi\left(
\I{\C}{E}(f)v
\right)
\right|
=
\left|
\int
f(\lambda)\,
d\,(\psi_{\phi,v}\circ E)(\lambda)
\right|
\leq
4M
\|f\|_{\sup}
\|\phi\|
\|v\|,
\end{equation}
where
$
M\doteqdot
\sup_{\delta\in\B(\C)}
\|E(\delta)\|
$,
$
\psi_{\phi,v}:
B(G)\ni A\mapsto 
\phi(Av)\in\C
$
and 
$\bb$
is the space
of 
all
totally
$\B(\C)-$measurable
complex
maps
on $\C$.
Next
we know
that
\begin{equation}
\label{02070103}
H(\C)
\subset
\bb.
\end{equation}
Here
$H(\C)$
is
the space
of all compactly supported 
complex 
continuous functions
on $\C$, 
with the direct limit topology,
of the
spaces
$H(\C;K)$
with
$K$
running
in
the class of all compact subsets of $\C$;
where
$H(\C;K)$
is 
the space
of
all
complex
continuous functions 
$f:\C\to\C$
such
that
$\textrm{supp}(f)\subset K$
with
the
topology of uniform convergence.
Let us set
$$
F_{w}^{\Df}
\doteqdot
\ov{B(\Df,G)}
\text{ in }
\Lf{w}{}(\Df,G),
$$
where
$
\Lf{w}{}(\Df,G)
$
is the Hausdorff
locally convex
space
of all linear operators
on $\Df$ at values in
$G$ with the topology
generated by the following
set of seminorms
$$
\{
\Lf{w}{}(\Df,G)
\ni 
B
\mapsto
|q_{\phi,v}(B)|
\mid
(\phi,v)\in 
G^{*}\times\Df
\},
$$
where
$
q_{\phi,v}(B)
\doteqdot
\phi(Bv)
$
for all
$(\phi,v)\in 
G^{*}\times\Df
$
and
$B\in\Lf{w}{}(\Df,G)$,
while
$B(\Df,G)$
is the space
of all bounded operators
belonging to
$\Lf{w}{}(\Df,G)$.
By 
\eqref{02070103}
we can define
$$
\mb{m}_{E}:
H(\C)
\ni
f
\mapsto
\left(
\I{\C}{E}(f)
\up\Df
\right)
\in
F_{w}^{\Df}
$$
Moreover
by
\eqref{02050103}
we have, 
with the notations
in \ref{13051406},
that
for all
compact
$K$
the operator
$
\mb{m}_{E}
\up
H(\C;K)
$
is 
continuous.
Therefore
as a corollary
of 
the
general
result
in
statement
$(ii)$
Proposition
$5$,
$n^{\circ} 4$,
\S $4$,
Ch $2$
of
\cite{BourTVS}
about
locally convex
final topologies,
so in particular for
the inductive limit topology,
we deduce that
$
\mb{m}_{E}
$
is continuous
on $H(\C)$
i.e.
$$
\mb{m}_{E}
\text{ 
is a 
vector
measure on $\C$
with vales
in 
$F_{w}^{\Df}$.
}
$$
Here,
by generalizing
to the complex case
the definition
$1$,
$n^{\circ} 1$,
\S $2$,
Ch $6$
of
\cite{IntBourb},
we call
a 
vector
measure on a locally
compact space
$X$
with values
in a complex 
Hausdorff
locally convex space
$Y$
any
$\C-$linear
continuous
map
$
\mb{m}:H(X)\to Y
$.
Furthermore for all
$(\phi,v)\in G^{*}\times\Df$
\begin{alignat*}{2}
q_{\phi,v}\circ\mb{m}_{E}
&
=
\psi_{\phi,v}
\circ 
\I{\C}{E}
\up H(\C)
\notag
\\
&
=
\I{\C}{\psi_{\phi,v}\circ E}
\up H(\C).
&
\text{ by \eqref{20522502}}
\end{alignat*}
Hence 
$$
\Lf{1}{}(\C;q_{\phi,v}\circ\mb{m}_{E})
=
\Lf{1}{}(\C;\psi_{\phi,v}\circ E),
$$
where
the left hand side 
it is intended
in the sense of Ch $4$
of \cite{IntBourb},
while
the right hand side
it is intended
in the standard sense,
and for all
$
f\in
\Lf{1}{}(\C;q_{\phi,v}\circ\mb{m}_{E})
$
\begin{equation}
\label{04260103}
\int
f(\lambda)\,d\,
(q_{\phi,v}\circ\mb{m}_{E})(\lambda)
=
\int
f(\lambda)\,d\,
(\psi_{\phi,v}\circ E)(\lambda)
\end{equation}
Finally
let
us
assume
that
$\Df$
is dense,
then
for all
$f\in Bor(\textrm{supp}\,E)$
such that
$Dom(f(E))=\Df$
by 
\eqref{II01041439}
we have
$$
f(E)\in F_{w}^{\Df},
$$
and
by
Theorem $18.2.11$ of \cite{ds}
for all
$(\phi,v)\in G^{*}\times\Df$
we have
$
f\in\Lf{1}{}(\C;\psi_{\phi,v}\circ E)
$
and
\begin{equation}
\label{04400103}
q_{\phi,v}(f(E))
=
\int
f(\lambda)\,d\,
(\psi_{\phi,v}\circ E)(\lambda).
\end{equation}
Therefore
by adopting
the definitions
in
$n^{\circ} 2$,
\S $2$,
Ch $6$
of
\cite{IntBourb},
we deduce 
by
\eqref{04260103}
that
each
$f\in Bor(\textrm{supp}\,E)$
such that
$Dom(f(E))=\Df$
is
\emph{
essentially
integrable
for
$\mb{m}_{E}$
}
and
$$
f(E)
=
\int
f(\lambda)
\,d\,
\mb{m}_{E}(\lambda).
$$
Here
$
\int
f(\lambda)
\,d\,
\mb{m}_{E}(\lambda)
$
is
the
\emph{
integral
of $f$
with respect to
$\mb{m}_{E}$}.
Thus
if
$R$
is
an
unbounded scalar type spectral operator
in $G$,
then
for all
$f\in Bor(\sigma(R))$
such that
$Dom(f(R))=\Df$
\emph{
$f$ is
essentially
integrable
for
$\mb{m}_{E}$
}
and
$$
f(R)
=
\int
f(\lambda)
\,d\,
\mb{m}_{E}(\lambda).
$$
\end{remark}
\section{
Generalization
of the
Newton-Leibnitz
formula
}
In this section we shall
apply the results
of the previous one
in order to prove
Newton-Leibnitz
formulas
for
integration with respect to the
$\sigma(B(G),\n)-$topology,
when
$\n$ is an
$E-$appropriate set
with the duality
property,
for 
integration
with respect to the
sigma-weak topology,
and
for 
integration
with respect to the
weak operator topology.
\begin{corollary}
[
\textbf{
$
\sigma(B(G),\n)-
$
Newton-Leibnitz
formula 1}
]
\label{20051321ta}
Let
$R$
be
a 
possibly unbounded 
scalar type spectral operator
in 
$G$,
$U$
an 
open
neighborhood of
$\sigma(R)$,
$S:U\to\C$
an analytic map
and
$\n$
an 
$E-$appropriate set
with the duality
property.
Assume
that
$S:U\to\C$
is
an analytic map
and there is
$L>0$
such that
$]-L,L[\,
\cdot
U\subseteq U$
and
\begin{enumerate}
\item
$
\w{S_{t}}
\in
\Lf{E}{\infty}(\sigma(R))
$,
$
\w{\left(
\frac{d\,S}{d\,\lambda}
\right)_{t}}
\in
\Lf{E}{\infty}(\sigma(R))
$
for all 
$t\in ]-L,L[$;
\item
$$
\int^{*}
\left\|
\w{
{\left(
\frac{d\,S}{d\,\lambda}
\right)_{t}}
}
\right\|_{\infty}^{E}
\,
dt
<
\infty
$$
(here
the upper
integral
is with respect
to 
the
Lebesgue measure on $]-L,L[$);
\item for all
$\omega\in\n$
the map
$
]-L,L[
\ni t
\mapsto
\omega\left(
\frac{d\,S}{d\,\lambda}(tR)
\right)
\in
\C
$
is Lebesgue measurable.
\end{enumerate}
Then
for all
$u_{1},u_{2}\in ]-L,L[$
$$
R
\int_{u_{1}}^{u_{2}}
\frac{d\,S}{d\,\lambda}(t R)
\,d\,t
=
S(u_{2} R)
-
S(u_{1} R)
\in B(G).
$$
Here
the
integral
is
the 
weak-integral
of the 
map
$
[u_{1},u_{2}]
\ni
t
\mapsto
\frac{d\,S}{d\,\lambda}(t R)
\in
B(G)
$
with
respect to the 
Lebesgue measure on
$[u_{1},u_{2}]$
and
with respect to the 
$
\sigma(B(G),\n)-
$
topology.
Moreover 
if
$\n$
is
an 
$E-$appropriate set
with the 
isometric
duality
property
and
$
M
\doteqdot
\sup_{\sigma\in\B(\C)}
\|E(\sigma)\|_{B(G)}
$
then
\begin{equation}
\label{17001501}
\left
\|
\int_{u_{1}}^{u_{2}}
\frac{d\,S}{d\,\lambda}
(tR)
\,d t
\right
\|_{B(G)}
\leq
4M
\int_{[u_{1},u_{2}]}^{*}
\left\|
\w{
\left(
\frac{d\,S}{d\,\lambda}
\right)_{t}
}
\right\|_{\infty}^{E}
\,
dt.
\end{equation}
\end{corollary}
\begin{proof}
Let 
$\mu$
the
Lebesgue measure on
$[u_{1},u_{2}]$.
By
\eqref{29051755},
the
hypotheses,
and
statement $(c)$
of
Theorem $18.2.11$
of \cite{ds}
we
have
\begin{description}
\item
[a]
$
(\forall t\in [u_{1},u_{2}])
(S
(t R)
\in B(G)
)
$;
\item
[b]
$
(\forall t\in [u_{1},u_{2}])
(
\frac{d\,S}{d\,\lambda}
(t R)
\in B(G)
)
$;
\item
[c]
$
([u_{1},u_{2}]\ni t
\mapsto
\|
\frac{d\,S}{d\,\lambda}
(t R)
\|_{B(G)})
\in
\F{1}{}([u_{1},u_{2}];\mu)
$,
\end{description}
So by
hypothesis 
$(3)$,
by 
$(c)$
and
Theorem \ref{19052218ta}
we have that
the map
\begin{equation}
\label{16581501}
[u_{1},u_{2}]
\ni
t
\mapsto
\frac{d\,S}{d\,\lambda}
(tR)
\in
\lr{B(G)}{\sigma(B(G),\n)}
\end{equation}
is scalarly essentially
$(\mu,B(G))-$
integrable
and
if
$\n$
is
an 
$E-$appropriate set
with the 
isometric
duality
property
then
its
weak-integral
satisfies
\eqref{17001501}.
This means that, 
made exception 
for
\eqref{19211001ta},
all
the hypotheses of 
Theorem 
\ref{18051958ta}
hold
for
$
X
\doteqdot 
[u_{1},u_{2}]
$,
$
h
\doteqdot 
(S_{u_{2}}
-
S_{u_{1}})
\up
\sigma(R)
$,
$
g:
\sigma(R)
\ni
\lambda
\mapsto
\lambda
\in
\C
$
and finally
for the 
map
$
[u_{1},u_{2}]
\ni
t
\mapsto
f_{t}
\doteqdot
\left(
\frac{d\,S}{d\,\lambda}
\right)_{t}
\up
\sigma(R)
$.
 
Next
let
$
\sigma\in\B(\C)
$
be
bounded,
so
by
Key
Lemma 
\ref{II31031834}
$R_{\sigma}\up G_{\sigma}$
is a scalar
type spectral operator
such that
$
R_{\sigma}\up G_{\sigma}
\in B(G_{\sigma})$,
moreover
by
\eqref{25052019}
$U$
is
an
open
neighborhood of
$\sigma(R_{\sigma}\up G_{\sigma})$.
Thus
we can 
apply
statement $(3)$
of Theorem
\ref{19500603}
to the
Banach space
$G_{\sigma}$,
the analytic
map
$S$
and
to the operator
$R_{\sigma}\up G_{\sigma}$.
In particular
the map
$
[u_{1},u_{2}]
\ni
t
\mapsto
\frac{d\,S}{d\,\lambda}
(t (R_{\sigma}\up G_{\sigma}))
\in
B(G_{\sigma})
$
is
Lebesgue integrable
in
$\|\cdot\|_{B(G_{\sigma})}-$topology,
that
is
in the meaning
of
Definition $2$,
$n^{\circ} 4$,
\S $3$,
Ch. $IV$
of
\cite{IntBourb}.
By
Lemma
\ref{13061030}
\,
$
\xi_{\sigma}\in 
B(B(G_{\sigma}),B(G))
$,
so
$$
\n_{\sigma}
\subset
B(G_{\sigma})^{*}.
$$
Therefore
we deduce by using
Theorem
$1$, $IV.35$
of the
\cite{IntBourb},
that for all
$\omega_{\sigma}\in\n_{\sigma}$
the map
$
[u_{1},u_{2}]
\ni
t
\mapsto
\omega_{\sigma}
\left(
\frac{d\,S}{d\,\lambda}
(t (R_{\sigma}\up G_{\sigma}))
\right)
\in
\C
)
$
is
Lebesgue integrable,
in addition for all
$\omega_{\sigma}\in\n_{\sigma}$
$$
\int_{u_{1}}^{u_{2}}
\omega_{\sigma}
\left(
\frac{d\,S}{d\,\lambda}
(t (R_{\sigma}\up G_{\sigma}))
\right)
\,d\,t
=
\omega_{\sigma}
\left(
\oint_{u_{1}}^{u_{2}}
\frac{d\,S}{d\,\lambda}
(t (R_{\sigma}\up G_{\sigma}))
\,d\,t
\right).
$$
Thus
we 
can state that
$
[u_{1},u_{2}]
\ni
t
\mapsto
\frac{d\,S}{d\,\lambda}
(t (R_{\sigma}\up G_{\sigma}))
\in
\lr{B(G_{\sigma})}
{\sigma(B(G_{\sigma}),\n_{\sigma})}
$
is scalarly essentially
$(\mu,B(G_{\sigma}))-$integrable
and
its
weak-integral
is such that
\begin{equation}
\label{W15052324ta}
\int_{u_{1}}^{u_{2}}
\frac{d\,S}{d\,\lambda}
(t (R_{\sigma}\up G_{\sigma}))
\,d\,t
=
\oint_{u_{1}}^{u_{2}}
\frac{d\,S}{d\,\lambda}
(t (R_{\sigma}\up G_{\sigma}))
\,d\,t.
\end{equation}
Moreover
by
statement $(3)$
of Theorem
\ref{19500603}
\begin{equation*}
(R_{\sigma}\up G_{\sigma})
\oint_{u_{1}}^{u_{2}}
\frac{d\,S}{d\,\lambda}
(t (R_{\sigma}\up G_{\sigma}))
\,d\,t
=
S(u_{2} (R_{\sigma}\up G_{\sigma}))
-
S(u_{1} (R_{\sigma}\up G_{\sigma})).
\end{equation*}
Thus
by
\eqref{W15052324ta}
\begin{equation}
\label{W15052333ta}
(R_{\sigma}\up G_{\sigma})
\int_{u_{1}}^{u_{2}}
\frac{d\,S}{d\,\lambda}
(t (R_{\sigma}\up G_{\sigma}))
\,d\,t
=
S(u_{2} (R_{\sigma}\up G_{\sigma}))
-
S(u_{1} (R_{\sigma}\up G_{\sigma})).
\end{equation}
This 
implies
exactly 
the hypothesis
\eqref{19211001ta}
of 
Theorem 
\ref{18051958ta},
by choosing
for example
$
\sigma_{n}
\doteqdot 
B_{n}(\ze)
$,
for all
$n\in\N$.
Therefore
by 
Theorem 
\ref{18051958ta}
we obtain
the
statement.
\end{proof}
A particular case of the previous result
has been obtained in a different context 
in \cite[Cor. $3.3$]{sil}.
\begin{corollary}
[
\textbf{
$
\sigma(B(G),\n)-
$
Newton-Leibnitz
formula 2}
]
\label{20051321taLOC}
Let
$R$
be
a 
possibly unbounded 
scalar type spectral operator
in 
$G$,
$U$
an 
open
neighborhood of
$\sigma(R)$,
$S:U\to\C$
an analytic map
and
$\n$
an 
$E-$appropriate set
with the duality
property.
Assume
that there exists
$L>0$
such that
$
]-L,L[\,
\cdot
U\subseteq U
$
and
for all
$t\in 
]-L,L[$,
$\w{S_{t}}
\in
\Lf{E}{\infty}(\sigma(R))$
and
there exists
$K_{0}\subset ]-L,L[$
such that
$
\complement
K_{0}
$
is a
Lebesgue negligible
set
and
for all
$t\in K_{0}$,
$\w{
\left(
\frac{d\,S}{d\,\lambda}
\right)_{t}
}\in\Lf{E}{\infty}(\sigma(R))$
moreover
\begin{enumerate}
\item
there is
$F:]-L,L[
\to B(G)$
extending
$K_{0}
\ni
t
\mapsto
\frac{d\,S}{d\,\lambda}(tR)
\in B(G)
$
such that
$$
\int^{*}
\left\|
F(t)
\right\|_{B(G)}
\,
dt
<
\infty
$$
(here
the upper
integral
is with respect
to 
the
Lebesgue measure on $]-L,L[$),
\item
for all
$\omega\in\n$
the map
$]-L,L[
\ni t
\mapsto
\omega(F(t))
\in\C
$
is Lebesgue
measurable.
\end{enumerate}
Then for all
$u_{1},u_{2}\in ]-L,L[$
$$
R
\int_{u_{1}}^{u_{2}}
\frac{d\,S}{d\,\lambda}(t R)
\,d\,t
=
S(u_{2} R)
-
S(u_{1} R)
\in B(G).
$$
Here
the
integral
is
the 
weak-integral
of the 
map
$
[u_{1},u_{2}]
\ni
t
\mapsto
\frac{d\,S}{d\,\lambda}(t R)
\in
B(G)
$
with
respect to the 
Lebesgue measure on
$[u_{1},u_{2}]$
and
with respect to the 
$
\sigma(B(G),\n)-
$
topology.
Moreover
if
$\n$
is
an 
$E-$appropriate set
with the 
isometric
duality
property
then
$$
\left\|
\int_{u_{1}}^{u_{2}}
\frac{d\,S}{d\,\lambda}(t R)
\,d\,t
\right\|_{B(G)}
\leq
\int_{[u_{1},u_{2}]}^{*}
\left\|
\frac{d\,S}{d\,\lambda}(t R)
\right\|_{B(G)}
\,d\,t.
$$
\end{corollary}
\begin{proof}
By Theorem
\ref{19052218ta}
and
\eqref{20171902}
$$
[u_{1},u_{2}]
\ni
t
\mapsto
\frac{d\,S}{d\,\lambda}(t R)
\in
\lr{B(G)}{\sigma(B(G),\n)}
$$
is scalarly essentially 
$(\mu,B(G))-$integrable
and 
if
$\n$
is
an 
$E-$appropriate set
with the 
isometric
duality
property
its weak integral 
satisfies
by
\eqref{04250403}
the inequality
in the statement.
Thus
the proof
goes on as
that in
Corollary
\ref{20051321ta}.
\end{proof}
\begin{corollary}
[Sigma-Weak
Newton-Leibnitz
formula
]
\label{20051321pd}
The statement of Corollary
\ref{20051321ta}
(respectively
Corollary
\ref{20051321taLOC})
holds
if 
$G$ is a complex Hilbert space
and everywhere
$\n$
is replaced
by
$
\n_{pd}(G)
$.
\end{corollary}
\begin{proof}
By 
Remark
\ref{14551501},
$\n_{pd}(G)$
is
an
$E-$appropriate
set
with
the 
isometric
duality property,
hence
the statement
by 
Corollary
\ref{20051321ta}
(respectively
Corollary
\ref{20051321taLOC}).
\end{proof}
\begin{corollary}
[Weak
Newton-Leibnitz
formula
]
\label{19371401}
The statement of 
Corollary
\ref{20051321ta}
(respectively
Corollary
\ref{20051321taLOC})
holds
if 
$G$ is a
reflexive 
complex 
Banach space
and everywhere
$\n$
is replaced
by
$
\n_{st}(G)
$.
\end{corollary}
\begin{proof}
By 
using 
Corollary
\ref{19052220}
instead
of
Theorem \ref{19052218ta},
we obtain 
\eqref{16581501}
and
\eqref{17001501}
by
replacing
$\n$
with
$\n_{st}(G)$.
Then the proof
procedes similarly
to that
of Corollary
\ref{20051321ta}
(respectively
Corollary
\ref{20051321taLOC}).
\end{proof}

\section*{Acknowledgments}
I am very grateful to 
Professor
Victor Burenkov
for having
supervised
this work
and to
Professors
Massimo Lanza de Cristoforis
and
Karl Michael Schmidt
for 
their
helpful
comments. 
Finally I thank Professor
W. Zelazko for his interest
in this work.

\bibliography{references}
\bibliographystyle{amsplain}

\end{document}